\documentclass[11pt]{article}

\usepackage[utf8]{inputenc}
\usepackage[T1]{fontenc}
\usepackage{amsmath, amssymb}
\usepackage{graphicx}
\usepackage{hyperref}
\usepackage{geometry,xcolor}

\usepackage{amsthm}

\usepackage{amsmath, amssymb, amsthm}

\theoremstyle{plain}
\newtheorem{theorem}{Theorem}[section]

\theoremstyle{definition}
\newtheorem{definition}{Definition}[section]

\theoremstyle{remark}
\newtheorem{remark}{Remark}[section]

\numberwithin{equation}{section}
\begin{document}

\title{Green's Function and Solution Representation for a Boundary Value Problem Involving the Prabhakar Fractional Derivative}

\author{
Erkinjon Karimov\\
Department of Mathematics: Analysis, Logics and Discrete Mathematics, Ghent University,\\
Krijgslaan 297, S8, Ghent 9000, Belgium\\
\texttt{erkinjon.karimov@ugent.be}
\and
Doniyor Usmonov\\
Department of Mathematical Analysis and Differential Equations, Fergana State University,\\
Murabbiylar Str. 19, 150100 Fergana, Uzbekistan\\
\texttt{dusmonov909@gmail.com}
\and
Maftuna Mirzaeva\\
Department of Mathematical Analysis and Differential Equations, Fergana State University,\\
Murabbiylar Str. 19, 150100 Fergana, Uzbekistan\\
\texttt{maftunamirzayeva2009@gmail.com}
}

\date{\today}

\maketitle

\begin{abstract}
We investigate a first boundary value problem for a second-order partial differential equation involving the Prabhakar fractional derivative in time. Using structural properties of the Prabhakar kernel and generalized Mittag–Leffler functions, we reduce the problem to a Volterra-type integral equation. This reduction enables the explicit construction of the corresponding Green’s function. Based on the obtained Green’s function, we derive a closed-form integral representation of the solution and prove its existence and uniqueness. The results extend classical Green-function techniques to a wider class of fractional operators and provide analytical tools for further study of boundary and inverse problems associated with Prabhakar-type fractional differential equations.
\end{abstract}

\vspace{0.3cm}

\section{Introduction}
Fractional differential equations extend classical ones to non-integer orders, enabling models with memory and hereditary effects \cite{Podlubny}. They are widely used in physics (anomalous diffusion, viscoelasticity) \cite{Suzuki}, engineering (control, signal processing), biology (population dynamics), medicine (pharmacokinetics), and finance \cite{Kilbas}. By capturing long-range dependencies, they often yield more accurate models than integer-order equations \cite{El-Sayed}.
Green’s functions provide integral representations for boundary value problems. While well established for classical equations, recent work has adapted them to fractional cases. For instance, Green’s functions for time-fractional diffusion-wave equations have been constructed using potential methods \cite{Pskhu}, and explicit forms for fractional PDEs have been derived via Laplace and Fourier transforms \cite{Sarah}, \cite{Odibat}.

Most of the existing literature is devoted to fractional models involving the Riemann-Liouville or Caputo derivatives, whose kernels exhibit a power-law behavior, corresponding to a single-scale memory effect. However, many complex systems -- such as heterogeneous materials, anomalous transport processes, and relaxation phenomena in viscoelastic media -- display multi-scale or non-exponential memory patterns, which cannot be adequately described by classical fractional operators \cite{wang}.

In this context, the Prabhakar fractional derivative has emerged as a natural and powerful generalization. This operator is defined in terms of the three-parameter Mittag–Leffler function and extends the classical Riemann-Liouville and Caputo derivatives as particular cases \cite{Giusti}, \cite{Garra}. The additional parameters allow one to model a wide spectrum of memory kernels, including stretched exponential and multi-term relaxation laws, thereby providing significantly greater flexibility in applications. As a result, the Prabhakar derivative has been successfully applied in the study of transport phenomena, renewal processes, and materials with complex internal structure \cite{guisti-col}.

Furthermore, fractional differential equations involving the Prabhakar derivative admit solutions expressed in terms of generalized Mittag–Leffler functions, which possess richer asymptotic behavior than classical power-law functions and better describe experimentally observed dynamics \cite{T22}. Recent studies have addressed analytical properties, solvability, and applications of such equations, highlighting the increasing importance of this operator in modern fractional calculus.

Despite these advances, the theory of boundary value problems for partial differential equations involving the Prabhakar fractional derivative remains relatively underdeveloped, particularly with regard to the explicit construction of Green’s functions. Such constructions are essential for obtaining integral representations of solutions and for further qualitative and numerical analysis. Recent studies address unique solvability of boundary-value problems and diffusion-wave equations involving Prabhakar derivatives, highlighting its growing role in fractional calculus \cite{Al-Refai}, \cite{Magar}, \cite{Karimov}, \cite{Turdiev 1}.

Motivated by the above considerations, in the present paper we study an initial-boundary value problem for a fractional diffusion equation involving the Prabhakar fractional derivative with respect to the time variable. Our main objective is to construct the corresponding Green’s function in explicit form and to derive a representation of the solution in terms of this function. This result provides a useful analytical tool for the investigation of more general fractional diffusion models with complex memory effects.


\section{Main result} \setcounter{equation}{0}\label{sec:2}
\subsection{Formulation of a problem and main statement.}
Let $D$ be a rectangular domain: $D=\left\{ \left( t,x \right):\,\,0<t<T,\,\,0<x<a \right\},$ $0<a,T<\infty .$ 
We formulate the first initial-boundary value problem for the following sub-diffusion equation:
\begin{equation}\label{eq2.1}
    {}^{PRL}D_{0t}^{\alpha ,\,\beta ,\,\gamma ,\,\delta }u\left( t,x \right)-{{u}_{xx}}\left( t,x \right)=f\left( t,x \right).
\end{equation}
Here $f\left( t,x \right)$ is a given function, $\alpha,\beta,\gamma,\delta$ are given real numbers such that $\alpha>0,$ $0<\beta \le 1$, ${ }^{PRL} D_{0 t}^{\alpha, \beta, \gamma, \delta} y(t)={\frac{d}{d t}}{}^{P}I_{0 t}^{\alpha, 1-\beta,-\gamma, \delta} y(t),\, t>0$ is the Prabhakar fractional derivative of Riemann-Liouville type,
${ }^P I_{0 t}^{\alpha, \beta, \gamma, \delta} y(t)=\int\limits_0^t(t-s)^{\beta-1} E_{\alpha, \beta}^\gamma\left[\mathcal{\delta}(t-s)^\alpha\right] y(s)ds$ is the Prabhakar fractional integral operator, whereas
$
E_{\alpha, \beta}^\gamma[z]=\sum\limits_{k=0}^{+\infty} \frac{(\gamma)_k z^k}{\Gamma(\alpha k+\beta) k!}
$
represents a generalized Mittag-Leffler (Prabhakar) function \cite{Prabhakar}. 

\begin{definition}
    A regular solution of the equation \eqref{eq2.1} in the domain $D$ is called a function $u\left( t,x \right)$ with the regularity ${}^{P}I_{0t}^{\alpha ,\,1-\beta ,\,-\gamma ,\,\delta }u\left( t,x \right)\in C\left( \overline{D} \right)$ $\left[t^{1-\beta}u(t,x)\in C( \overline{D})\right],$ ${{u}_{xx}}\left( t,x \right)$, ${}^{PRL}D_{0t}^{\alpha ,\,\beta ,\,\gamma ,\,\delta }u\left( t,x \right)\in C\left( D \right)$ that satisfies the equation \eqref{eq2.1} at all points $\left( t,x \right)\in D.$
\end{definition}
\begin{remark}
	The assumptions $t^{1-\beta}u(t,x)\in C( \overline{D})$ and ${}^{PRL}D_{0t}^{\alpha ,\,\beta ,\,\gamma ,\,\delta }u\left( t,x \right)\in C\left( D \right)$ constitute regularity conditions with respect to the time variable in the fractional sense. They replace classical time differentiability, which is not compatible with Riemann-Liouville-type operators.
	\end{remark}

	\textbf{Problem.} Find a regular solution $u\left( t,x \right)$ of the equation \eqref{eq2.1} in the domain $D$, satisfying the following boundary and initial conditions:
\begin{equation}\label{eq2.2}
    u\left( t,0 \right)={{\varphi }_{0}}\left( t \right),\,\,\,   u\left( t,a \right)={{\varphi }_{1}}\left( t \right), \,\,\,   0<t<T,	
\end{equation}
\begin{equation}\label{eq2.3}
    \underset{t\to 0}{\mathop{\lim }}\,{}^{P}I_{0t}^{\alpha ,\,1-\beta ,\,-\gamma ,\,\delta }u\left( t,x \right)=\tau \left( x \right),\,\,\,  0\le x\le a,			 
\end{equation}
where ${{\varphi }_{0}}\left( t \right),$ ${{\varphi }_{1}}\left( t \right),$ $\tau \left( x \right)$ are  assumed to possess sufficient smoothness so that the conditions of the following theorem are satisfied.

\begin{theorem} 
 Let\, ${{t}^{1-\beta }}{{\varphi }_{0}}\left( t \right),$ ${{t}^{1-\beta }}{{\varphi }_{1}}\left( t \right)\in C\left[ 0;T \right],$ $\tau \left( x \right)\in C\left[ 0;a \right],$ ${{t}^{1-\beta }}f\left( t,x \right)\in C\left( \overline{D} \right),$ $f\left( t,x \right)$ satisfies the H{\"o}lder condition with respect to $x$ and the conditions
 $$\underset{t\to 0}{\mathop{\lim }}\,{}^{P}I_{0t}^{\alpha ,\,1-\beta ,\,-\gamma ,\,\delta }{{\varphi }_{0}}\left( t \right)=\tau \left( 0 \right),\,\,\, \underset{t\to 0}{\mathop{\lim }}\,{}^{P}I_{0t}^{\alpha ,\,1-\beta ,\,-\gamma ,\,\delta }{{\varphi }_{1}}\left( t \right)=\tau \left( a \right).$$
Then there exists a unique regular solution of the equation \eqref{eq2.1} in the domain $D,$ satisfying the conditions \eqref{eq2.2} and \eqref{eq2.3}, represented as
$$u\left( t,x \right)=\int\limits_{0}^{t}{{{\varphi }_{0}}\left( \eta  \right){{G}_{s}}\left( t,x,\eta ,0 \right)d\eta }-\int\limits_{0}^{t}{{{\varphi }_{1}}\left( \eta  \right){{G}_{s}}\left( t,x,\eta ,a \right)d\eta }+$$
\begin{equation}\label{eq2.4}
    +\,\int\limits_{0}^{a}{\tau \left( s \right)G\left( t,x,0,s \right)ds}+\int\limits_{0}^{t}{\int\limits_{0}^{a}{f\left( \eta ,s \right)G\left( t,x,\eta ,s \right)dsd\eta }}.
\end{equation}
\end{theorem}
Here 
$$G\left( t,x,\eta ,s \right)=\frac{{{\left( t-\eta  \right)}^{{{\beta }_{1}}-1}}}{2}\times$$
$$\times\sum\limits_{n=-\infty }^{\infty }{\left[ {{E}_{12}}\left( \left. \begin{matrix}
   -{{\gamma }_{1}},1,{{\gamma }_{1}};\,\,\,\,\,\,\,\,\,\,\,\,\,\,\,\,\,\,\,\,\,\,\,\,\,\,\,  \\
   -{{\beta }_{1}},\alpha ,{{\beta }_{1}};-{{\gamma }_{1}},{{\gamma }_{1}};1,1;1,1  \\
\end{matrix} \right|\begin{matrix}
   -\left| x-s+2an \right|{{\left( t-\eta  \right)}^{-{{\beta }_{1}}}}  \\
   \delta {{\left( t-\eta  \right)}^{\alpha }}  \\
\end{matrix} \right) \right.}-$$
\begin{equation}\label{eq2.5}
    \left. -{{E}_{12}}\left( \left. \begin{matrix}
   -{{\gamma }_{1}},1,{{\gamma }_{1}};\,\,\,\,\,\,\,\,\,\,\,\,\,\,\,\,\,\,\,\,\,\,\,\,\,\,\,  \\
   -{{\beta }_{1}},\alpha ,{{\beta }_{1}};-{{\gamma }_{1}},{{\gamma }_{1}};1,1;1,1  \\
\end{matrix} \right|\begin{matrix}
   -\left| x+s+2an \right|{{\left( t-\eta  \right)}^{-{{\beta }_{1}}}}  \\
   \delta {{\left( t-\eta  \right)}^{\alpha }}  \\
\end{matrix} \right) \right],
\end{equation}
$${{E}_{12}}\left( \left. \begin{matrix}
   {{\alpha }_{1}},{{\widetilde\beta }_{1}},{{\delta }_{1}};\,\,\,\,\,\,\,\,\,\,\,\,\,\,\,\,\,\,\,\,\,\,\,\,\,\,\,\,\,\,\,\,\,\,\,\,\,\,\,\,\,  \\
   {{\alpha }_{2}},{{\widetilde\beta }_{2}},{{\delta }_{2}};{{\alpha }_{3}},{{\delta }_{3}};{{\alpha }_{4}},{{\delta }_{4}};{{\widetilde\beta }_{3}},{{\delta }_{5}}  \\
\end{matrix} \right|\begin{matrix}
   x  \\
   y  \\
\end{matrix} \right)=$$
$$=\sum\limits_{n=0}^{+\infty }{\sum\limits_{m=0}^{+\infty }{\frac{\Gamma \left( {{\alpha }_{1}}n+{{\widetilde\beta }_{1}}m+{{\delta }_{1}} \right){{x}^{n}}{{y}^{m}}}{\Gamma \left( {{\alpha }_{2}}n+{{\widetilde\beta }_{2}}m+{{\delta }_{2}} \right)\Gamma \left( {{\alpha }_{3}}n+{{\delta }_{3}} \right)\Gamma \left( {{\alpha }_{4}}n+{{\delta }_{4}} \right)\Gamma \left( {{\widetilde\beta }_{3}}m+{{\delta }_{5}} \right)}}},$$
$\left( x,y,{{\alpha }_{l}},{{\widetilde\beta }_{i}},{{\delta }_{j}}\in \mathbb{R};\min \left\{ {{\alpha }_{l}},{{\widetilde\beta }_{i}} \right\}>0;\left( l=\left\{ 1,...,4 \right\},i=\left\{ 1,2,3 \right\}, j=\left\{ 1,...,5 \right\} \right) \right),$
in which the double series converges for $x,y\in \mathbb{R}$, if ${{\Delta }_{1}}>0$, and ${{\Delta }_{2}}>0$, whereas ${{\Delta }_{1}}={{\alpha }_{2}}+{{\alpha }_{3}}+{{\alpha }_{4}}-{{\alpha }_{1}}$, ${{\Delta }_{2}}={{\widetilde\beta }_{2}}+{{\widetilde\beta }_{3}}-{{\widetilde\beta }_{1}}$, ${{\beta }_{1}}=\frac{\beta }{2},\,\,\,{{\gamma }_{1}}=\frac{\gamma }{2}$ \cite{Turdiev 2}.

\begin{remark}
If $\delta=0$ or $\gamma=0$ in the problem considered, the result of \cite[pp. 107--112]{Pskhu 2} will follow.
\end{remark}

\begin{remark}
	The set of assumptions of Theorem 2.1 is nonempty. Indeed, the functions
	$$\varphi_0(t)=t^{\beta-1},\,\varphi_1(t)=t^{\beta-1},\,\tau(x)=\sin x+\Gamma(\beta),\,f(t,x)=t^{\beta-1} \sin x$$
	with $a=\pi$ so that $x\in[0,\pi],\,t\in[0,T]$, satisfy all the imposed conditions.
	\end{remark}

    Solution \eqref{eq2.4} satisfies the conditions of Theorem 2.1 and Definition 2.1 \textcolor{red}{(see for details Section \ref{secA1})}.

\subsection{Proof of Theorem 2.1.}

\subsubsection{Reduction to the system of first-order differential equations and corresponding initial-boundary problems.}

Let $u\left( t,x \right)$ satisfy the equation \eqref{eq2.1} and the conditions \eqref{eq2.2} and \eqref{eq2.3}. Since ${{t}^{1-\beta }}u\left( t,x \right)\in C\left( \overline{D} \right)$, it follows that \textcolor{red}{(see for details Section \ref{secA2})}
\begin{equation}\label{eq2.6}
\underset{t\to 0}{\mathop{\lim }}\,{}^PI_{0t}^{\alpha ,\,1-\frac{\beta }{2},\,-\frac{\gamma }{2},\,\delta }u\left( t,x \right)=0.
\end{equation}

\begin{remark}
	Unlike the Caputo-Prabhakar case, for the Riemann-Liouville-type Prabhakar derivative continuity of $u(t,x)$ does not imply vanishing of the corresponding Prabhakar integral at $t=0$, which necessitates additional weighted regularity assumptions.
	\end{remark}
Since ${}^{PRL}D_{0t}^{\alpha ,\,\beta ,\,\gamma ,\,\delta }u\left( t,x \right)={}^{PRL}D_{0t}^{\alpha ,\,\frac{\beta }{2},\,\frac{\gamma }{2},\,\delta }\bigg({}^{PRL}D_{0t}^{\alpha ,\,\frac{\beta }{2},\,\frac{\gamma }{2},\,\delta }u\left( t,x \right)\bigg)$ \textcolor{red}{(see for details Section \ref{secA3})}, we can rewrite the equation \eqref{eq2.1} as follows:
$$\left( {}^{PRL}D_{0t}^{\alpha ,\,\frac{\beta }{2},\,\frac{\gamma }{2},\,\delta }-\frac{\partial }{\partial x} \right)\left( {}^{PRL}D_{0t}^{\alpha ,\,\frac{\beta }{2},\,\frac{\gamma }{2},\,\delta }+\frac{\partial }{\partial x} \right)u\left( t,x \right)=f\left( t,x \right).$$
By denoting 
\begin{equation}\label{eq2.7}
    v\left( t,x \right)={}^{PRL}D_{0t}^{\alpha ,\,{{\beta }_{1}},\,{{\gamma }_{1}},\,\delta }u\left( t,x \right)+{{u}_{x}}\left( t,x \right),
\end{equation}
we obtain that the function $u\left( t,x \right)$ is a solution of the following system:
\begin{equation}\label{eq2.8}
  \left\{ \begin{matrix}
   {}^{PRL}D_{0t}^{\alpha ,\,{{\beta }_{1}},\,{{\gamma }_{1}},\,\delta }u\left( t,x \right)+{{u}_{x}}\left( t,x \right)=v\left( t,x \right),  \\
   {}^{PRL}D_{0t}^{\alpha ,\,{{\beta }_{1}},\,{{\gamma }_{1}},\delta }v\left( t,x \right)-{{v}_{x}}\left( t,x \right)=f\left( t,x \right).  \\
\end{matrix} \right.	 
\end{equation}
From \eqref{eq2.3}, \eqref{eq2.6}, and \eqref{eq2.7}, we can easily arrive at the following condition \cite{Usmonov}:
$$\underset{t\to 0}{\mathop{\lim }}\,{}^PI_{0t}^{\alpha ,\,1-{{\beta }_{1}},\,-{{\gamma }_{1}},\,\delta }v\left( t,x \right)=\tau \left( x \right),\,\,\,\,0\le x\le a.$$
Firstly, we see the first equation of the system \eqref{eq2.8}. According to condition \eqref{eq2.6}, we obtain the following problem:
\begin{equation}\label{eq2.9}
    {}^{PRL}D_{0t}^{\alpha ,\,{{\beta }_{1}},\,{{\gamma }_{1}},\,\delta }u\left( t,x \right)+{{u}_{x}}\left( t,x \right)=v\left( t,x \right),		
\end{equation}
\begin{equation}\label{eq2.10}
    \underset{t\to 0}{\mathop{\lim }}\,{}^P I_{0t}^{\alpha ,\,1-{{\beta }_{1}},\,-{{\gamma }_{1}},\,\delta }u\left( t,x \right)=0,\,\,\,\,   0\le x\le a,
\end{equation}
\begin{equation}\label{eq2.11}
    u\left( t,0 \right)={{\varphi }_{0}}\left( t \right),   0<t<T.		
\end{equation}

The solution of the problem \eqref{eq2.9}-\eqref{eq2.10}-\eqref{eq2.11} will be represented as follows \textcolor{red}{(see for details Section \ref{secA4})}:
\begin{equation}\label{eq2.12}
u\left( t,x \right)=\int\limits_{0}^{t}{{{\varphi }_{0}}\left( \eta  \right)\omega \left( t-\eta ,x \right)d\eta +}\int\limits_{0}^{t}{\int\limits_{0}^{x}{v\left( \eta ,\xi  \right)\omega \left( t-\eta ,x-\xi  \right)d\xi d\eta }},
\end{equation}
where
$$\omega \left( t,x \right)=\sum\limits_{n=0}^{+\infty }{\frac{{{\left( -1 \right)}^{n}}{{x}^{n}}}{n!}{{t}^{-{{\beta }_{1}}n-1}}E_{\alpha ,-{{\beta }_{1}}n}^{-{{\gamma }_{1}}n}\left[ \delta {{t}^{\alpha }} \right]}=$$
\begin{equation}\label{eq2.13}
    ={{t}^{-1}}{{E}_{12}}\left( \left. \begin{matrix}
   -{{\gamma }_{1}},1,0;\,\,\,\,\,\,\,\,\,\,\,\,\,\,\,\,\,\,\,\,\,\,\,\,\,\,\,  \\
   -{{\beta }_{1}},\alpha ,0;-{{\gamma }_{1}},0;1,1;1,1  \\
\end{matrix} \right|\begin{matrix}
   -x{{t}^{-{{\beta }_{1}}}}  \\
   \delta {{t}^{\alpha }}  \\
\end{matrix} \right).
\end{equation}
The function $\omega(t,x)$ does not have a singularity at $t=0$ \textcolor{red}{(see for details Section \ref{secA5})}.

Similarly, we will get the initial-boundary problem for the second equation of system \eqref{eq2.8}:
\begin{equation}\label{eq2.14}
    {}^{PRL}D_{0t}^{\alpha ,\,{{\beta }_{1}},\,{{\gamma }_{1}},\delta }v\left( t,x \right)-{{v}_{x}}\left( t,x \right)=f\left( t,x \right),
\end{equation}
\begin{equation}\label{eq2.15}
    \underset{t\to 0}{\mathop{\lim }}\,{}^P I_{0t}^{\alpha ,\,1-{{\beta }_{1}},\,-{{\gamma }_{1}},\,\delta }v\left( t,x \right)=\tau \left( x \right),  \,\,\,  0\le x\le a,
\end{equation}
\begin{equation}\label{eq2.16}
    v\left( t,a \right)=\psi \left( t \right),\,\,\,   0<t<T,
\end{equation}
$\psi \left( t \right)$ is an unknown function.

The solution to the problem \eqref{eq2.14}-\eqref{eq2.15}-\eqref{eq2.16} has the following form \textcolor{red}{(see for details Section \ref{secA6})}: 
$$v\left( t,x \right)=\int\limits_{0}^{t}{\psi \left( \eta  \right)\omega \left( t-\eta ,a-x \right)d\eta }+$$
\begin{equation}\label{eq2.17}
+\int\limits_{x}^{a}{\tau \left( \xi  \right)\omega \left( t,\xi -x \right)d\xi }+\int\limits_{0}^{t}{\int\limits_{x}^{a}{f\left( \eta ,\xi  \right)\omega \left( t-\eta ,\xi -x \right)d\xi d\eta }},
\end{equation}
where $\omega \left( t,x \right)$ is defined by \eqref{eq2.13}.

\subsubsection{Reduction to the integral equation.}

We substitute \eqref{eq2.17} into \eqref{eq2.12}:
$$u\left( t,x \right)=\int\limits_{0}^{t}{{{\varphi }_{0}}\left( \eta  \right)\omega \left( t-\eta ,x \right)d\eta +}\int\limits_{0}^{t}{\int\limits_{0}^{x}{\left[ \int\limits_{0}^{\eta }{\psi \left( y \right)\omega \left( \eta -y,a-\xi  \right)dy+} \right.}}$$
$$\left. +\int\limits_{\xi }^{a}{\tau \left( s \right)\omega \left( \eta ,s-\xi  \right)ds+}\int\limits_{0}^{\eta }{\int\limits_{\xi }^{a}{f\left( y,s \right)\omega \left( \eta -y,s-\xi  \right)dsdy}} \right]\times$$
    $$\times\omega \left( t-\eta ,x-\xi  \right)d\xi d\eta .$$

 \begin{remark}
	In the sequel, interchanges between integration and infinite summation are justified by the absolute and uniform convergence of the series representations of the involved special functions under the imposed assumptions. We therefore omit repetitive technical details.
\end{remark}

Let us now change the order of integration in the multiple integrals:

1)	$\int\limits_{0}^{t}{d\eta \int\limits_{0}^{x}{d\xi \int\limits_{0}^{\eta }{\psi \left( y \right)\omega \left( \eta -y,a-\xi  \right)\omega \left( t-\eta ,x-\xi  \right)dy=}}}$

$=\int\limits_{0}^{t}{\psi \left( y \right)dy\int\limits_{0}^{x}{d\xi \int\limits_{y}^{t}{\omega \left( \eta -y,a-\xi  \right)\omega \left( t-\eta ,x-\xi  \right)d\eta ;}}}$

2)  $\int\limits_{0}^{t}{d\eta \int\limits_{0}^{x}{d\xi \int\limits_{\xi }^{a}{\tau \left( s \right)\omega \left( \eta ,s-\xi  \right)\omega \left( t-\eta ,x-\xi  \right)ds=}}}$

 $=\left( \int\limits_{0}^{x}{ds\int\limits_{0}^{s}{d\xi +\int\limits_{x}^{a}{ds\int\limits_{0}^{x}{d\xi }}}} \right)\tau \left( s \right)\int\limits_{0}^{t}{\omega \left( \eta ,s-\xi  \right)\omega \left( t-\eta ,x-\xi  \right)d\eta ;}$
 
3)  $\int\limits_{0}^{t}{d\eta }\int\limits_{0}^{x}{d\xi }\int\limits_{\xi }^{a}{ds\int\limits_{0}^{\eta }{f\left( y,s \right)\omega \left( \eta -y,s-\xi  \right)\omega \left( t-\eta ,x-\xi  \right)dy=}}$

   $=\left( \int\limits_{0}^{x}{ds\int\limits_{0}^{s}{d\xi +\int\limits_{x}^{a}{ds\int\limits_{0}^{x}{d\xi }}}} \right)\int\limits_{0}^{t}{f\left( y,s \right)dy}\int\limits_{y}^{t}{\omega \left( \eta -y,s-\xi  \right)\omega \left( t-\eta ,x-\xi  \right)d\eta .}$

As a result, we obtain following:
$$u\left( t,x \right)=\int\limits_{0}^{t}{{{\varphi }_{0}}\left( \eta  \right)\omega \left( t-\eta ,x \right)d\eta +}$$
$$+\int\limits_{0}^{t}{\psi \left( y \right)dy\int\limits_{0}^{x}{d\xi \int\limits_{y}^{t}{\omega \left( \eta -y,a-\xi  \right)\omega \left( t-\eta ,x-\xi  \right)d\eta +}}}$$
$$+\left( \int\limits_{0}^{x}{ds\int\limits_{0}^{s}{d\xi +\int\limits_{x}^{a}{ds\int\limits_{0}^{x}{d\xi }}}} \right)\tau \left( s \right)\int\limits_{0}^{t}{\omega \left( \eta ,s-\xi  \right)\omega \left( t-\eta ,x-\xi  \right)d\eta +}$$
$$+\left( \int\limits_{0}^{x}{ds\int\limits_{0}^{s}{d\xi +\int\limits_{x}^{a}{ds\int\limits_{0}^{x}{d\xi }}}} \right)\int\limits_{0}^{t}{f\left( y,s \right)dy}\times$$
\begin{equation}\label{eq2.18}
\times\int\limits_{y}^{t}{\omega \left( \eta -y,s-\xi  \right)\omega \left( t-\eta ,x-\xi  \right)d\eta .}
\end{equation}
The following formula holds for the function $\omega \left( t,x \right)$ \textcolor{red}{(see for details Section \ref{secA7})}:
\begin{equation}\label{eq2.19}
\int\limits_{0}^{y}{\omega \left( y-t,{{x}_{1}} \right)}\,\omega \left( t,{{x}_{2}} \right)dt=\omega \left( y,{{x}_{1}}+{{x}_{2}} \right).
\end{equation}

Now we apply the formula \eqref{eq2.19} for the second integral kernel on the right-hand side of the equality \eqref{eq2.18}:
$$\int\limits_{0}^{x} d\xi \int\limits_{y}^{t}\omega \left( \eta -y,a-\xi  \right)\omega \left( t-\eta ,x-\xi  \right)d\eta =\int\limits_{0}^{x}{\omega \left( t-y,x+a-2\xi  \right)d\xi }.$$
In the resulting integral, we perform the substitution $x+a-2\xi =s$, then we use the following equality \textcolor{red}{(see for details Section \ref{secA8})}:
$$\int\limits_{a-x}^{a+x}{\omega \left( t-y,s \right)ds}={}^P I_{yt}^{\alpha ,\,{{\beta }_{1}},\,{{\gamma }_{1}},\,\delta }\left[ \omega \left( t-y,a-x \right)-\omega \left( t-y,a+x \right) \right]$$
and deduce
$$\int\limits_{0}^{x}{\omega \left( t-y,x+a-2\xi  \right)d\xi }=\frac{1}{2}{}^P I_{yt}^{\alpha ,\,{{\beta }_{1}},\,{{\gamma }_{1}},\,\delta }\left[ \omega \left( t-y,a-x \right)-\omega \left( t-y,a+x \right) \right].$$
Similarly, by transforming the kernels of the remaining terms from \eqref{eq2.18}, we obtain an expression for the function $u\left( t,x \right)$:
$$u\left( t,x \right)=\int\limits_{0}^{t}{{{\varphi }_{0}}\left( \eta  \right)\omega \left( t-\eta ,x \right)d\eta +}\int\limits_{0}^{t}{\psi \left( y \right)W\left( t-y,a-x,a+x \right)dy+}$$
$$+\int\limits_{0}^{a}{\tau \left( s \right)W\left( t,\left| s-x \right|,s+x \right)ds+}$$
\begin{equation}\label{eq2.20}
    +\int\limits_{0}^{a}{\int\limits_{0}^{t}{f\left( y,s \right)W\left( t-y,\left| s-x \right|,s+x \right)dyds,}}
\end{equation}
here $W\left( t,{{x}_{1}},{{x}_{2}} \right)=\frac{1}{2}{}^P I_{0t}^{\alpha ,\,{{\beta }_{1}},\,{{\gamma }_{1}},\,\delta }\left[ \omega \left( {t},{{x}_{1}} \right)-\omega \left( {t},{{x}_{2}} \right) \right].$

In the obtained equation, we let $x$ tend to $a$. Hence, we have for the second term of the right-hand side of \eqref{eq2.20}
$$\underset{x\to a}{\mathop{\lim }}\,\int\limits_{0}^{t}{\psi \left( y \right)W\left( t-y,a-x,a+x \right)dy=}\frac{1}{2}\underset{x\to a}{\mathop{\lim }}\,\int\limits_{0}^{t}{\psi \left( y \right)dy\times }$$
$$\times \int\limits_{y}^{t}{{{\left( t-s \right)}^{{{\beta }_{1}}-1}}}E_{\alpha ,\,{{\beta }_{1}}}^{{{\gamma }_{1}}}\left[ \delta {{\left( t-s \right)}^{\alpha }} \right]\left[ \omega \left( s-y,a-x \right)-\omega \left( s-y,a+x \right) \right]ds.$$
We perform the substitution $s-y=z$ and change the order of integration:
$$\frac{1}{2}\underset{x\to a}{\mathop{\lim }}\,\int\limits_{0}^{t}{\psi \left( y \right)dy\times }$$
$$\times \int\limits_{0}^{t-y}{{{\left( \left( t-y \right)-z \right)}^{{{\beta }_{1}}-1}}}E_{\alpha ,\,{{\beta }_{1}}}^{{{\gamma }_{1}}}\left[ \delta {{\left( \left( t-y \right)-z \right)}^{\alpha }} \right]\left[ \omega \left( z,a-x \right)-\omega \left( z,a+x \right) \right]dz=$$
$$=\frac{1}{2}\underset{x\to a}{\mathop{\lim }}\,\int\limits_{0}^{t}{\left[ \omega \left( z,a-x \right)-\omega \left( z,a+x \right) \right]{}^P I_{0\left( t-z \right)}^{\alpha ,\,{{\beta }_{1}},\,{{\gamma }_{1}},\,\delta }\psi \left( t-z \right)dz}.$$
Now we make the substitution $t-z=y$ and write the result as follows:
$$\frac{1}{2}\underset{x\to a}{\mathop{\lim }}\,\int\limits_{0}^{t}{\left[ \omega \left( t-y,a-x \right)-\omega \left( t-y,a+x \right) \right]{}^P I_{0y}^{\alpha ,\,{{\beta }_{1}},\,{{\gamma }_{1}},\,\delta }\psi \left( y \right)dy}=$$
$$=\frac{1}{2}\underset{x\to a}{\mathop{\lim }}\,\int\limits_{0}^{t}{\omega \left( t-y,a-x \right){}^P I_{0y}^{\alpha ,\,{{\beta }_{1}},\,{{\gamma }_{1}},\,\delta }\psi \left( y \right)dy}-
$$
$$
-\frac{1}{2}\int\limits_{0}^{t}{\omega \left( t-y,2a \right){}^P I_{0y}^{\alpha ,\,{{\beta }_{1}},\,{{\gamma }_{1}},\,\delta }\psi \left( y \right)dy}.$$
The following holds:
$$\frac{1}{2}\underset{x\to a}{\mathop{\lim }}\,\int\limits_{0}^{t}{\omega \left( t-y,a-x \right){}^P I_{0y}^{\alpha ,\,{{\beta }_{1}},\,{{\gamma }_{1}},\,\delta }\psi \left( y \right)dy}=\frac{1}{2}{}^P I_{0t}^{\alpha ,\,{{\beta }_{1}},\,{{\gamma }_{1}},\,\delta }\psi \left( t \right).$$
Hence, the second integral of \eqref{eq2.20} can be written in the following form:
$$\underset{x\to a}{\mathop{\lim }}\,\int\limits_{0}^{t}{\psi \left( y \right)W\left( t-y,a-x,a+x \right)dy=}$$
\begin{equation}\label{eq2.21}
    ={{\psi }_{1}}\left( t \right)-\int\limits_{0}^{t}{{{\psi }_{1}}\left( t \right)\omega \left( t-y,2a \right)dy},
\end{equation}
where 
\begin{equation}\label{eq2.22}
{{\psi }_{1}}\left( t \right)=\frac{1}{2}{}^P I_{0t}^{\alpha ,\,{{\beta }_{1}},\,{{\gamma }_{1}},\,\delta }\psi \left( t \right).
\end{equation}
Since, according to the condition, $u\left( t,a \right)={{\varphi }_{1}}\left( t \right),$ it follows from \eqref{eq2.20} and \eqref{eq2.21} that
\begin{equation}\label{eq2.23}
{{\psi }_{1}}\left( t \right)-\int\limits_{0}^{t}{\omega }\left( t-y,2a \right){{\psi }_{1}}\left( y \right)dy=F\left( t \right),
\end{equation}
where 
$$F\left( t \right)={{\varphi }_{1}}\left( t \right)-\int\limits_{0}^{t}{{{\varphi }_{0}}\left( \eta  \right)\omega \left( t-\eta ,a \right)d\eta -}$$
$$-\int\limits_{0}^{a}{\tau \left( s \right)W\left( t,a-x,a+x \right)ds-}\int\limits_{0}^{a}{\int\limits_{0}^{t}{f\left( y,s \right)W\left( t-y,a-x,a+x \right)dyds.}}$$
Thus, the function ${{\psi }_{1}}\left( t \right)$ represents a solution of the Volterra equation \eqref{eq2.23}. This solution can be write in the form of the Neumann series \cite{Kolmogorov}:
$${{\psi }_{1}}\left( t \right)=F\left( t \right)+\sum\limits_{n=1}^{\infty }{{{K}^{n}}F\left( t \right)},$$
where $KF\left( t \right)=\int\limits_{0}^{t}{F\left( y \right)\omega \left( t-y,2a \right)}\,dy.$ Using \eqref{eq2.19}, it is easy to obtain that \textcolor{red}{(see for details Section \ref{secA9})}
$${{K}^{n}}F\left( t \right)=\int\limits_{0}^{t}{F\left( y \right)\omega \left( t-y,2na \right)}\,dy.$$
Therefore, the solution to \eqref{eq2.23} has the following form: 
\begin{equation}\label{eq2.24}
    {{\psi }_{1}}\left( t \right)=F\left( t \right)+\int\limits_{0}^{t}{F\left( y \right)\sum\limits_{n=1}^{\infty }{\omega \left( t-y,2na \right)}}\,dy.
\end{equation}
According to \eqref{eq2.22} and \eqref{eq2.24}, we get
\begin{equation}\label{eq2.25}
\psi \left( t \right)=2{}^{PRL}D_{0t}^{\alpha ,\,{{\beta }_{1}},\,{{\gamma }_{1}},\,\delta }\left[ F\left( t \right)+\int\limits_{0}^{t}{F\left( y \right)\sum\limits_{n=1}^{\infty }{\omega \left( t-y,2na \right)}}\,dy \right].
\end{equation}
Taking \eqref{eq2.25} and the form of the function $F\left( t \right)$ into account, \eqref{eq2.20} can be reformulated as follows:
$$u\left( t,x \right)=\int\limits_{0}^{t}{{{\varphi }_{0}}\left( \eta  \right)\omega \left( t-\eta ,x \right)d\eta +}$$
$$+\int\limits_{0}^{t}{2{}^{PRL}D_{0y}^{\alpha ,\,{{\beta }_{1}},\,{{\gamma }_{1}},\,\delta }\left\{ \left[ {{\varphi }_{1}}\left( y \right)-\int\limits_{0}^{y}{{{\varphi }_{0}}\left( \eta  \right)\omega \left( y-\eta ,a \right)d\eta -} \right. \right.}$$
$$\left. -\int\limits_{0}^{a}{\tau \left( s \right)W\left( y,a-s,a+s \right)ds-}\int\limits_{0}^{a}{\int\limits_{0}^{y}{f\left( z,s \right)W\left( y-z,a-s,a+s \right)dzds}} \right]+$$
$$+\int\limits_{0}^{y}{\left[ {{\varphi }_{1}}\left( \eta  \right)-\int\limits_{0}^{\eta }{{{\varphi }_{0}}\left( z \right)\omega \left( \eta -z,a \right)dz-} \right.}\int\limits_{0}^{a}{\tau \left( s \right)W\left( \eta ,a-s,a+s \right)ds-}$$
$$\left. \left. -\int\limits_{0}^{a}{\int\limits_{0}^{\eta }{f\left( y,s \right)W\left( \eta -y,a-s,a+s \right)dyds}} \right]\sum\limits_{n=1}^{\infty }{\omega \left( y-\eta ,2na \right)}\,d\eta  \right\}\times $$
$$\times W\left( t-y,a-x,a+x \right)dy+\int\limits_{0}^{a}{\tau \left( s \right)W\left( t,\left| s-x \right|,s+x \right)ds+}$$
\begin{equation}\label{eq2.26}
+\int\limits_{0}^{a}{\int\limits_{0}^{t}{f\left( y,s \right)W\left( t-y,\left| s-x \right|,s+x \right)dyds.}}
\end{equation}
Now, in \eqref{eq2.26}, we will perform some transformations. We begin by calculating the  function $W\left( t-y,a-x,a+x \right)$:
$$W\left( t-y,a-x,a+x \right)=\frac{1}{2}{}^P I_{yt}^{\alpha ,\,{{\beta }_{1}},\,{{\gamma }_{1}},\,\delta }\left[ \omega \left( t-y,a-x \right)-\omega \left( t-y,a+x \right) \right]=$$
$$=\frac{1}{2}\int\limits_{y}^{t}{{{\left( t-s \right)}^{{{\beta }_{1}}-1}}E_{\alpha ,{{\beta }_{1}}}^{{{\gamma }_{1}}}\left[ \delta {{\left( t-s \right)}^{\alpha }} \right]}\left[ \omega \left( s-y,a-x \right)-\omega \left( s-y,a+x \right) \right]ds=$$
$$=\frac{1}{2}\sum\limits_{n=0}^{+\infty }{\frac{{{\left( -1 \right)}^{n}}\left[ {{\left( a-x \right)}^{n}}-{{\left( a+x \right)}^{n}} \right]}{n!}\times }$$
$$\times \int\limits_{y}^{t}{{{\left( t-s \right)}^{{{\beta }_{1}}-1}}{{\left( s-y \right)}^{{{\beta }_{1}}n-1}}E_{\alpha ,{{\beta }_{1}}}^{{{\gamma }_{1}}}\left[ \delta {{\left( t-s \right)}^{\alpha }} \right]}E_{\alpha ,-{{\beta }_{1}}n}^{-{{\gamma }_{1}}n}\left[ \delta {{\left( s-y \right)}^{\alpha }} \right]ds.$$
We utilize \cite{Knopp}
\begin{equation}\label{eq2.27}
\sum\limits_{k=0}^{+\infty }{{{a}_{k}}}\sum\limits_{m=0}^{+\infty }{{{b}_{m}}=\sum\limits_{k=0}^{+\infty }{\sum\limits_{m=0}^{k}{{{a}_{m}}}{{b}_{k-m}}}}
\end{equation}
for $E_{\alpha ,{{\beta }_{1}}}^{{{\gamma }_{1}}}\left[ \delta {{\left( t-s \right)}^{\alpha }} \right]$ and $E_{\alpha ,-{{\beta }_{1}}n}^{-{{\gamma }_{1}}n}\left[ \delta {{\left( s-y \right)}^{\alpha }} \right]$, then we make the substitution $s=\left( t-y \right)z+y$ and apply \cite{Prudnikov}
\begin{equation}\label{eq2.28}
\sum\limits_{m=0}^{k}{\frac{{{\left( \delta  \right)}_{m}}{{\left( \gamma  \right)}_{k-m}}}{m!\left( k-m \right)!}=\frac{{{\left( \delta +\gamma  \right)}_{k}}}{k!}}
\end{equation}
to get
$$W\left( t-y,a-x,a+x \right)=\frac{1}{2}\sum\limits_{n=0}^{+\infty }{\frac{{{\left( -1 \right)}^{n}}\left[ {{\left( a-x \right)}^{n}}-{{\left( a+x \right)}^{n}} \right]}{n!}\times}$$
$${\times\sum\limits_{k=0}^{+\infty }{\sum\limits_{m=0}^{k}{\frac{{{\left( {{\gamma }_{1}} \right)}_{m}}{{\left( -{{\gamma }_{1}}n \right)}_{k-m}}{{\delta }^{k}}}{m!\left( k-m \right)!\Gamma \left( \alpha m+{{\beta }_{1}} \right)\Gamma \left( \alpha k-\alpha m-{{\beta }_{1}}n \right)}}}\times }$$
$$\times \int\limits_{y}^{t}{{{\left( t-s \right)}^{\alpha m+{{\beta }_{1}}-1}}{{\left( s-y \right)}^{\alpha k-\alpha m+{{\beta }_{1}}n-1}}}ds=$$
$$=\frac{1}{2}\sum\limits_{n=0}^{+\infty }{\frac{{{\left( -1 \right)}^{n}}\left[ {{\left( a-x \right)}^{n}}-{{\left( a+x \right)}^{n}} \right]}{n!}\times}$$
$${\times\sum\limits_{k=0}^{+\infty }{\frac{{{\delta }^{k}}{{\left( t-y \right)}^{\alpha k+{{\beta }_{1}}-{{\beta }_{1}}n-1}}}{\Gamma \left( \alpha k+{{\beta }_{1}}-{{\beta }_{1}}n \right)}\sum\limits_{m=0}^{k}{\frac{{{\left( {{\gamma }_{1}} \right)}_{m}}{{\left( -{{\gamma }_{1}}n \right)}_{k-m}}}{m!\left( k-m \right)!}}}=}$$
$$=\frac{1}{2}\sum\limits_{n=0}^{+\infty }{\frac{{{\left( -1 \right)}^{n}}\left[ {{\left( a-x \right)}^{n}}-{{\left( a+x \right)}^{n}} \right]}{n!}{{\left( t-y \right)}^{{{\beta }_{1}}-{{\beta }_{1}}n-1}}\times}$$
$${\times\sum\limits_{k=0}^{+\infty }{\frac{{{\left( {{\gamma }_{1}}-{{\gamma }_{1}}n \right)}_{k}}{{\delta }^{k}}{{\left( t-y \right)}^{\alpha k}}}{k!\Gamma \left( \alpha k+{{\beta }_{1}}-{{\beta }_{1}}n \right)}}=}$$
$$=\frac{1}{2}\sum\limits_{n=0}^{+\infty }{\frac{{{\left( -1 \right)}^{n}}\left[ {{\left( a-x \right)}^{n}}-{{\left( a+x \right)}^{n}} \right]}{n!}{{\left( t-y \right)}^{{{\beta }_{1}}-{{\beta }_{1}}n-1}}E_{\alpha ,{{\beta }_{1}}-{{\beta }_{1}}n}^{{{\gamma }_{1}}-{{\gamma }_{1}}n}\left[ \delta {{\left( t-y \right)}^{\alpha }} \right].}$$
Next, in \eqref{eq2.26}, we denote the terms that involve the function ${{\varphi }_{0}}\left( \eta  \right)$ by ${{F}_{1}}$ and proceed with their evaluation:
$${{F}_{1}}=\int\limits_{0}^{t}{{{\varphi }_{0}}\left( \eta  \right)\omega \left( t-\eta ,x \right)d\eta -}$$
$$-\int\limits_{0}^{t}{2{}^{PRL} D_{0y}^{\alpha ,\,{{\beta }_{1}},\,{{\gamma }_{1}},\,\delta }\left[ \int\limits_{0}^{y}{{{\varphi }_{0}}\left( \eta  \right)\omega \left( y-\eta ,a \right)d\eta } \right]}\,W\left( t-y,a-x,a+x \right)dy-$$
$$-\int\limits_{0}^{t}{2{}^{PRL}D_{0y}^{\alpha ,\,{{\beta }_{1}},\,{{\gamma }_{1}},\,\delta }\left\{ \int\limits_{0}^{y}{\left[ \int\limits_{0}^{\eta }{{{\varphi }_{0}}\left( z \right)\omega \left( \eta -z,a \right)dz} \right]}\sum\limits_{n=1}^{\infty }{\omega \left( y-\eta ,2na \right)}\,d\eta  \right\}}\times $$
$$\times W\left( t-y,a-x,a+x \right)dy.$$
Let us denote the integral containing $\omega \left( y-\eta ,a \right)$ as ${{F}_{1,1}}$ and the one containing $\omega \left( y-\eta ,2na \right)$ as ${{F}_{1,2}}$. Let us start by evaluating ${{F}_{1,1}}$. Firstly, we apply the operator $D_{0y}^{\alpha ,\,{{\beta }_{1}},\,{{\gamma }_{1}},\,\delta }$ and interchange the order of integration:
$${{F}_{1,1}}=\int\limits_{0}^{t}{2{}^{PRL}D_{0y}^{\alpha ,\,{{\beta }_{1}},\,{{\gamma }_{1}},\,\delta }\left[ \int\limits_{0}^{y}{{{\varphi }_{0}}\left( \eta  \right)\omega \left( y-\eta ,a \right)d\eta } \right]}\,W\left( t-y,a-x,a+x \right)dy=$$
$$=2\int\limits_{0}^{t}{\frac{\partial }{\partial y}\left\{ \int\limits_{0}^{y}{{{\left( y-s \right)}^{-{{\beta }_{1}}}}E_{\alpha ,1-{{\beta }_{1}}}^{-{{\gamma }_{1}}}\left[ \delta {{\left( y-s \right)}^{\alpha }} \right]}\,ds\times  \right.}$$
$$\times \left. \int\limits_{0}^{s}{{{\varphi }_{0}}\left( \eta  \right)\omega \left( s-\eta ,a \right)d\eta } \right\}W\left( t-y,a-x,a+x \right)dy=$$
$$=2\int\limits_{0}^{t}{\frac{\partial }{\partial y}\left\{ \int\limits_{0}^{y}{{{\varphi }_{0}}\left( \eta  \right)d\eta \int\limits_{\eta }^{y}{{{\left( y-s \right)}^{-{{\beta }_{1}}}}}E_{\alpha ,1-{{\beta }_{1}}}^{-{{\gamma }_{1}}}\left[ \delta {{\left( y-s \right)}^{\alpha }} \right]}\,\omega \left( s-\eta ,a \right)ds \right\}}\times $$
$$\times W\left( t-y,a-x,a+x \right)dy.$$
Now, for convenience, we calculate only the part of the integral with respect to $s$ that involves $\omega \left( s-\eta ,a \right)$:
$${\int\limits_{\eta }^{y}{{{\left( y-s \right)}^{-{{\beta }_{1}}}}}E_{\alpha ,1-{{\beta }_{1}}}^{-{{\gamma }_{1}}}\left[ \delta {{\left( y-s \right)}^{\alpha }}\right]\omega \left( s-\eta ,a \right)ds=\sum\limits_{n=0}^{+\infty }{\frac{{{\left( -1 \right)}^{n}}{{a}^{n}}}{n!}}\times}$$
$${\times\int\limits_{\eta }^{y}{{{\left( y-s \right)}^{-{{\beta }_{1}}}}{{\left( s-\eta  \right)}^{-{{\beta }_{1}}n-1}}}E_{\alpha ,1-{{\beta }_{1}}}^{-{{\gamma }_{1}}}\left[ \delta {{\left( y-s \right)}^{\alpha }} \right]E_{\alpha ,-{{\beta }_{1}}n}^{-{{\gamma }_{1}}n}\left[ \delta {{\left( s-\eta  \right)}^{\alpha }} \right]ds=}$$
$$=\sum\limits_{n=0}^{+\infty }{\frac{{{\left( -1 \right)}^{n}}{{a}^{n}}}{n!}\sum\limits_{k=0}^{+\infty }{\sum\limits_{m=0}^{k}{\frac{{{\left( -{{\gamma }_{1}} \right)}_{m}}{{\left( -{{\gamma }_{1}}n \right)}_{k-m}}{{\delta }^{k}}}{m!\left( k-m \right)!\Gamma \left( \alpha m+1-{{\beta }_{1}} \right)\Gamma \left( \alpha k-\alpha m-{{\beta }_{1}}n \right)}}}}\times $$
$$\times \int\limits_{\eta }^{y}{{{\left( y-s \right)}^{\alpha m-{{\beta }_{1}}}}{{\left( s-\eta  \right)}^{\alpha k-\alpha m-{{\beta }_{1}}n-1}}}ds=\left\{ s=\left( y-\eta  \right)z+\eta  \right\}=$$
$$=\sum\limits_{n=0}^{+\infty }{\frac{{{\left( -1 \right)}^{n}}{{a}^{n}}}{n!}{{\left( y-\eta  \right)}^{-{{\beta }_{1}}-{{\beta }_{1}}n}}E_{\alpha ,1-{{\beta }_{1}}-{{\beta }_{1}}n}^{-{{\gamma }_{1}}-{{\gamma }_{1}}n}\left[ \delta {{\left( y-\eta  \right)}^{\alpha }} \right]}.$$
We substitute the obtained result into ${{F}_{1,1}}$  and integrate by parts:
$${{F}_{1,1}}=2\int\limits_{0}^{t}{\frac{\partial }{\partial y}\left\{ \int\limits_{0}^{y}{{{\varphi }_{0}}\left( \eta  \right)\sum\limits_{n=0}^{+\infty }{\frac{{{\left( -1 \right)}^{n}}{{a}^{n}}}{n!}{{\left( y-\eta  \right)}^{-{{\beta }_{1}}-{{\beta }_{1}}n}}}\times}\right.}$$
$$\left.{{\times E_{\alpha ,1-{{\beta }_{1}}-{{\beta }_{1}}n}^{-{{\gamma }_{1}}-{{\gamma }_{1}}n}\left[ \delta {{\left( y-\eta  \right)}^{\alpha }} \right]}\,d\eta } \right\}W\left( t-y,a-x,a+x \right)dy=$$
$$=\left[ W\left( t-y,a-x,a+x \right)\times  \right.$$
$$\left. \left. \times \int\limits_{0}^{y}{{{\varphi }_{0}}\left( \eta  \right)\sum\limits_{n=0}^{+\infty }{\frac{{{\left( -1 \right)}^{n}}{{a}^{n}}}{n!}{{\left( y-\eta  \right)}^{-{{\beta }_{1}}-{{\beta }_{1}}n}}E_{\alpha ,1-{{\beta }_{1}}-{{\beta }_{1}}n}^{-{{\gamma }_{1}}-{{\gamma }_{1}}n}\left[ \delta {{\left( y-\eta  \right)}^{\alpha }} \right]}\,d\eta } \right] \right|_{y=0}^{y=t}-$$
$$-2\int\limits_{0}^{t}{{{{{W}'}}_{y}}\left( t-y,a-x,a+x \right)dy\times }$$
$$\times \int\limits_{0}^{y}{{{\varphi }_{0}}\left( \eta  \right)\sum\limits_{n=0}^{+\infty }{\frac{{{\left( -1 \right)}^{n}}{{a}^{n}}}{n!}{{\left( y-\eta  \right)}^{-{{\beta }_{1}}-{{\beta }_{1}}n}}E_{\alpha ,1-{{\beta }_{1}}-{{\beta }_{1}}n}^{-{{\gamma }_{1}}-{{\gamma }_{1}}n}\left[ \delta {{\left( y-\eta  \right)}^{\alpha }} \right]}\,d\eta }.$$
The following holds here \textcolor{red}{(see for details Section \ref{secA10})}:
$$\underset{y\to t}{\mathop{\lim }}\,\left[ W\left( t-y,a-x,a+x \right)\times  \right.$$
$$\left. \times \int\limits_{0}^{y}{{{\varphi }_{0}}\left( \eta  \right)\sum\limits_{n=0}^{+\infty }{\frac{{{\left( -1 \right)}^{n}}{{a}^{n}}}{n!}{{\left( y-\eta  \right)}^{-{{\beta }_{1}}-{{\beta }_{1}}n}}E_{\alpha ,1-{{\beta }_{1}}-{{\beta }_{1}}n}^{-{{\gamma }_{1}}-{{\gamma }_{1}}n}\left[ \delta {{\left( y-\eta  \right)}^{\alpha }} \right]}\,d\eta } \right]=0$$
and 
$$\underset{y\to 0}{\mathop{\lim }}\,\left[ W\left( t-y,a-x,a+x \right)\times  \right.$$
$$\left. \times \int\limits_{0}^{y}{{{\varphi }_{0}}\left( \eta  \right)\sum\limits_{n=0}^{+\infty }{\frac{{{\left( -1 \right)}^{n}}{{a}^{n}}}{n!}{{\left( y-\eta  \right)}^{-{{\beta }_{1}}-{{\beta }_{1}}n}}E_{\alpha ,1-{{\beta }_{1}}-{{\beta }_{1}}n}^{-{{\gamma }_{1}}-{{\gamma }_{1}}n}\left[ \delta {{\left( y-\eta  \right)}^{\alpha }} \right]}\,d\eta } \right]=0.$$
We continue by computing ${{F}_{1,1}}$, that is, we change the order of integration in the multiple integral:
$${{F}_{1,1}}=-2\int\limits_{0}^{t}{{{\varphi }_{0}}\left( \eta  \right)d\eta }\int\limits_{\eta }^{t}{{{{{W}'}}_{y}}\left( t-y,a-x,a+x \right)\times \,}$$
$$\times \sum\limits_{n=0}^{+\infty }{\frac{{{\left( -1 \right)}^{n}}{{a}^{n}}}{n!}{{\left( y-\eta  \right)}^{-{{\beta }_{1}}-{{\beta }_{1}}n}}E_{\alpha ,1-{{\beta }_{1}}-{{\beta }_{1}}n}^{-{{\gamma }_{1}}-{{\gamma }_{1}}n}\left[ \delta {{\left( y-\eta  \right)}^{\alpha }} \right]}\,dy.$$
We define the required value of ${{{W}'}_{y}}\left( t-y,a-x,a+x \right):$
$${{{W}'}_{y}}\left( t-y,a-x,a+x \right)=-\frac{1}{2}\sum\limits_{m=0}^{+\infty }{\frac{{{\left( -1 \right)}^{m}}\left[ {{\left( a-x \right)}^{m}}-{{\left( a+x \right)}^{m}} \right]}{m!}\times}$$
$${\times{{\left( t-y \right)}^{{{\beta }_{1}}-{{\beta }_{1}}m-2}}E_{\alpha ,{{\beta }_{1}}-{{\beta }_{1}}m-1}^{{{\gamma }_{1}}-{{\gamma }_{1}}m}\left[ \delta {{\left( t-y \right)}^{\alpha }} \right]}.$$
We substitute the calculated value in place of ${{{W}'}_{y}}\left( t-y,a-x,a+x \right)$, and continue with the computation:
$${{F}_{1,1}}=\int\limits_{0}^{t}{{{\varphi }_{0}}\left( \eta  \right)d\eta \int\limits_{\eta }^{t}{\sum\limits_{n=0}^{+\infty }{\frac{{{\left( -1 \right)}^{n}}{{a}^{n}}}{n!}{{\left( y-\eta  \right)}^{-{{\beta }_{1}}-{{\beta }_{1}}n}}E_{\alpha ,1-{{\beta }_{1}}-{{\beta }_{1}}n}^{-{{\gamma }_{1}}-{{\gamma }_{1}}n}\left[ \delta {{\left( y-\eta  \right)}^{\alpha }} \right]}\times \,}}$$
$$\times \sum\limits_{m=0}^{+\infty }{\frac{{{\left( -1 \right)}^{m}}\left[ {{\left( a-x \right)}^{m}}-{{\left( a+x \right)}^{m}} \right]}{m!}\times}$$
$${\times{{\left( t-y \right)}^{{{\beta }_{1}}-{{\beta }_{1}}m-2}}E_{\alpha ,{{\beta }_{1}}-{{\beta }_{1}}m-1}^{{{\gamma }_{1}}-{{\gamma }_{1}}m}\left[ \delta {{\left( t-y \right)}^{\alpha }} \right]}\,dy.$$
First, we apply the formula \eqref{eq2.27} to the sums over $n$ and $m$, and then to the generalized Mittag-Leffler functions:
$${{F}_{1,1}}=\sum\limits_{n=0}^{+\infty }{\sum\limits_{m=0}^{n}{\frac{{{\left( -1 \right)}^{n}}{{a}^{m}}\left[ {{\left( a-x \right)}^{n-m}}-{{\left( a+x \right)}^{n-m}} \right]}{m!\left( n-m \right)!}}}\int\limits_{0}^{t}{{{\varphi }_{0}}\left( \eta  \right)d\eta \times }$$
$$\times \int\limits_{\eta }^{t}{{{\left( y-\eta  \right)}^{-{{\beta }_{1}}-{{\beta }_{1}}m}}{{\left( t-y \right)}^{{{\beta }_{1}}-{{\beta }_{1}}n+{{\beta }_{1}}m-2}}\times }$$
$$\times E_{\alpha ,1-{{\beta }_{1}}-{{\beta }_{1}}m}^{-{{\gamma }_{1}}-{{\gamma }_{1}}m}\left[ \delta {{\left( y-\eta  \right)}^{\alpha }} \right]E_{\alpha ,{{\beta }_{1}}-{{\beta }_{1}}n+{{\beta }_{1}}m-1}^{{{\gamma }_{1}}-{{\gamma }_{1}}n+{{\gamma }_{1}}m}\left[ \delta {{\left( t-y \right)}^{\alpha }} \right]dy=$$
$$=\sum\limits_{n=0}^{+\infty }{\sum\limits_{m=0}^{n}{\frac{{{\left( -1 \right)}^{n}}{{a}^{m}}\left[ {{\left( a-x \right)}^{n-m}}-{{\left( a+x \right)}^{n-m}} \right]}{m!\left( n-m \right)!}}}\sum\limits_{k=0}^{+\infty }{\sum\limits_{z=0}^{k}{\frac{{\delta }^{k}}{z!\left( k-z \right)!}}}\times $$
$${\times\frac{{{\left( -{{\gamma }_{1}}-{{\gamma }_{1}}m \right)}_{z}}{{\left( {{\gamma }_{1}}-{{\gamma }_{1}}n+{{\gamma }_{1}}m \right)}_{k-z}}}{\Gamma \left( \alpha z+1-{{\beta }_{1}}-{{\beta }_{1}}m \right)\Gamma \left( \alpha k-\alpha z+{{\beta }_{1}}-{{\beta }_{1}}n+{{\beta }_{1}}m-1 \right)}\times}$$
$$\times \int\limits_{0}^{t}{{{\varphi }_{0}}\left( \eta  \right)d\eta }\int\limits_{\eta }^{t}{{{\left( y-\eta  \right)}^{\alpha z-{{\beta }_{1}}-{{\beta }_{1}}m}}{{\left( t-y \right)}^{\alpha k-\alpha z+{{\beta }_{1}}-{{\beta }_{1}}n+{{\beta }_{1}}m-2}}dy.}$$
Now we make the substitution $y=\left( t-\eta  \right)s+\eta $ and will use formulas \eqref{eq2.28} and \cite{Graham}
\begin{equation}\label{eq2.29}
 \sum\limits_{m=0}^{n}{\frac{{{a}^{m}}{{b}^{n-m}}}{m!\left( n-m \right)!}=\frac{{{\left( a+b \right)}^{n}}}{n!}}
\end{equation} 
to get
$${{F}_{1,1}}=\sum\limits_{n=0}^{+\infty }{\sum\limits_{m=0}^{n}{\frac{{{\left( -1 \right)}^{n}}{{a}^{m}}\left[ {{\left( a-x \right)}^{n-m}}-{{\left( a+x \right)}^{n-m}} \right]}{m!\left( n-m \right)!}}}\times $$
$$\times \sum\limits_{k=0}^{+\infty }{\frac{{{\left( -{{\gamma }_{1}}n \right)}_{k}}{{\delta }^{k}}}{k!\Gamma \left( \alpha k-{{\beta }_{1}}n \right)}}\int\limits_{0}^{t}{{{\varphi }_{0}}\left( \eta  \right){{\left( t-\eta  \right)}^{\alpha k-{{\beta }_{1}}n-1}}d\eta }=$$
$$=\sum\limits_{n=0}^{+\infty }{\frac{{{\left( -1 \right)}^{n}}\left[ {{\left( 2a-x \right)}^{n}}-{{\left( 2a+x \right)}^{n}} \right]}{n!}}\times $$
$$\times \int\limits_{0}^{t}{{{\varphi }_{0}}\left( \eta  \right){{\left( t-\eta  \right)}^{-{{\beta }_{1}}n-1}}E_{\alpha ,-{{\beta }_{1}}n}^{-{{\gamma }_{1}}n}\left[ \delta {{\left( t-\eta  \right)}^{\alpha }} \right]d\eta }=$$
$$=\int\limits_{0}^{t}{{{\varphi }_{0}}\left( \eta  \right)\left[ \omega \left( t-\eta ,2a-x \right)-\omega \left( t-\eta ,2a+x \right) \right]d\eta }.$$
We now turn our attention to the evaluation of the integral ${{F}_{1,2}}$. Firstly, the order of integration is interchanged, and the formula \eqref{eq2.19} is applied:
$${{F}_{1,2}}=2\int\limits_{0}^{t}{{}^{PRL}D_{0y}^{\alpha ,\,{{\beta }_{1}},\,{{\gamma }_{1}},\,\delta }\left[ \int\limits_{0}^{y}{\sum\limits_{n=1}^{\infty }{\omega \left( y-\eta ,2na \right)}\,d\eta \int\limits_{0}^{\eta }{{{\varphi }_{0}}\left( z \right)\omega \left( \eta -z,a \right)dz}} \right]}\times $$
$$\times W\left( t-y,a-x,a+x \right)dy=$$
$$=2\int\limits_{0}^{t}{{}^{PRL}D_{0y}^{\alpha ,\,{{\beta }_{1}},\,{{\gamma }_{1}},\,\delta }\left[ \int\limits_{0}^{y}{{{\varphi }_{0}}\left( z \right)dz\int\limits_{z}^{y}{\omega \left( \eta -z,a \right)\sum\limits_{n=1}^{\infty }{\omega \left( y-\eta ,2na \right)}\,d\eta }} \right]}\times $$
$$\times W\left( t-y,a-x,a+x \right)dy=$$
$$=2\int\limits_{0}^{t}{{}^{PRL}D_{0y}^{\alpha ,\,{{\beta }_{1}},\,{{\gamma }_{1}},\,\delta }\left\{ \int\limits_{0}^{y}{{{\varphi }_{0}}\left( \eta  \right)\sum\limits_{n=1}^{\infty }{\omega \left( y-\eta ,2na+a \right)}\,d\eta } \right\}}\times $$
$$\times W\left( t-y,a-x,a+x \right)dy.$$
By repeating the calculations performed for determining ${{F}_{1,1}}$, the following result is obtained for ${{F}_{1,2}}$:
$${{F}_{1,2}}=\int\limits_{0}^{t}{{{\varphi }_{0}}\left( \eta  \right)\sum\limits_{n=1}^{\infty }{\left[ \omega \left( t-\eta ,\left( 2n+2 \right)a-x \right)-\omega \left( t-\eta ,\left( 2n+2 \right)a+x \right) \right]\,}d\eta }.$$
Based on the results obtained for ${{F}_{1,1}}$ and ${{F}_{1,2}}$, ${{F}_{1}}$ can be expressed as follows:
$${{F}_{1}}=\int\limits_{0}^{t}{{{\varphi }_{0}}\left( \eta  \right)\omega \left( t-\eta ,x \right)d\eta -}$$
$$-\int\limits_{0}^{t}{{{\varphi }_{0}}\left( \eta  \right)\left[ \omega \left( t-\eta ,2a-x \right)-\omega \left( t-\eta ,2a+x \right) \right]d\eta }-$$
$$-\int\limits_{0}^{t}{{{\varphi }_{0}}\left( \eta  \right)\sum\limits_{n=1}^{\infty }{\left[ \omega \left( t-\eta ,\left( 2n+2 \right)a-x \right)-\omega \left( t-\eta ,\left( 2n+2 \right)a+x \right) \right]\,}d\eta }=$$
$$=\int\limits_{0}^{t}{{{\varphi }_{0}}\left( \eta  \right)\omega \left( t-\eta ,x \right)d\eta -}$$
$$-\int\limits_{0}^{t}{{{\varphi }_{0}}\left( \eta  \right)\sum\limits_{n=0}^{\infty }{\left[ \omega \left( t-\eta ,\left( 2n+2 \right)a-x \right)-\omega \left( t-\eta ,\left( 2n+2 \right)a+x \right) \right]\,}d\eta }=$$
$$=\int\limits_{0}^{t}{{{\varphi }_{0}}\left( \eta  \right)\omega \left( t-\eta ,x \right)d\eta -}$$
$$-\int\limits_{0}^{t}{{{\varphi }_{0}}\left( \eta  \right)\sum\limits_{n=1}^{\infty }{\left[ \omega \left( t-\eta ,2na-x \right)-\omega \left( t-\eta ,2na+x \right) \right]\,}d\eta }=$$
$$=\int\limits_{0}^{t}{{{\varphi }_{0}}\left( \eta  \right)\omega \left( t-\eta ,x \right)d\eta }-\int\limits_{0}^{t}{{{\varphi }_{0}}\left( \eta  \right)\sum\limits_{n=-\infty }^{-1}{\omega \left( t-\eta ,-2na-x \right)\,}d\eta }+$$
$$+\int\limits_{0}^{t}{{{\varphi }_{0}}\left( \eta  \right)\sum\limits_{n=1}^{\infty }{\omega \left( t-\eta ,2na+x \right)\,}d\eta }=$$
$$=\int\limits_{0}^{t}{{{\varphi }_{0}}\left( \eta  \right)\sum\limits_{n=-\infty }^{\infty }{sign\left( x+2na \right)\omega \left( t-\eta ,\left| x+2na \right| \right)\,}d\eta }.$$
By denoting the terms involving ${{\varphi }_{1}}\left( \eta  \right)$ in equation \eqref{eq2.26} as ${{F}_{2}}$, repeating the same steps we applied to ${{F}_{1}}$ and we do some calculations:
$${{F}_{2}}=\int\limits_{0}^{t}{2{}^{PRL}D_{0y}^{\alpha ,\,{{\beta }_{1}},\,{{\gamma }_{1}},\,\delta }\left\{ {{\varphi }_{1}}\left( y \right)+\int\limits_{0}^{y}{{{\varphi }_{1}}\left( \eta  \right)\sum\limits_{n=1}^{\infty }{\omega \left( y-\eta ,2na \right)}\,d\eta } \right\}}\times$$
$$\times W\left( t-y,a-x,a+x \right)dy=\int\limits_{0}^{t}{{{\varphi }_{1}}\left( \eta  \right)\left[ \omega \left( t-\eta ,a-x \right)-\omega \left( t-\eta ,a+x \right) \right]d\eta }+$$
$$+\int\limits_{0}^{t}{{{\varphi }_{1}}\left( \eta  \right)\sum\limits_{n=1}^{\infty }{\left[ \omega \left( t-\eta ,\left( 2n+1 \right)a-x \right)-\omega \left( t-\eta ,\left( 2n+1 \right)a+x \right) \right]\,}d\eta }=$$
$$=-\int\limits_{0}^{t}{{{\varphi }_{1}}\left( \eta  \right)\omega \left( t-\eta ,a+x \right)d\eta }+\int\limits_{0}^{t}{{{\varphi }_{1}}\left( \eta  \right)\sum\limits_{n=0}^{\infty }{\omega \left( t-\eta ,\left( 2n+1 \right)a-x \right)\,}d\eta }-$$
$$-\int\limits_{0}^{t}{{{\varphi }_{1}}\left( \eta  \right)\sum\limits_{n=1}^{\infty }{\omega \left( t-\eta ,\left( 2n+1 \right)a+x \right)\,}d\eta }=-\int\limits_{0}^{t}{{{\varphi }_{1}}\left( \eta  \right)\omega \left( t-\eta ,a+x \right)d\eta }+$$
$$+\int\limits_{0}^{t}{{{\varphi }_{1}}\left( \eta  \right)\sum\limits_{n=1}^{\infty }{\omega \left( t-\eta ,\left( 2n-1 \right)a-x \right)\,}d\eta }-$$
$$-\int\limits_{0}^{t}{{{\varphi }_{1}}\left( \eta  \right)\sum\limits_{n=1}^{\infty }{\omega \left( t-\eta ,\left( 2n+1 \right)a+x \right)\,}d\eta }=$$
$$=-\int\limits_{0}^{t}{{{\varphi }_{1}}\left( \eta  \right)\omega \left( t-\eta ,a+x \right)d\eta }+$$
$$+\int\limits_{0}^{t}{{{\varphi }_{1}}\left( \eta  \right)\sum\limits_{n=-\infty }^{-1}{\omega \left( t-\eta ,\left( -2n-1 \right)a-x \right)\,}d\eta }-$$
$$-\int\limits_{0}^{t}{{{\varphi }_{1}}\left( \eta  \right)\sum\limits_{n=1}^{\infty }{\omega \left( t-\eta ,\left( 2n+1 \right)a+x \right)\,}d\eta }=$$
$$=-\int\limits_{0}^{t}{{{\varphi }_{1}}\left( \eta  \right)\sum\limits_{n=-\infty }^{+\infty }{sign\left( x+\left( 2n+1 \right)a \right)\omega \left( t-\eta ,\left| x+\left( 2n+1 \right)a \right| \right)\,}d\eta }.$$

Let the terms involving $\tau \left( s \right)$ be denoted by ${{F}_{3}}$ in \eqref{eq2.26}:
$${{F}_{3}}=-2\int\limits_{0}^{t}{{}^{PRL}D_{0y}^{\alpha ,\,{{\beta }_{1}},\,{{\gamma }_{1}},\,\delta }}\left\{ \int\limits_{0}^{a}{\tau \left( s \right)W\left( y,a-s,a+s \right)ds+} \right.$$
$$\left. +\int\limits_{0}^{y}{\sum\limits_{n=1}^{\infty }{\omega \left( y-\eta ,2na \right)}\,d\eta \int\limits_{0}^{a}{\tau \left( s \right)W\left( \eta ,a-s,a+s \right)ds}} \right\}\times$$
$$\times W\left( t-y,a-x,a+x \right)dy+\int\limits_{0}^{a}{\tau \left( s \right)W\left( t,\left| s-x \right|,s+x \right)ds.}$$
Now, we denote the first integral as ${{F}_{3,1}}$ and proceed with its evaluation:
$${{F}_{3,1}}=-\int\limits_{0}^{t}{2{}^{PRL}D_{0y}^{\alpha ,\,{{\beta }_{1}},\,{{\gamma }_{1}},\,\delta }}\left\{ \int\limits_{0}^{a}{\tau \left( s \right)W\left( y,a-s,a+s \right)ds} \right\}\times$$
$$\times W\left( t-y,a-x,a+x \right)dy=-2\int\limits_{0}^{t}{W\left( t-y,a-x,a+x \right)\times}$$
$${\times\frac{\partial }{\partial y}{}^{P}I_{0y}^{\alpha ,\,1-{{\beta }_{1}},\,-{{\gamma }_{1}},\,\delta }}\left\{ \int\limits_{0}^{a}{\tau \left( s \right)W\left( y,a-s,a+s \right)ds} \right\}dy=$$
$$=-2\int\limits_{0}^{t}{W\left( t-y,a-x,a+x \right)\times }$$
$$\times \frac{\partial }{\partial y}\left[ \int\limits_{0}^{y}{{{\left( y-z \right)}^{-{{\beta }_{1}}}}E_{\alpha ,\,1-{{\beta }_{1}}}^{-{{\gamma }_{1}}}}\left[ \delta {{\left( y-z \right)}^{\alpha }} \right]dz\int\limits_{0}^{a}{\tau \left( s \right)W\left( z,a-s,a+s \right)ds} \right]dy.$$
Next, we interchange the order of integration with respect to the variables $s$ and $z$, substitute the value of $W\left( z,a-s,a+s \right)$:
$$\int\limits_{0}^{y}{{{\left( y-z \right)}^{-{{\beta }_{1}}}}E_{\alpha ,\,1-{{\beta }_{1}}}^{-{{\gamma }_{1}}}}\left[ \delta {{\left( y-z \right)}^{\alpha }} \right]dz\int\limits_{0}^{a}{\tau \left( s \right)W\left( z,a-s,a+s \right)ds}=$$
$$=\int\limits_{0}^{a}{\tau \left( s \right)ds}\int\limits_{0}^{y}{W\left( z,a-s,a+s \right){{\left( y-z \right)}^{-{{\beta }_{1}}}}E_{\alpha ,\,1-{{\beta }_{1}}}^{-{{\gamma }_{1}}}\left[ \delta {{\left( y-z \right)}^{\alpha }} \right]dz}=$$
$$=\frac{1}{2}\int\limits_{0}^{a}{\tau \left( s \right)ds}\int\limits_{0}^{y}{\sum\limits_{n=0}^{+\infty }{\frac{{{\left( -1 \right)}^{n}}\left[ {{\left( a-s \right)}^{n}}-{{\left( a+s \right)}^{n}} \right]}{n!}\times }}$$
$$\times {{z}^{{{\beta }_{1}}-{{\beta }_{1}}n-1}}E_{\alpha ,{{\beta }_{1}}-{{\beta }_{1}}n}^{{{\gamma }_{1}}-{{\gamma }_{1}}n}\left[ \delta {{z}^{\alpha }} \right]{{\left( y-z \right)}^{-{{\beta }_{1}}}}E_{\alpha ,\,1-{{\beta }_{1}}}^{-{{\gamma }_{1}}}\left[ \delta {{\left( y-z \right)}^{\alpha }} \right]\,dz.$$
We continue the calculation by using \eqref{eq2.27} and \eqref{eq2.28}:
$$\frac{1}{2}\int\limits_{0}^{a}{\sum\limits_{n=0}^{+\infty }{\frac{{{\left( -1 \right)}^{n}}\left[ {{\left( a-s \right)}^{n}}-{{\left( a+s \right)}^{n}} \right]}{n!}}\,}\times $$
$$\times \sum\limits_{k=0}^{+\infty }{\sum\limits_{m=0}^{k}{\frac{{{\left( {{\gamma }_{1}}-{{\gamma }_{1}}n \right)}_{m}}{{\left( -{{\gamma }_{1}} \right)}_{k-m}}{{\delta }^{k}}}{k!\left( m-k \right)!\Gamma \left( \alpha m+{{\beta }_{1}}-{{\beta }_{1}}n \right)\Gamma \left( \alpha k-\alpha m+1-{{\beta }_{1}} \right)}}}\,\tau \left( s \right)ds\times $$
$$\times \int\limits_{0}^{y}{{{z}^{\alpha m+{{\beta }_{1}}-{{\beta }_{1}}n-1}}{{\left( y-z \right)}^{\alpha k-\alpha m-{{\beta }_{1}}}}\,dz}=\int\limits_{0}^{a}{\sum\limits_{n=0}^{+\infty }{\frac{{{\left( -1 \right)}^{n}}\left[ {{\left( a-s \right)}^{n}}-{{\left( a+s \right)}^{n}} \right]}{2n!}}\,}\times $$
$$\times \sum\limits_{k=0}^{+\infty }{\frac{{{\delta }^{k}}{{y}^{\alpha k-{{\beta }_{1}}n}}}{\Gamma \left( \alpha k+1-{{\beta }_{1}}n \right)}\sum\limits_{m=0}^{k}{\frac{{{\left( {{\gamma }_{1}}-{{\gamma }_{1}}n \right)}_{m}}{{\left( -{{\gamma }_{1}} \right)}_{k-m}}}{k!\left( m-k \right)!}}}\,\tau \left( s \right)ds=$$
$$=\frac{1}{2}\int\limits_{0}^{a}{\sum\limits_{n=0}^{+\infty }{\frac{{{\left( -1 \right)}^{n}}\left[ {{\left( a-s \right)}^{n}}-{{\left( a+s \right)}^{n}} \right]}{n!}}\,}\sum\limits_{k=0}^{+\infty }{\frac{{{\left( -{{\gamma }_{1}}n \right)}_{k}}{{\delta }^{k}}{{y}^{\alpha k-{{\beta }_{1}}n}}}{k!\Gamma \left( \alpha k+1-{{\beta }_{1}}n \right)}}\,\tau \left( s \right)ds.$$
We take the derivative with respect to $y$ from the last obtained result and change the order of integration in ${{F}_{3,1}}$:
$${{F}_{3,1}}=-\int\limits_{0}^{t}{W\left( t-y,a-x,a+x \right)dy\times }$$
$$\times \int\limits_{0}^{a}{\sum\limits_{n=0}^{+\infty }{\frac{{{\left( -1 \right)}^{n}}\left[ {{\left( a-s \right)}^{n}}-{{\left( a+s \right)}^{n}} \right]}{n!}}\,}\sum\limits_{k=0}^{+\infty }{\frac{{{\left( -{{\gamma }_{1}}n \right)}_{k}}{{\delta }^{k}}{{y}^{\alpha k-{{\beta }_{1}}n-1}}}{k!\Gamma \left( \alpha k-{{\beta }_{1}}n \right)}}\,\tau \left( s \right)ds=$$
$$=-\frac{1}{2}\int\limits_{0}^{a}{\tau \left( s \right)ds\int\limits_{0}^{t}{\sum\limits_{n=0}^{+\infty }{\frac{{{\left( -1 \right)}^{n}}\left[ {{\left( a-s \right)}^{n}}-{{\left( a+s \right)}^{n}} \right]}{n!}}{{y}^{-{{\beta }_{1}}n-1}}E_{\alpha ,-{{\beta }_{1}}n}^{-{{\gamma }_{1}}n}\left[ \delta {{y}^{\alpha }} \right]\times }}$$
$$\times \sum\limits_{m=0}^{+\infty }{\frac{{{\left( -1 \right)}^{m}}\left[ {{\left( a-x \right)}^{m}}-{{\left( a+x \right)}^{m}} \right]}{m!}{{\left( t-y \right)}^{{{\beta }_{1}}-{{\beta }_{1}}m-1}}E_{\alpha ,{{\beta }_{1}}-{{\beta }_{1}}m}^{{{\gamma }_{1}}-{{\gamma }_{1}}m}\left[ \delta {{\left( t-y \right)}^{\alpha }} \right]dy.}$$
We then apply formula \eqref{eq2.27} repeatedly, twice and get
$${{F}_{3,1}}=-\frac{1}{2}\int\limits_{0}^{a}{\tau \left( s \right)ds\int\limits_{0}^{t}{\sum\limits_{n=0}^{+\infty }{\sum\limits_{m=0}^{n}{\frac{{{\left( -1 \right)}^{n}}\left[ {{\left( a-s \right)}^{m}}-{{\left( a+s \right)}^{m}} \right]}{m!}}}\times}}$$
$${\times{{y}^{-{{\beta }_{1}}m-1}}E_{\alpha ,-{{\beta }_{1}}m}^{-{{\gamma }_{1}}m}\left[ \delta {{y}^{\alpha }} \right] } \frac{\left[ {{\left( a-x \right)}^{n-m}}-{{\left( a+x \right)}^{n-m}} \right]}{\left( n-m \right)!}\times$$
$$\times{{\left( t-y \right)}^{{{\beta }_{1}}-{{\beta }_{1}}\left( n-m \right)-1}}E_{\alpha ,{{\beta }_{1}}-{{\beta }_{1}}\left( n-m \right)}^{{{\gamma }_{1}}-{{\gamma }_{1}}\left( n-m \right)}\left[ \delta {{\left( t-y \right)}^{\alpha }} \right]dy=$$
$$=-\frac{1}{2}\int\limits_{0}^{a}{\tau \left( s \right)\sum\limits_{n=0}^{+\infty }{\sum\limits_{m=0}^{n}{\frac{{{\left( -1 \right)}^{n}}\left[ {{\left( a-s \right)}^{m}}-{{\left( a+s \right)}^{m}} \right]}{m!}}}\times}$$
$${\times\frac{\left[ {{\left( a-x \right)}^{n-m}}-{{\left( a+x \right)}^{n-m}} \right]}{\left( n-m \right)!}ds } \int\limits_{0}^{t}{{{y}^{-{{\beta }_{1}}m-1}}{{\left( t-y \right)}^{{{\beta }_{1}}-{{\beta }_{1}}n+{{\beta }_{1}}m-1}}\times }$$
$${\times E_{\alpha ,-{{\beta }_{1}}m}^{-{{\gamma }_{1}}m}\left[ \delta {{y}^{\alpha }} \right]}E_{\alpha ,{{\beta }_{1}}-{{\beta }_{1}}n+{{\beta }_{1}}m}^{{{\gamma }_{1}}-{{\gamma }_{1}}n+{{\gamma }_{1}}m}\left[ \delta {{\left( t-y \right)}^{\alpha }} \right]dy=-\frac{1}{2}\int\limits_{0}^{a}{\tau \left( s \right)\times}$$
$${\times\sum\limits_{n=0}^{+\infty }{\sum\limits_{m=0}^{n}{\frac{{{\left( -1 \right)}^{n}}\left[ {{\left( a-s \right)}^{m}}-{{\left( a+s \right)}^{m}} \right]}{m!}}}\frac{\left[ {{\left( a-x \right)}^{n-m}}-{{\left( a+x \right)}^{n-m}} \right]}{\left( n-m \right)!}\times }$$
$$\times \sum\limits_{k=0}^{+\infty }{\sum\limits_{l=0}^{k}{\frac{{{\left( -{{\gamma }_{1}}m \right)}_{l}}{{\left( {{\gamma }_{1}}-{{\gamma }_{1}}n+{{\gamma }_{1}}m \right)}_{k-l}}{{\delta }^{k}}}{l!\left( k-l \right)!\Gamma \left( \alpha l-{{\beta }_{1}}m \right)\Gamma \left( \alpha k-\alpha l+{{\beta }_{1}}-{{\beta }_{1}}n+{{\beta }_{1}}m \right)}}}\,ds\times $$
$$\times \int\limits_{0}^{t}{{{y}^{\alpha l-{{\beta }_{1}}m-1}}{{\left( t-y \right)}^{\alpha k-\alpha l+{{\beta }_{1}}-{{\beta }_{1}}n+{{\beta }_{1}}m-1}}}dy.$$
After performing the substitution $y=tz$, we use the formulas \eqref{eq2.28} and \eqref{eq2.29}:
$${{F}_{3,1}}=-\frac{1}{2}\int\limits_{0}^{a}{\tau \left( s \right)\sum\limits_{n=0}^{+\infty }{\sum\limits_{m=0}^{n}{\frac{{{\left( -1 \right)}^{n}}\left[ {{\left( a-s \right)}^{m}}-{{\left( a+s \right)}^{m}} \right]}{m!}}}\times}$$
$${\times\frac{\left[ {{\left( a-x \right)}^{n-m}}-{{\left( a+x \right)}^{n-m}} \right]}{\left( n-m \right)!}}\sum\limits_{k=0}^{+\infty }{\frac{{{\delta }^{k}}{{t}^{\alpha k+{{\beta }_{1}}-{{\beta }_{1}}n-1}}}{\Gamma \left( \alpha k+{{\beta }_{1}}-{{\beta }_{1}}n \right)}\times} $$
$$\times \sum\limits_{l=0}^{k}{\frac{{{\left( -{{\gamma }_{1}}m \right)}_{l}}{{\left( {{\gamma }_{1}}-{{\gamma }_{1}}n+{{\gamma }_{1}}m \right)}_{k-l}}}{l!\left( k-l \right)!}}\,ds=-\frac{1}{2}\int\limits_{0}^{a}{\tau \left( s \right)\times}$$
$${\times\sum\limits_{n=0}^{+\infty }{\sum\limits_{m=0}^{n}{\frac{{{\left( -1 \right)}^{n}}\left[ {{\left( a-s \right)}^{m}}-{{\left( a+s \right)}^{m}} \right]}{m!}}}\frac{\left[ {{\left( a-x \right)}^{n-m}}-{{\left( a+x \right)}^{n-m}} \right]}{\left( n-m \right)!}\times }$$
$$\times \sum\limits_{k=0}^{+\infty }{\frac{{{\delta }^{k}}{{\left( {{\gamma }_{1}}-{{\gamma }_{1}}n \right)}_{k}}{{t}^{\alpha k+{{\beta }_{1}}-{{\beta }_{1}}n-1}}}{k!\Gamma \left( \alpha k+{{\beta }_{1}}-{{\beta }_{1}}n \right)}}\,ds=-\frac{1}{2}\int\limits_{0}^{a}{\tau \left( s \right)\times}$$
$${\times \sum\limits_{n=0}^{+\infty }{\frac{{{\left( -1 \right)}^{n}}\left[ {{\left( 2a-x-s \right)}^{n}}-{{\left( 2a-x+s \right)}^{n}}-{{\left( 2a+x-s \right)}^{n}}+{{\left( 2a+x+s \right)}^{n}} \right]}{n!}}\times }$$
$$\times \sum\limits_{k=0}^{+\infty }{\frac{{{\delta }^{k}}{{\left( {{\gamma }_{1}}-{{\gamma }_{1}}n \right)}_{k}}{{t}^{\alpha k+{{\beta }_{1}}-{{\beta }_{1}}n-1}}}{k!\Gamma \left( \alpha k+{{\beta }_{1}}-{{\beta }_{1}}n \right)}}\,ds.$$
We express the final result in terms of the generalized Mittag-Leffler function and then represent it using the Prabhakar integral:
$${{F}_{3,1}}=\frac{1}{2}\int\limits_{0}^{a}{\tau \left( s \right)\times}$$
$$\times{\sum\limits_{n=0}^{+\infty }{\frac{{{\left( -1 \right)}^{n}}\left[ {{\left( 2a-x+s \right)}^{n}}-{{\left( 2a-x-s \right)}^{n}} \right]}{n!}{{t}^{{{\beta }_{1}}-{{\beta }_{1}}n-1}}}}E_{\alpha ,{{\beta }_{1}}-{{\beta }_{1}}n}^{{{\gamma }_{1}}-{{\gamma }_{1}}n}\left[ \delta {{t}^{\alpha }} \right]\,ds=$$
$$+\frac{1}{2}\int\limits_{0}^{a}{\tau \left( s \right)\sum\limits_{n=0}^{+\infty }{\frac{{{\left( -1 \right)}^{n}}\left[ {{\left( 2a+x-s \right)}^{n}}-{{\left( 2a+x+s \right)}^{n}} \right]}{n!}}}\times$$
$$\times{{t}^{{{\beta }_{1}}-{{\beta }_{1}}n-1}}E_{\alpha ,{{\beta }_{1}}-{{\beta }_{1}}n}^{{{\gamma }_{1}}-{{\gamma }_{1}}n}\left[ \delta {{t}^{\alpha }} \right]\,ds=$$
$$=\frac{1}{2}\int\limits_{0}^{a}{\tau \left( s \right){}^{P}I_{0t}^{\alpha ,\,{{\beta }_{1}},\,{{\gamma }_{1}},\,\delta }\left[ \omega \left( t,2a-x+s \right)-\omega \left( t,2a-x-s \right) \right]}\,ds+$$
$$+\frac{1}{2}\int\limits_{0}^{a}{\tau \left( s \right){}^{P}I_{0t}^{\alpha ,\,{{\beta }_{1}},\,{{\gamma }_{1}},\,\delta }\left[ \omega \left( t,2a+x-s \right)-\omega \left( t,2a+x+s \right) \right]}\,ds.$$
We designate the term containing $\sum\limits_{n=1}^{\infty }{\omega \left( y-\eta ,2na \right)}$ as ${{F}_{3,2}}$:
$${{F}_{3,2}}=-\int\limits_{0}^{t}{2{}^{PRL}D_{0y}^{\alpha ,\,{{\beta }_{1}},\,{{\gamma }_{1}},\,\delta }}\left\{ \int\limits_{0}^{y}{\sum\limits_{n=1}^{\infty }{\omega \left( y-\eta ,2na \right)}\,d\eta\times}\right.$$
$$\left.\times \int\limits_{0}^{a}{\tau \left( s \right)W\left( \eta ,a-s,a+s \right)ds} \right\}W\left( t-y,a-x,a+x \right)dy.$$
Below, we first change the order of integration and then substitute the values of the functions $\omega \left( y-\eta ,2na \right)$ and $W\left( \eta ,a-s,a+s \right)$ to perform the calculation:
$$\int\limits_{0}^{y}{\sum\limits_{n=1}^{\infty }{\omega \left( y-\eta ,2na \right)}\,d\eta \int\limits_{0}^{a}{\tau \left( s \right)W\left( \eta ,a-s,a+s \right)ds}}=$$
$$=\int\limits_{0}^{a}{\tau \left( s \right)ds}\int\limits_{0}^{y}{\sum\limits_{n=1}^{\infty }{\omega \left( y-\eta ,2na \right)}\,W\left( \eta ,a-s,a+s \right)d\eta }=$$
$$=\frac{1}{2}\int\limits_{0}^{a}{\tau \left( s \right)ds}\int\limits_{0}^{y}{\sum\limits_{n=1}^{\infty }{\sum\limits_{k=0}^{\infty }{\frac{{{\left( -1 \right)}^{k}}{{\left( 2na \right)}^{k}}}{k!}{{\left( y-\eta  \right)}^{-{{\beta }_{1}}k-1}}E_{\alpha ,-{{\beta }_{1}}k}^{-{{\gamma }_{1}}k}\left[ \delta {{\left( y-\eta  \right)}^{\alpha }} \right]}}\,\times }$$
$$\,\times \sum\limits_{m=0}^{\infty }{\frac{{{\left( -1 \right)}^{m}}\left[ {{\left( a-s \right)}^{m}}-{{\left( a+s \right)}^{m}} \right]}{m!}{{\eta }^{{{\beta }_{1}}-{{\beta }_{1}}m-1}}E_{\alpha ,{{\beta }_{1}}-{{\beta }_{1}}m}^{{{\gamma }_{1}}-{{\gamma }_{1}}m}\left[ \delta {{\eta }^{\alpha }} \right]}\,d\eta .$$
Next, we apply formulas \eqref{eq2.27}-\eqref{eq2.29} and arrive at the following result:
$$\frac{1}{2}\int\limits_{0}^{a}{\tau \left( s \right)\sum\limits_{n=1}^{\infty }{\sum\limits_{k=0}^{\infty }{\sum\limits_{m=0}^{k}{\frac{{{\left( -1 \right)}^{k}}{{\left( 2na \right)}^{m}}\left[ {{\left( a-s \right)}^{k-m}}-{{\left( a+s \right)}^{k-m}} \right]}{m!\left( k-m \right)!}}\,}}ds}\times $$
$$\times \int\limits_{0}^{y}{{{\left( y-\eta  \right)}^{-{{\beta }_{1}}m-1}}{{\eta }^{{{\beta }_{1}}-{{\beta }_{1}}k+{{\beta }_{1}}m-1}}}\times$$
$${\times E_{\alpha ,-{{\beta }_{1}}m}^{-{{\gamma }_{1}}m}\left[ \delta {{\left( y-\eta  \right)}^{\alpha }} \right]}\, E_{\alpha ,{{\beta }_{1}}-{{\beta }_{1}}k+{{\beta }_{1}}m}^{{{\gamma }_{1}}-{{\gamma }_{1}}k+{{\gamma }_{1}}m}\left[ \delta {{\eta }^{\alpha }} \right]\,d\eta=$$ 
$$=\frac{1}{2}\int\limits_{0}^{a}{\tau \left( s \right)\sum\limits_{n=1}^{\infty }{\sum\limits_{k=0}^{\infty }{\frac{{{\left( -1 \right)}^{k}}\left[ {{\left( 2na+a-s \right)}^{k}}-{{\left( 2na+a+s \right)}^{k}} \right]}{k!}\,}}}\times $$
$$\times {{y}^{{{\beta }_{1}}-{{\beta }_{1}}k-1}}E_{\alpha ,{{\beta }_{1}}-{{\beta }_{1}}k}^{{{\gamma }_{1}}-{{\gamma }_{1}}k}\left[ \delta {{y}^{\alpha }} \right]ds=$$
$$=\int\limits_{0}^{a}{\tau \left( s \right)\sum\limits_{n=1}^{\infty }{W\left( y,\left( 2n+1 \right)a-s,\left( 2n+1 \right)a+s \right)ds}}.$$
We substitute the obtained final result into ${{F}_{3,2}}$, and since it is similar to ${{F}_{3,1}}$, we express its result using ${{F}_{3,1}}$:
$${{F}_{3,2}}=-\int\limits_{0}^{t}{2{}^{PRL}D_{0y}^{\alpha ,\,{{\beta }_{1}},\,{{\gamma }_{1}},\,\delta }}\left\{ \int\limits_{0}^{a}\tau \left( s \right)\times \right.$$
$$\left.\times\sum\limits_{n=1}^{\infty }{W\left( y,\left( 2n+1 \right)a-s,\left( 2n+1 \right)a+s \right)ds} \right\} W\left( t-y,a-x,a+x \right)dy=$$
$$=\frac{1}{2}\int\limits_{0}^{a}{\tau \left( s \right)\times}$$
$$\times\sum\limits_{n=1}^{+\infty }{{}^{P}I_{0t}^{\alpha ,\,{{\beta }_{1}},\,{{\gamma }_{1}},\,\delta }\left[ \omega \left( t,\left( 2n+2 \right)a-x+s \right)-\omega \left( t,\left( 2n+2 \right)a-x-s \right) \right]}\,ds+$$
$$+\frac{1}{2}\int\limits_{0}^{a}{\tau \left( s \right)\times}$$
$${\times\sum\limits_{n=1}^{+\infty }{{}^{P}I_{0t}^{\alpha ,\,{{\beta }_{1}},\,{{\gamma }_{1}},\,\delta }\left[ \omega \left( t,\left( 2n+2 \right)a+x-s \right)-\omega \left( t,\left( 2n+2 \right)a+x+s \right) \right]}}\,ds.$$
Next, we determine the value of ${{F}_{3}}$ based on the results of ${{F}_{3,1}}$ and ${{F}_{3,2}}$:
$${{F}_{3}}=\frac{1}{2}\int\limits_{0}^{a}{\tau \left( s \right){}^{P}I_{0t}^{\alpha ,\,{{\beta }_{1}},\,{{\gamma }_{1}},\,\delta }\left[ \omega \left( t,2a-x+s \right)-\omega \left( t,2a-x-s \right) \right]}\,ds+$$
$$+\frac{1}{2}\int\limits_{0}^{a}{\tau \left( s \right){}^{P}I_{0t}^{\alpha ,\,{{\beta }_{1}},\,{{\gamma }_{1}},\,\delta }\left[ \omega \left( t,2a+x-s \right)-\omega \left( t,2a+x+s \right) \right]}\,ds+$$
$$+\frac{1}{2}\int\limits_{0}^{a}{\tau \left( s \right)\sum\limits_{n=1}^{+\infty }{{}^{P}I_{0t}^{\alpha ,\,{{\beta }_{1}},\,{{\gamma }_{1}},\,\delta }\left[ \omega \left( t,\left( 2n+2 \right)a-x+s \right)-\omega \left( t,\left( 2n+2 \right)a-x-s \right) \right]}}\,ds+$$
$$+\frac{1}{2}\int\limits_{0}^{a}{\tau \left( s \right)\sum\limits_{n=1}^{+\infty }{{}^{P}I_{0t}^{\alpha ,\,{{\beta }_{1}},\,{{\gamma }_{1}},\,\delta }\left[ \omega \left( t,\left( 2n+2 \right)a+x-s \right)-\omega \left( t,\left( 2n+2 \right)a+x+s \right) \right]}}\,ds+$$
$$+\frac{1}{2}\int\limits_{0}^{a}{\tau \left( s \right){}^{P}I_{0t}^{\alpha ,\,{{\beta }_{1}},\,{{\gamma }_{1}},\,\delta }\left[ \omega \left( t,\left| s-x \right| \right)-\omega \left( t,s+x \right) \right]}\,ds=$$
$$=\frac{1}{2}\int\limits_{0}^{a}{\tau \left( s \right)\sum\limits_{n=0}^{+\infty }{{}^{P}I_{0t}^{\alpha ,\,{{\beta }_{1}},\,{{\gamma }_{1}},\,\delta }\left[ \omega \left( t,\left( 2n+2 \right)a-x+s \right)-\omega \left( t,\left( 2n+2 \right)a-x-s \right) \right]}}\,ds+$$
$$+\frac{1}{2}\int\limits_{0}^{a}{\tau \left( s \right)\sum\limits_{n=0}^{+\infty }{{}^{P}I_{0t}^{\alpha ,\,{{\beta }_{1}},\,{{\gamma }_{1}},\,\delta }\left[ \omega \left( t,\left( 2n+2 \right)a+x-s \right)-\omega \left( t,\left( 2n+2 \right)a+x+s \right) \right]}}\,ds+$$
$$+\frac{1}{2}\int\limits_{0}^{a}{\tau \left( s \right){}^{P}I_{0t}^{\alpha ,\,{{\beta }_{1}},\,{{\gamma }_{1}},\,\delta }\left[ \omega \left( t,\left| s-x \right| \right)-\omega \left( t,s+x \right) \right]}\,ds=$$
$$=\frac{1}{2}\int\limits_{0}^{a}{\tau \left( s \right)\sum\limits_{n=1}^{+\infty }{{}^{P}I_{0t}^{\alpha ,\,{{\beta }_{1}},\,{{\gamma }_{1}},\,\delta }\left[ \omega \left( t,2na-x+s \right)-\omega \left( t,2na-x-s \right) \right]}}\,ds+$$
$$+\frac{1}{2}\int\limits_{0}^{a}{\tau \left( s \right)\sum\limits_{n=1}^{+\infty }{{}^{P}I_{0t}^{\alpha ,\,{{\beta }_{1}},\,{{\gamma }_{1}},\,\delta }\left[ \omega \left( t,2na+x-s \right)-\omega \left( t,2na+x+s \right) \right]}}\,ds+$$
$$+\frac{1}{2}\int\limits_{0}^{a}{\tau \left( s \right){}^{P}I_{0t}^{\alpha ,\,{{\beta }_{1}},\,{{\gamma }_{1}},\,\delta }\left[ \omega \left( t,\left| s-x \right| \right)-\omega \left( t,s+x \right) \right]}\,ds=$$
$$=\frac{1}{2}\int\limits_{0}^{a}{\tau \left( s \right)\sum\limits_{n=-\infty }^{-1}{{}^{P}I_{0t}^{\alpha ,\,{{\beta }_{1}},\,{{\gamma }_{1}},\,\delta }\left[ \omega \left( t,-2na-x+s \right)-\omega \left( t,-2na-x-s \right) \right]}}\,ds+$$
$$+\frac{1}{2}\int\limits_{0}^{a}{\tau \left( s \right)\sum\limits_{n=1}^{+\infty }{{}^{P}I_{0t}^{\alpha ,\,{{\beta }_{1}},\,{{\gamma }_{1}},\,\delta }\left[ \omega \left( t,2na+x-s \right)-\omega \left( t,2na+x+s \right) \right]}}\,ds+$$
$$+\frac{1}{2}\int\limits_{0}^{a}{\tau \left( s \right){}^{P}I_{0t}^{\alpha ,\,{{\beta }_{1}},\,{{\gamma }_{1}},\,\delta }\left[ \omega \left( t,\left| s-x \right| \right)-\omega \left( t,s+x \right) \right]}\,ds=$$
$$=\frac{1}{2}\int\limits_{0}^{a}{\tau \left( s \right)\sum\limits_{n=-\infty }^{+\infty }{{}^{P}I_{0t}^{\alpha ,\,{{\beta }_{1}},\,{{\gamma }_{1}},\,\delta }\left[ \omega \left( t,\left| x-s+2na \right| \right)-\omega \left( t,\left| x+s+2na \right| \right) \right]}}\,ds.$$
Finally, we label the terms containing $f\left( y,s \right)$ in \eqref{eq2.26} as ${{F}_{4}}$:
$${{F}_{4}}=-\int\limits_{0}^{t}{2{}^{PRL}D_{0y}^{\alpha ,\,{{\beta }_{1}},\,{{\gamma }_{1}},\,\delta }\left\{ \int\limits_{0}^{a}{\int\limits_{0}^{y}{f\left( z,s \right)W\left( y-z,a-s,a+s \right)dzds}} \right.}+$$
$$\left. +\int\limits_{0}^{y}{\sum\limits_{n=1}^{\infty }{\omega \left( y-\eta ,2na \right)}\,d\eta }\int\limits_{0}^{a}{\int\limits_{0}^{\eta }{f\left( y,s \right)W\left( \eta -y,a-s,a+s \right)dyds}} \right\}\times $$
$$\times W\left( t-y,a-x,a+x \right)dy+\int\limits_{0}^{a}{\int\limits_{0}^{t}{f\left( y,s \right)W\left( t-y,\left| s-x \right|,s+x \right)dyds.}}$$
If we calculate this expression in the same way as ${{F}_{1}}$, we obtain the following result.
$${{F}_{4}}=\frac{1}{2}\int\limits_{0}^{a}{\int\limits_{0}^{t}{f\left( z,s \right)}\times}$$
$${\times\sum\limits_{n=-\infty }^{+\infty }{{}^{P}I_{zt}^{\alpha ,\,{{\beta }_{1}},\,{{\gamma }_{1}},\,\delta }\left[ \omega \left( t-z,\left| x-s+2na \right| \right)-\omega \left( t-z,\left| x+s+2na \right| \right) \right]}}\,dzds.$$
In ${{F}_{4}}$, we substitute $t-z=\eta $ and obtain the final result:
$${{F}_{4}}=\frac{1}{2}\int\limits_{0}^{t}{\int\limits_{0}^{a}{f\left( t-\eta ,s \right)}}\times$$
$${\times\sum\limits_{n=-\infty }^{+\infty }{{}^{P}I_{0\eta }^{\alpha ,\,{{\beta }_{1}},\,{{\gamma }_{1}},\,\delta }\left[ \omega \left( \eta ,\left| x-s+2na \right| \right)-\omega \left( \eta ,\left| x+s+2na \right| \right) \right]}}\,dsd\eta .$$
Now, by substituting the results of ${{F}_{1}},$ ${{F}_{2}},$ ${{F}_{3}}$ and ${{F}_{4}}$ into \eqref{eq2.26}, we obtain the following equality:
$$u\left( t,x \right)=\int\limits_{0}^{t}{{{\varphi }_{0}}\left( \eta  \right)\sum\limits_{n=-\infty }^{\infty }{sign\left( x+2na \right)\omega \left( t-\eta ,\left| x+2na \right| \right)\,}d\eta }-$$
$$-\int\limits_{0}^{t}{{{\varphi }_{1}}\left( \eta  \right)\sum\limits_{n=-\infty }^{+\infty }{sign\left( x+\left( 2n+1 \right)a \right)\omega \left( t-\eta ,\left| x+\left( 2n+1 \right)a \right| \right)\,}d\eta }+$$
$$+\frac{1}{2}\int\limits_{0}^{a}{\tau \left( s \right)\sum\limits_{n=-\infty }^{+\infty }{{}^{P}I_{0t}^{\alpha ,\,{{\beta }_{1}},\,{{\gamma }_{1}},\,\delta }\left[ \omega \left( t,\left| x-s+2na \right| \right)-\omega \left( t,\left| x+s+2na \right| \right) \right]}}\,ds+$$
$$+\frac{1}{2}\int\limits_{0}^{t}{\int\limits_{0}^{a}{f\left( t-\eta ,s \right)}\times}$$
\begin{equation}\label{eq2.30}
    {\times\sum\limits_{n=-\infty }^{+\infty }{{}^{P}I_{0\eta }^{\alpha ,\,{{\beta }_{1}},\,{{\gamma }_{1}},\,\delta }\left[ \omega \left( \eta ,\left| x-s+2na \right| \right)-\omega \left( \eta ,\left| x+s+2na \right| \right) \right]}}\,dsd\eta.
\end{equation}
We express the kernel of the fourth integral in \eqref{eq2.30} in terms of the function ${{E}_{12}}$:
$$\frac{1}{2}\sum\limits_{n=-\infty }^{+\infty }{{}^{P}I_{0\eta }^{\alpha ,\,{{\beta }_{1}},\,{{\gamma }_{1}},\,\delta }\left[ \omega \left( \eta ,\left| x-s+2na \right| \right)-\omega \left( \eta ,\left| x+s+2na \right| \right) \right]}=$$
$$=\sum\limits_{n=-\infty }^{+\infty }{W\left( \eta ,\left| x-s+2na \right|,\left| x+s+2na \right| \right)}=$$
$$=\frac{1}{2}{{\eta}^{{{\beta }_{1}}-1}}\sum\limits_{n=-\infty }^{+\infty }{\sum\limits_{k=0}^{+\infty }{\frac{{{\left( -1 \right)}^{k}}\left[ {{\left| x-s+2na \right|}^{k}}-{{\left| x+s+2na \right|}^{k}} \right]}{k!}\times }}$$
$$\times {{\eta}^{-{{\beta }_{1}}k}}E_{\alpha ,{{\beta }_{1}}-{{\beta }_{1}}k}^{{{\gamma }_{1}}-{{\gamma }_{1}}k}\left[ \delta {{\eta}^{\alpha }} \right]=\frac{1}{2}{{\eta}^{{{\beta }_{1}}-1}}\times$$
$$\times\sum\limits_{n=-\infty }^{+\infty }{\sum\limits_{k=0}^{+\infty }{\sum\limits_{m=0}^{+\infty }}}\left[{{{\frac{\Gamma \left( {{\gamma }_{1}}-{{\gamma }_{1}}k+m \right){{\left[ -\left| x-s+2na \right|{{\eta}^{-{{\beta }_{1}}}} \right]}^{k}}{{\left[ \delta {{\eta}^{\alpha }} \right]}^{m}}}{\Gamma \left( k+1 \right)\Gamma \left( {{\gamma }_{1}}-{{\gamma }_{1}}k \right)\Gamma \left( \alpha m+{{\beta }_{1}}-{{\beta }_{1}}k \right)\Gamma \left( m+1 \right)}}-}}\right.$$
$$\left.-{{{\frac{\Gamma \left( {{\gamma }_{1}}-{{\gamma }_{1}}k+m \right){{\left[ -\left| x+s+2na \right|{{\eta}^{-{{\beta }_{1}}}} \right]}^{k}}{{\left[ \delta {{\eta}^{\alpha }} \right]}^{m}}}{\Gamma \left( k+1 \right)\Gamma \left( {{\gamma }_{1}}-{{\gamma }_{1}}k \right)\Gamma \left( \alpha m+{{\beta }_{1}}-{{\beta }_{1}}k \right)\Gamma \left( m+1 \right)}}}}\right]=$$
$$=\frac{{{\eta}^{{{\beta }_{1}}-1}}}{2}\sum\limits_{n=-\infty }^{\infty }{\left[ {{E}_{12}}\left( \left. \begin{matrix}
   -{{\gamma }_{1}},1,{{\gamma }_{1}};\,\,\,\,\,\,\,\,\,\,\,\,\,\,\,\,\,\,\,\,\,\,\,\,\,\,\,  \\
   -{{\beta }_{1}},\alpha ,{{\beta }_{1}};-{{\gamma }_{1}},{{\gamma }_{1}};1,1;1,1  \\
\end{matrix} \right|\begin{matrix}
   -\left| x-s+2an \right|{{\eta}^{-{{\beta }_{1}}}}  \\
   \delta {{\eta}^{\alpha }}  \\
\end{matrix} \right) \right.}-$$
$$\left. -{{E}_{12}}\left( \left. \begin{matrix}
   -{{\gamma }_{1}},1,{{\gamma }_{1}};\,\,\,\,\,\,\,\,\,\,\,\,\,\,\,\,\,\,\,\,\,\,\,\,\,\,\,  \\
   -{{\beta }_{1}},\alpha ,{{\beta }_{1}};-{{\gamma }_{1}},{{\gamma }_{1}};1,1;1,1  \\
\end{matrix} \right|\begin{matrix}
   -\left| x+s+2an \right|{{\eta}^{-{{\beta }_{1}}}}  \\
   \delta {{\eta}^{\alpha }}  \\
\end{matrix} \right) \right].$$
After replacing $\eta$ with $t-\eta$ in the last integral of \eqref{eq2.30}, the kernel of this integral takes the form \eqref{eq2.5}.
Here the conditions for the convergence of the series are satisfied, that are ${{\Delta }_{1}}=-\beta_1-\gamma_1+1+\gamma_1=1-{{\beta }_{1}}>0,$ ${{\Delta }_{2}}=\alpha+1-1=\alpha >0.$

Taking expression \eqref{eq2.5} and the following relations 
$$\underset{s\to 0}{\mathop{\lim }}\,{{G}_{s}}\left( t,x,s,\eta  \right)=\sum\limits_{n=-\infty }^{\infty }{sign\left( x+2na \right)\omega \left( t-\eta ,\left| x+2na \right| \right)\,}$$
and
$$\underset{s\to a}{\mathop{\lim }}\,{{G}_{s}}\left( t,x,s,\eta  \right)=\sum\limits_{n=-\infty }^{+\infty }{sign\left( x+\left( 2n+1 \right)a \right)\omega \left( t-\eta ,\left| x+\left( 2n+1 \right)a \right| \right)\,}$$
into account, the equality \eqref{eq2.30} can be written in the form of \eqref{eq2.4}. From \eqref{eq2.4} it follows that the boundary value problem \eqref{eq2.2}-\eqref{eq2.3} for the equation \eqref{eq2.1} has a unique solution.

 \smallskip

\section{Appendices}

\subsection{\bf First appendix} \label{secA1}

First, we verify that  ${}^{P}I_{0t}^{\alpha,\,1-\beta,\,-\gamma,\,\delta}u(t,x)\in C(\overline{D}),$
that is, we check that  $t^{1-\beta}u(t,x)\in C(\overline{D}).$ Therefore, we multiply both sides of the solution by $t:$
$$t^{1-\beta}u\left( t,x \right)=t^{1-\beta}\left[\int\limits_{0}^{t}{{{\varphi }_{0}}\left( \eta  \right){{G}_{s}}\left( t,x,\eta ,0 \right)d\eta }-\int\limits_{0}^{t}{{{\varphi }_{1}}\left( \eta  \right){{G}_{s}}\left( t,x,\eta ,a \right)d\eta }+\right.$$
\begin{equation}\label{eq3}
   \left.+\,\int\limits_{0}^{a}{\tau \left( s \right)G\left( t,x,0,s \right)ds}+\int\limits_{0}^{t}{\int\limits_{0}^{a}{f\left( \eta ,s \right)G\left( t,x,\eta ,s \right)dsd\eta }}\right]. 
\end{equation}

Now, we consider the situation as $t \to 0$. Let us consider the first integral:
$$\lim_{t\to 0}\left[t^{1-\beta}\int\limits_{0}^{t}{\eta ^{1-\beta}{{\varphi }_{0}}\left( \eta  \right)\eta ^{\beta-1}{{G}_{s}}\left( t,x,\eta ,0 \right)d\eta }\right]=$$
$$=\lim_{t\to 0}\left[t^{1-\beta}\int\limits_{0}^{t}\eta ^{1-\beta}{{\varphi }_{0}}\left( \eta  \right)\eta ^{\beta-1}\sum\limits_{n=-\infty }^{+\infty }{sign\left( x+2na \right)\omega \left( t-\eta ,\left| x+2na \right| \right)\,}d\eta\right]=$$
$$=\lim_{t\to 0}\left[t^{1-\beta}\int\limits_{0}^{t}\eta ^{1-\beta}{{\varphi }_{0}}\left( \eta  \right)\eta ^{\beta-1}\sum\limits_{n=-\infty }^{+\infty }{sign\left( x+2na \right){{(t-\eta)}^{-1}}}\times\right.$$
$$\times \left.{{{E}_{12}}\left( \left. \begin{matrix}
   -{{\gamma }_{1}},1,0;\,\,\,\,\,\,\,\,\,\,\,\,\,\,\,\,\,\,\,\,\,\,\,\,\,\,\,  \\
   -{{\beta }_{1}},\alpha ,0;-{{\gamma }_{1}},0;1,1;1,1  \\
\end{matrix} \right|\begin{matrix}
   -\left| x+2na \right|{{(t-\eta)}^{-{{\beta }_{1}}}}  \\
   \delta {{(t-\eta)}^{\alpha }}  \\
\end{matrix} \right)}d\eta\right]=$$
$$=\lim _{t \to 0}\left[t^{1-\beta}\int\limits_{0}^{t}\eta ^{1-\beta}{{\varphi }_{0}}\left( \eta  \right)\eta ^{\beta-1}\sum\limits_{n=-\infty }^{+\infty }sign\left( x+2na \right)\times\right.$$
$$\left.\times\sum\limits_{k=0}^{+\infty }\frac{{{\left( -1 \right)}^{k}}{{\left| x+2na \right|}^{k}}}{k!}{{(t-\eta)}^{-{{\beta }_{1}}k-1}}E_{\alpha ,-{{\beta }_{1}}k}^{-{{\gamma }_{1}}k}\left[ \delta {{(t-\eta)}^{\alpha }} \right]d\eta\right].$$
Let us make the substitution $\eta=tz$. Then, by evaluating the limit of $E_{\alpha ,-{{\beta }_{1}}k}^{-{{\gamma }_{1}}k}\left[ \delta {{t^{\alpha}(1-z)}^{\alpha }} \right]$
as $t \to 0$, the remaining expression can be represented in terms of the Wright function $e_{\alpha,{\beta}}^{\mu,{\delta }}\left( z \right)$ \cite{Pskhu 2}:
$$\lim _{t \to 0}\left|\int\limits_{0}^{1}(tz)^{1-\beta}{{\varphi }_{0}}\left( tz \right)t^{1-\beta}(tz)^{\beta-1}\times\right.$$
$$\left.\times\sum\limits_{n=-\infty }^{+\infty }sign\left( x+2na \right)\sum\limits_{k=0}^{+\infty }\frac{{{\left( -1 \right)}^{k}}{{\left|x+2na \right|}^{k}}}{k!}{{t^{{-{{\beta }_{1}}k}}(1-z)}^{-{{\beta }_{1}}k-1}}E_{\alpha ,-{{\beta }_{1}}k}^{-{{\gamma }_{1}}k}\left[ \delta {{t^{\alpha}(1-z)}^{\alpha }} \right]dz\right|\le$$
$$\le\left\|(tz)^{1-\beta}{{\varphi }_{0}}\left( tz \right)\right\|\lim _{t \to 0}\left|\int\limits_{0}^{1}\frac{z^{\beta-1}}{1-z}\sum\limits_{n=-\infty }^{+\infty }sign\left( x+2na \right)\right.\times$$
$$\left.\times\sum\limits_{k=0}^{+\infty }\frac{1}{\Gamma(k+1)\Gamma(-{{\beta }_{1}}k)}\left(-\frac{\left| x+2na \right|}{\left[t(1-z)\right]^{\beta_1}}\right)^{k}dz\right|=$$
$$=\left\|(tz)^{1-\beta}{{\varphi }_{0}}\left( tz \right)\right\|\lim _{t \to 0}\left|\int\limits_{0}^{1}\frac{z^{\beta-1}}{1-z}\sum\limits_{n=-\infty }^{+\infty }sign\left( x+2na \right)e^{1,0}_{1,\beta_1}\left(-\frac{\left|x+2na \right|}{\left[t(1-z)\right]^{\beta_1}}\right)dz\right|.$$
Considering that \cite{Pskhu 2}
$$\lim _{t \to 0}e^{1,0}_{1,\beta_1}\left(-\frac{\left| x+2na \right|}{\left[t(1-z)\right]^{\beta_1}}\right)=0,$$
we arrive at the following result:
\begin{equation}\label{eq4}
    \lim_{t\to 0}\left[t^{1-\beta}\int\limits_{0}^{t}{\eta ^{1-\beta}{{\varphi }_{0}}\left( \eta  \right)\eta ^{\beta-1}{{G}_{s}}\left( t,x,\eta ,0 \right)d\eta }\right]=0.
\end{equation}

By the same method, if we also analyze the integral involving $\varphi_1(\eta)$ at t=0, we obtain the following:
\begin{equation}\label{eq5}
    \lim_{t\to 0}\left[t^{1-\beta}\int\limits_{0}^{t}\eta^{1-\beta}{{{\varphi }_{1}}\left( \eta  \right)\eta^{\beta-1}{{G}_{s}}\left( t,x,\eta ,a \right)d\eta }\right]=0.
\end{equation}

Now we analyze the third integral of \eqref{eq3}.
Based on the structure of the Green’s function, we can write the Green's function as follows
$${G}\left( t,x,\eta ,s\right)=$$
\begin{equation}\label{eq6*}
    =\frac{{(t-\eta)^{{{\beta }_{1}}-1}}}{2}{ {{E}_{12}}\left( \left. \begin{matrix}
   -{{\gamma }_{1}},1,{{\gamma }_{1}};\,\,\,\,\,\,\,\,\,\,\,\,\,\,\,\,\,\,\,\,\,\,\,\,\,\,\,  \\
   -{{\beta }_{1}},\alpha ,{{\beta }_{1}};-{{\gamma }_{1}},{{\gamma }_{1}};1,1;1,1  \\
\end{matrix} \right|\begin{matrix}
   -\left| x-s \right|{(t-\eta)^{-{{\beta }_{1}}}}  \\
   \delta {(t-\eta)^{\alpha }}  \\
\end{matrix} \right) }+K(t,x,\eta,s),
\end{equation}
where
$$K(t,x,\eta,s)=\frac{{{(t-\eta)}^{{{\beta }_{1}}-1}}}{2}\sum\limits_{n=-\infty }^{\infty }{}' { {{E}_{12}}\left( \left. \begin{matrix}
   -{{\gamma }_{1}},1,{{\gamma }_{1}};\,\,\,\,\,\,\,\,\,\,\,\,\,\,\,\,\,\,\,\,\,\,\,\,\,\,\,  \\
   -{{\beta }_{1}},\alpha ,{{\beta }_{1}};-{{\gamma }_{1}},{{\gamma }_{1}};1,1;1,1  \\
\end{matrix} \right|\begin{matrix}
   -\left| x-s+2an \right|{{(t-\eta)}^{-{{\beta }_{1}}}}  \\
   \delta {{(t-\eta)}^{\alpha }}  \\
\end{matrix} \right) }-$$
$$-\frac{{{(t-\eta)}^{{{\beta }_{1}}-1}}}{2}\sum\limits_{n=-\infty }^{\infty }{E}_{12}\left( \left. \begin{matrix}
   -{{\gamma }_{1}},1,{{\gamma }_{1}};\,\,\,\,\,\,\,\,\,\,\,\,\,\,\,\,\,\,\,\,\,\,\,\,\,\,\,  \\
   -{{\beta }_{1}},\alpha ,{{\beta }_{1}};-{{\gamma }_{1}},{{\gamma }_{1}};1,1;1,1  \\
\end{matrix} \right|\begin{matrix}
   -\left| x+s+2an \right|{{(t-\eta)}^{-{{\beta }_{1}}}}  \\
   \delta {{(t-\eta)}^{\alpha }}  \\
\end{matrix} \right).$$
The symbol $\sum\limits_{n=-\infty }^{\infty }{}'$ indicates that the summation is taken over all nonzero integer values of $n$.

Hence, the third integral of \eqref{eq3} can be expressed as 
$$I(t,x)=t^{1-\beta}\int\limits_{0}^{a}\tau \left(s \right)G\left( t,x,0,s \right)ds=$$
$$=\frac{1}{2}\int\limits_{0}^{a}\tau \left( s \right){t^{{{-\beta }_{1}}}}{ {{E}_{12}}\left( \left. \begin{matrix}
   -{{\gamma }_{1}},1,{{\gamma }_{1}};\,\,\,\,\,\,\,\,\,\,\,\,\,\,\,\,\,\,\,\,\,\,\,\,\,\,\,  \\
   -{{\beta }_{1}},\alpha ,{{\beta }_{1}};-{{\gamma }_{1}},{{\gamma }_{1}};1,1;1,1  \\
\end{matrix} \right|\begin{matrix}
   -\left| x-s \right|{t^{-{{\beta }_{1}}}}  \\
   \delta {t^{\alpha }}  \\
\end{matrix} \right) }ds+$$
$$+\int\limits_{0}^{a}\tau \left( s \right)t^{1-\beta}K(t,x,0,s)ds.$$

Now, we examine the first integral of $I(t,x)$ as $t\to 0.$ For this purpose, we first remove the modulus:
$$\frac{1}{2}\lim_{t\to 0}\int\limits_{0}^{a}\tau \left( s \right){t^{{{-\beta }_{1}}}}{ {{E}_{12}}\left( \left. \begin{matrix}
   -{{\gamma }_{1}},1,{{\gamma }_{1}};\,\,\,\,\,\,\,\,\,\,\,\,\,\,\,\,\,\,\,\,\,\,\,\,\,\,\,  \\
   -{{\beta }_{1}},\alpha ,{{\beta }_{1}};-{{\gamma }_{1}},{{\gamma }_{1}};1,1;1,1  \\
\end{matrix} \right|\begin{matrix}
   -\left| x-s \right|{t^{-{{\beta }_{1}}}}  \\
   \delta {t^{\alpha }}  \\
\end{matrix} \right) }ds=$$
$$=\frac{1}{2}\lim_{t\to 0}\int\limits_{0}^{x}\tau \left( s\right){t^{{{-\beta }_{1}}}}{ {{E}_{12}}\left( \left. \begin{matrix}
   -{{\gamma }_{1}},1,{{\gamma }_{1}};\,\,\,\,\,\,\,\,\,\,\,\,\,\,\,\,\,\,\,\,\,\,\,\,\,\,\,  \\
   -{{\beta }_{1}},\alpha ,{{\beta }_{1}};-{{\gamma }_{1}},{{\gamma }_{1}};1,1;1,1  \\
\end{matrix} \right|\begin{matrix}
   -(x-s){t^{-{{\beta }_{1}}}}  \\
   \delta {t^{\alpha }}  \\
\end{matrix} \right) }ds+$$
$$+\frac{1}{2}\lim_{t\to 0}\int\limits_{x}^{a}\tau \left(s \right){t^{{{-\beta }_{1}}}}{ {{E}_{12}}\left( \left. \begin{matrix}
   -{{\gamma }_{1}},1,{{\gamma }_{1}};\,\,\,\,\,\,\,\,\,\,\,\,\,\,\,\,\,\,\,\,\,\,\,\,\,\,\,  \\
   -{{\beta }_{1}},\alpha ,{{\beta }_{1}};-{{\gamma }_{1}},{{\gamma }_{1}};1,1;1,1  \\
\end{matrix} \right|\begin{matrix}
   -(s-x){t^{-{{\beta }_{1}}}}  \\
   \delta {t^{\alpha }}  \\
\end{matrix} \right) }ds=$$
$$=\frac{1}{2}\lim_{t\to 0}\left[{t^{{{-\beta }_{1}}}}\sum_{k=0}^{+\infty}\frac{1}{k!} E_{\alpha, -\beta_1k+\beta_1}^{-\gamma_1k+\gamma_1}[\delta t^\alpha]\int\limits_{0}^{x}\tau \left( s \right)\left(-\frac{x-s}{t^{\beta_1}}\right)^kds\right]+$$
$$+\frac{1}{2}\lim_{t\to 0}\left[{t^{{{-\beta }_{1}}}}\sum_{k=0}^{+\infty}\frac{1}{k!} E_{\alpha, -\beta_1k+\beta_1}^{-\gamma_1k+\gamma_1}[\delta t^\alpha]\int\limits_{x}^{a}\tau \left( s \right)\left(-\frac{s-x}{t^{\beta_1}}\right)^kds\right].$$
By substituting $z_1=-\frac{x-s}{t^{\beta_1}}$ in the first integral and 
$z_2=-\frac{s-x}{t^{\beta_1}}$
 in the second integral of the last equality, then we calculate the limit of both integrals as $t$ tends to $0$:
 $$\frac{1}{2}\lim_{t\to 0}\left[\sum_{k=0}^{+\infty}\frac{1}{k!} E_{\alpha, -\beta_1k+\beta_1}^{-\gamma_1k+\gamma_1}[\delta t^\alpha]\int\limits_{-\frac{x}{t^{\beta_1}}}^{0}\tau \left( x+t^{\beta_1}z_1 \right){z_1}^kdz_1\right]+$$
$$+\frac{1}{2}\lim_{t\to 0}\left[\sum_{k=0}^{+\infty}\frac{1}{k!} E_{\alpha, -\beta_1k+\beta_1}^{-\gamma_1k+\gamma_1}[\delta t^\alpha]\int\limits_{\frac{x-a}{t^{\beta_1}}}^{0}\tau \left( x-t^{\beta_1}z_2 \right)z_2^kdz_2\right]=$$
$$=\frac{1}{2}\lim_{t\to 0}\left[\sum_{k=0}^{+\infty}\frac{\tau \left( x \right)}{k!} E_{\alpha, -\beta_1k+\beta_1}^{-\gamma_1k+\gamma_1}[\delta t^\alpha]\int\limits_{-\infty}^{0}{z_1}^kdz_1\right]+$$
$$+\frac{1}{2}\lim_{t\to 0}\left[\sum_{k=0}^{+\infty}\frac{\tau \left( x \right)}{k!} E_{\alpha, -\beta_1k+\beta_1}^{-\gamma_1k+\gamma_1}[\delta t^\alpha]\int\limits_{-\infty}^{0}z_2^kdz_2\right].$$
Now we use the Wright function  and its property \cite{Pskhu 2}
\begin{equation}\label{eq6}
    \lim_{\left|z\right| \to \infty}z e_{\alpha,{\beta}}^{\mu,{\delta }}\left( z \right)=-\frac{1}{\Gamma(\mu-\alpha)\Gamma(\delta+\beta)}:
\end{equation}
$$\frac{\tau \left( x \right)}{2}\left[\sum_{k=0}^{+\infty}\left.\frac{{z_1}^{k+1}}{(k+1)!\Gamma(-\beta_1k+\beta_1)}\right|_{-\infty}^{0}+\sum_{k=0}^{+\infty}\left.\frac{{z_2}^{k+1}}{(k+1)!\Gamma(-\beta_1k+\beta_1)}\right|_{-\infty}^{0}\right]=$$
$$=\frac{\tau \left( x \right)}{2}\left[\left.z_1e^{2,\beta_1}_{1,\beta_1}(z_1)\right|_{-\infty}^{0}+\left.z_2e^{2,\beta_1}_{1,\beta_1}(z_2)\right|_{-\infty}^{0}\right]=\frac{\tau \left( x \right)}{\Gamma(\beta)}.$$

Next, we consider the following limit:
$$\lim_{t\to 0}\left[t^{1-\beta}K(t,x,0,s)\right]=$$
$$=\lim_{t\to 0}\left[\frac{{{t}^{{-{\beta }_{1}}}}}{2}\sum\limits_{n=-\infty }^{\infty }{}' { {{E}_{12}}\left( \left. \begin{matrix}
   -{{\gamma }_{1}},1,{{\gamma }_{1}};\,\,\,\,\,\,\,\,\,\,\,\,\,\,\,\,\,\,\,\,\,\,\,\,\,\,\,  \\
   -{{\beta }_{1}},\alpha ,{{\beta }_{1}};-{{\gamma }_{1}},{{\gamma }_{1}};1,1;1,1  \\
\end{matrix} \right|\begin{matrix}
   -\left| x-s+2an \right|{{t}^{-{{\beta }_{1}}}}  \\
   \delta {{t}^{\alpha }}  \\
\end{matrix} \right) }-\right.$$
$$\left.-\frac{{{t}^{{-{\beta }_{1}}}}}{2}\sum\limits_{n=-\infty }^{\infty }{E}_{12}\left( \left. \begin{matrix}
   -{{\gamma }_{1}},1,{{\gamma }_{1}};\,\,\,\,\,\,\,\,\,\,\,\,\,\,\,\,\,\,\,\,\,\,\,\,\,\,\,  \\
   -{{\beta }_{1}},\alpha ,{{\beta }_{1}};-{{\gamma }_{1}},{{\gamma }_{1}};1,1;1,1  \\
\end{matrix} \right|\begin{matrix}
   -\left| x+s+2an \right|{{t}^{-{{\beta }_{1}}}}  \\
   \delta {{t}^{\alpha }}  \\
\end{matrix} \right)\right]=$$
$$=\lim_{t\to 0}\left[\frac{{{t}^{{-{\beta }_{1}}}}}{2}\sum\limits_{n=-\infty }^{+\infty }{}'{\sum\limits_{k=0}^{+\infty }{\frac{{{\left( -1 \right)}^{k}} {{\left| x-s+2na \right|}^{k}}}{k!}}}{{t}^{-{{\beta }_{1}}k}}E_{\alpha ,{{\beta }_{1}}-{{\beta }_{1}}k}^{{{\gamma }_{1}}-{{\gamma }_{1}}k}\left[ \delta {{t}^{\alpha }} \right]\right.-$$
$$-\left.\frac{{{t}^{{-{\beta }_{1}}}}}{2}\sum\limits_{n=-\infty }^{+\infty }{\sum\limits_{k=0}^{+\infty }{\frac{{{\left( -1 \right)}^{k}} {{\left| x+s+2na \right|}^{k}}}{k!}}}{{t}^{-{{\beta }_{1}}k}}E_{\alpha ,{{\beta }_{1}}-{{\beta }_{1}}k}^{{{\gamma }_{1}}-{{\gamma }_{1}}k}\left[ \delta {{t}^{\alpha }} \right]\right]=$$
$$=\lim_{t\to 0}\left[\frac{{{t}^{{-{\beta }_{1}}}}}{2}\sum\limits_{n=-\infty }^{+\infty }{}'{\sum\limits_{k=0}^{+\infty }{\frac{{{\left( -1 \right)}^{k}} {{\left| x-s+2na \right|}^{k}{{t}^{-{{\beta }_{1}}k}}}}{k!}}}-\right.$$
$$\left.-\frac{{{t}^{{-{\beta }_{1}}}}}{2}\sum\limits_{n=-\infty }^{+\infty }{\sum\limits_{k=0}^{+\infty }{\frac{{{\left( -1 \right)}^{k}} {{\left| x+s+2na \right|}^{k}{{t}^{-{{\beta }_{1}}k}}}}{k!}}}\right]\times$$
$$\times\lim_{t\to 0}\left[\frac{1}{\Gamma({{\beta }_{1}}-{{\beta }_{1}}k)}+\frac{(\gamma_1-\gamma_1k)_1}{\Gamma(\alpha +{{\beta }_{1}}-\beta_1k)}\delta t^\alpha+\frac{(\gamma_1-\gamma_1k)_2}{2\Gamma(2\alpha +{{\beta }_{1}}-\beta_1k)}\delta^2t^{2\alpha}+...\right].$$
By evaluating the last limit in the final equality, the remaining expression can be represented in terms of the Wright function:
$$\lim_{t\to 0}\left[t^{1-\beta}K(t,x,0,s)\right]=$$
$$=\lim_{t\to 0}\left[\frac{{{t}^{{-{\beta }_{1}}}}}{2}\sum\limits_{n=-\infty }^{+\infty }{}'\sum\limits_{k=0}^{+\infty }\frac{1}{\Gamma(k+1)\Gamma({{\beta }_{1}}-{{\beta }_{1}}k)}\left(-\frac{\left| x-s+2na \right|}{t^{\beta_1}}\right)^k-\right.$$
$$-\left.\frac{{{t}^{{-{\beta }_{1}}}}}{2}\sum\limits_{n=-\infty }^{+\infty }\sum\limits_{k=0}^{+\infty }\frac{1}{\Gamma(k+1)\Gamma({{\beta }_{1}}-{{\beta }_{1}}k)}\left(-\frac{\left| x+s+2na \right|}{t^{\beta_1}}\right)^k\right]=$$
$$=\frac{1}{2}\lim_{t\to 0}\left[\sum\limits_{n=-\infty }^{+\infty }{}'\left({{t}^{{-{\beta }_{1}}}}e^{1,\beta_1}_{1,\beta_1}\left(-\frac{\left| x-s+2na \right|}{t^{\beta_1}}\right)\right)-\sum\limits_{n=-\infty }^{+\infty }\left({{t}^{{-{\beta }_{1}}}}e^{1,\beta_1}_{1,\beta_1}\left(-\frac{\left| x+s+2na \right|}{t^{\beta_1}}\right)\right)\right].$$
We express the last equality in the following form and apply the formula \eqref{eq6}: 
$$\frac{1}{2}\lim_{t\to 0}\left[\sum\limits_{n=-\infty }^{+\infty }{}'\left(-\frac{1}{\left| x-s+2na \right|}\right)\left(-\frac{\left| x-s+2na \right|}{t^{\beta_1}}\right) e^{1,\beta_1}_{1,\beta_1}\left(-\frac{\left| x-s+2na \right|}{t^{\beta_1}}\right)-\right.$$
$$\left.-\sum\limits_{n=-\infty }^{+\infty }\left(-\frac{1}{\left| x+s+2na \right|}\right)\left(-\frac{\left| x+s+2na \right|}{t^{\beta_1}}\right)e^{1,\beta_1}_{1,\beta_1}\left(-\frac{\left| x+s+2na \right|}{t^{\beta_1}}\right)\right]=0.$$

Thus,  
$$\lim_{t\to 0}\left[t^{1-\beta}K(t,x,0,s)\right]=0.$$

We conclude that 
\begin{equation}\label{eq7}
    \lim_{t\to 0}I(t,x)=\frac{\tau(x)}{\Gamma(\beta)}.
\end{equation}

Next, we analyze the final integral of \eqref{eq3} as $t\to 0:$
$$\lim_{t\to 0}\left[t^{1-\beta}\int\limits_{0}^{t}{\int\limits_{0}^{a}{\eta ^{1-\beta}f\left( \eta ,s \right)\eta^{\beta-1}G\left( t,x,\eta ,s \right)dsd\eta }}\right]=$$
$$=\frac{1}{2}\lim_{t\to 0}\left[t^{1-\beta}\int\limits_{0}^{t}{\int\limits_{0}^{a}{\eta ^{1-\beta}f\left( \eta , s \right)\eta^{\beta-1}}}{{\left( t-\eta  \right)}^{{{\beta }_{1}}-1}}\times\right.$$
$$\left.\times\sum\limits_{n=-\infty }^{\infty }{\sum\limits_{k=0}^{+\infty }}\frac{ (-1)^k\left[\left| x-s+2na \right|^k-\left| x+s+2na \right|^k\right]}{k!}(t-\eta)^{-\beta_{1}{k}}E_{\alpha,\beta_1-\beta_1 k}^{{{\gamma }_{1}}-{{\gamma }_{1}}k}\left[\delta {{(t-\eta)}^{\alpha }} \right]ds d\eta\right].$$
We perform the substitution $\eta=tz.$ Following the same procedure as above, we find the limit of the Prabhakar function and write the remaining expression through the Wright function:
$$\frac{1}{2}\lim_{t\to 0}\left|t^{\beta_1}\int\limits_{0}^{1}{\int\limits_{0}^{a}{(tz) ^{1-\beta}f\left( tz ,s \right)z^{\beta-1}}}{{\left( 1-z  \right)}^{{{\beta }_{1}}-1}}\times\right.$$
$$\times\sum\limits_{n=-\infty }^{\infty }{\sum\limits_{k=0}^{+\infty }}\frac{ (-1)^k\left[\left| x-s+2na \right|^k-\left| x+s+2na \right|^k\right]}{k!}t^{-\beta_{1}{k}}(1-z)^{-\beta_{1}{k}}\times$$
$$\left.\times E_{\alpha,\beta_1-\beta_1 k}^{{{\gamma }_{1}}-{{\gamma }_{1}}k}\left[\delta {{t^\alpha(1-z)}^{\alpha }} \right]ds dz\right|\le$$
$$\le\frac{1}{2}\left\|(tz) ^{1-\beta}f\left( tz ,s \right)\right\|\lim_{t\to 0}\left|t^{\beta_1}\int\limits_{0}^{1}{\int\limits_{0}^{a}{z^{\beta-1}}}{{\left( 1-z  \right)}^{{{\beta }_{1}}-1}}\times\right.$$
$$\left.\times\sum\limits_{n=-\infty }^{\infty }\left(e^{1,\beta_1}_{1,\beta_1}\left(-\frac{\left| x-s+2na \right|}{\left[t(1-z)\right]^{\beta_1}}\right)-e^{1,\beta_1}_{1,\beta_1}\left(-\frac{\left| x+s+2na \right|}{\left[t(1-z)\right]^{\beta_1}}\right)\right)ds dz\right|.$$
Since
$$\lim_{t\to 0}e^{1,\beta_1}_{1,\beta_1}\left(-\frac{\left| x\pm s+2na \right|}{\left[t(1-z)\right]^{\beta_1}}\right)=0$$
and the exponent of $t$ is positive, it follows that this limit equals zero:
\begin{equation}\label{eq8}
    \lim_{t\to 0}\left[t^{1-\beta}\int\limits_{0}^{t}{\int\limits_{0}^{a}{\eta ^{1-\beta}f\left( \eta ,s \right)\eta^{\beta-1}G\left( t,x,\eta ,s \right)ds d\eta }}\right]=0.
\end{equation}

According to ${{t}^{1-\beta }}{{\varphi }_{0}}\left( t \right),$ ${{t}^{1-\beta }}{{\varphi }_{1}}\left( t \right)\in C\left[ 0;T \right],$ $\tau \left( x \right)\in C\left[ 0;a \right],$ ${{t}^{1-\beta }}f\left( t,x \right)\in C\left( \overline{D} \right)$ and the equalities \eqref{eq4}, \eqref{eq5}, \eqref{eq7}, \eqref{eq8}, the expression \eqref{eq3} is continuous at the point $t=0.$

Now, we examine the case at $\eta \to t$. Let us begin by examining the term in \eqref{eq3} that involves the function $\varphi_0:$ 

$$\lim_{\eta \to t}\left[{{\varphi }_{0}}\left( \eta  \right){{G}_{s}}\left( t,x,\eta ,0 \right)\right]=$$
$$=\lim _{\eta \to t}\left[{{\varphi }_{0}}\left( \eta  \right)\sum\limits_{n=-\infty }^{+\infty }sign\left( x+2na \right)\sum\limits_{k=0}^{+\infty }\frac{{{\left( -1 \right)}^{k}}{{\left| x+2na \right|}^{k}}}{k!}{{(t-\eta)}^{-{{\beta }_{1}}k-1}}E_{\alpha ,-{{\beta }_{1}}k}^{-{{\gamma }_{1}}k}\left[ \delta {{(t-\eta)}^{\alpha }}\right]\right]$$
$$=\lim _{\eta \to t}\left[{{\varphi }_{0}}\left( \eta  \right)\sum\limits_{n=-\infty }^{+\infty }sign\left( x+2na \right)\sum\limits_{k=0}^{+\infty }\frac{{{\left( -1 \right)}^{k}}{{\left| x+2na \right|}^{k}}}{k!}{{(t-\eta)}^{-{{\beta }_{1}}k-1}}\right]\times$$
$$\times\lim _{\eta \to t}\left[\frac{1}{\Gamma(-{{\beta }_{1}}k)}+\frac{(-\gamma_1k)_1}{\Gamma(\alpha -\beta_1k)}\delta(t-\eta)^\alpha+\frac{(-\gamma_1k)_2}{2\Gamma(2\alpha -\beta_1k)}\delta^2(t-\eta)^{2\alpha}+...\right]=$$
$$=\lim _{\eta \to t}\left[{{\varphi }_{0}}\left( \eta  \right)\sum\limits_{n=-\infty }^{+\infty }sign\left( x+2na \right)\sum\limits_{k=0}^{+\infty }\frac{{{\left( -1 \right)}^{k}}{{\left| x+2na \right|}^{k}}}{k!\Gamma(-{{\beta }_{1}}k)}{{(t-\eta)}^{-{{\beta }_{1}}k-1}}\right].$$

We rewrite the final expression using the Wright function and turn to the following lemma \cite{Pskhu 2}:

\textbf{Lemma 3.1.}
   \textit{If ${\delta }_{1} <1$ and for arbitrary positive $x$ and $y$, ${{x}_{0}}\in \left( 0;x \right),$ ${{\beta }_{0}}\in \left[ \beta ;1 \right],$ $\theta \in \left( \frac{1}{2},\min \left\{ 1,\frac{1}{2\beta } \right\} \right),$ the following inequality holds:
$$\left| {{y}^{{\delta }_{1} -1}}e_{1,\beta }^{1,{\delta }_{1} }\left( -\frac{x}{{{y}^{\beta }}} \right) \right|\le \frac{1}{\beta \pi }{{\left( \frac{\cos \beta \theta \pi }{\beta \pi } \right)}^{\frac{{\delta }_{1} -1}{\beta }}}\Gamma \left( \frac{1-{\delta }_{1} }{\beta } \right){{\left( x-{{x}_{0}} \right)}^{\frac{{\delta }_{1} -1}{\beta }}}e_{1,\beta }^{1,1}\left( -\frac{{{x}_{0}}}{{{y}^{\beta }}} \right).$$}

According to Lemma 3.1, we find that its value is equal to zero:
$$\lim _{\eta \to t}\left[{{\varphi }_{0}}\left( \eta  \right)\sum\limits_{n=-\infty }^{+\infty }sign\left( x+2na \right){{(t-\eta)}^{-1}}e_{1,{\beta_1}}^{1,{0}}\left(-\frac{{\left| x+2na \right|}}{{{(t-\eta)}^{{{\beta }_{1}}}}}\right)\right]=0.$$

Hence, 
\begin{equation}\label{eq9}
    \lim_{\eta \to t}\left[{{\varphi }_{0}}\left( \eta  \right){{G}_{s}}\left( t,x,\eta ,0 \right)\right]=0.
\end{equation}

By the same method,
\begin{equation}\label{eq10}
    \lim_{\eta \to t}\left[{{\varphi }_{1}}\left( \eta  \right){{G}_{s}}\left( t,x,\eta ,a\right)\right]=0.
\end{equation}

Next, we examine the final integrand of \eqref{eq3} the case $\eta \to t:$

$$\lim_{\eta \to t}\left[f\left( \eta ,s \right)G\left( t,x,\eta ,s \right)\right]=\lim_{\eta \to t}\left[f\left( \eta ,s \right){{\left( t-\eta  \right)}^{{{\beta }_{1}}-1}}\times\right.$$
$$\left.\times\sum\limits_{n=-\infty }^{\infty }{\sum\limits_{k=0}^{+\infty }}\frac{ (-1)^k\left[\left| x-s+2na \right|^k-\left| x+s+2na \right|^k\right]}{k!}(t-\eta)^{-\beta_{1}{k}}E_{\alpha,\beta_1-\beta_1 k}^{{{\gamma }_{1}}-{{\gamma }_{1}}k}\left[\delta {{(t-\eta)}^{\alpha }} \right]\right].$$
Using the equality 
$$\lim_{\eta \to t}E_{\alpha,\beta_1-\beta_1 k}^{{{\gamma }_{1}}-{{\gamma }_{1}}k}\left[\delta {{(t-\eta)}^{\alpha }} \right]=\frac{1}{\Gamma(\beta_1-\beta_1 k)},$$ 
the sum over $k$ can be represented again in terms of the Wright function. By Lemma 3.1, it is not difficult to verify that this limit is equal to zero:
$$\frac{1}{2}\lim_{\eta \to t}\left[f\left( \eta ,\xi \right){{\left( t-\eta  \right)}^{{{\beta }_{1}}-1}}\sum\limits_{n=-\infty }^{\infty }\left[e^{1,\beta_1}_{1,\beta_1}\left(-\frac{\left| \xi_0-\xi+2na \right|}{(t-\eta)^{\beta_{1}}}\right)-e^{1,\beta_1}_{1,\beta_1}\left(-\frac{\left| \xi_0+\xi+2na \right|}{(t-\eta)^{\beta_{1}}}\right)\right]\right]=0.$$

Thus, 
\begin{equation}\label{eq11}
    \lim_{\eta \to t}\left[f\left( \eta ,s \right)G\left( t,x,\eta ,s \right)\right]=0.
\end{equation}

From the equalities \eqref{eq9}, \eqref{eq10}, and \eqref{eq11}, it follows that expression \eqref{eq3} is also continuous at the point $\eta=t$.

Now we verify that ${{u}_{xx}}\left( t,x \right)\in C\left( D \right).$ To do this, we first compute the second derivative of the solution with respect to $x$:
$$u_{xx}\left( t,x \right)=\int\limits_{0}^{t}{{{\varphi }_{0}}\left( \eta  \right)\frac{\partial^2}{\partial x^2}{{G}_{s}}\left( t,x,\eta ,0 \right)d\eta }-\int\limits_{0}^{t}{{{\varphi }_{1}}\left( \eta  \right)\frac{\partial^2}{\partial x^2}{{G}_{s}}\left( t,x,\eta ,a \right)d\eta }+$$
\begin{equation}\label{eq12}
    +\,\int\limits_{0}^{a}{\tau \left( s \right)\frac{\partial^2}{\partial x^2}G\left( t,x,0,s \right)ds}+\int\limits_{0}^{t}{\int\limits_{0}^{a}{f\left( \eta ,s \right)\frac{\partial^2}{\partial x^2}G\left( t,x,\eta ,s \right)dsd\eta }}.
\end{equation}
We first compute $\frac{\partial^2}{\partial x^2}{{G}_{s}}\left( t,x,\eta ,0 \right)$ and $\frac{\partial^2}{\partial x^2}{G}\left( t,x,\eta ,s \right):$

$$1)\,\,\,\frac{\partial^2}{\partial x^2}{{G}_{s}}\left( t,x,\eta ,0 \right)=\frac{\partial^2}{\partial x^2}\sum_{n=-\infty}^{+\infty}sign(x+2na)\omega(t-\eta,\left|x+2na\right|)=$$
$$=\frac{\partial^2}{\partial x^2}\left[-\sum_{n=-\infty}^{-1}\omega(t-\eta,-x-2na)+\sum_{n=0}^{+\infty}\omega(t-\eta,x+2na)\right]=$$
$$=\frac{\partial^2}{\partial x^2}\sum_{n=-\infty}^{-1}\left(-(t-\eta)^{-1}E_{\alpha, 0}^{0}\left[ \delta {{(t-\eta)}^{\alpha }} \right]-\sum\limits_{k=1}^{+\infty }{\frac{{}{{(x+2na)}^{k}}}{k!}{{(t-\eta)}^{-{{\beta }_{1}}k-1}}E_{\alpha ,-{{\beta }_{1}}k}^{-{{\gamma }_{1}}k}\left[ \delta {{(t-\eta)}^{\alpha }} \right]}\right)+$$
$$+\frac{\partial^2}{\partial x^2}\sum_{n=0}^{+\infty}\left((t-\eta)^{-1}E_{\alpha, 0}^{0}\left[ \delta {{(t-\eta)}^{\alpha }} \right]+\sum\limits_{k=1}^{+\infty }{\frac{{{\left( -1 \right)}^{k}}{{(x+2na)}^{k}}}{k!}{{(t-\eta)}^{-{{\beta }_{1}}k-1}}E_{\alpha ,-{{\beta }_{1}}k}^{-{{\gamma }_{1}}k}\left[ \delta {{(t-\eta)}^{\alpha }} \right]}\right)=$$
$$=\frac{\partial}{\partial x}\sum_{n=-\infty}^{-1}\left(-\sum\limits_{k=0}^{+\infty }{\frac{{}{{(x+2na)}^{k}}}{k!}{{(t-\eta)}^{-{{\beta }_{1}}k-\beta_1-1}}E_{\alpha ,-{{\beta }_{1}}k-\beta_1}^{-{{\gamma }_{1}}k-\gamma_1}\left[ \delta {{(t-\eta)}^{\alpha }} \right]}\right)+$$
$$+\frac{\partial}{\partial x}\sum_{n=0}^{+\infty}\left(-\sum\limits_{k=0}^{+\infty }{\frac{{{\left( -1 \right)}^{k}}{{(x+2na)}^{k}}}{k!}{{(t-\eta)}^{-{{\beta }_{1}}k-\beta_1-1}}E_{\alpha ,-{{\beta }_{1}}k-\beta_1}^{-{{\gamma }_{1}}k-\gamma_1}\left[ \delta {{(t-\eta)}^{\alpha }} \right]}\right)=$$
$$=\frac{\partial}{\partial x}\sum_{n=-\infty}^{-1}\left(-{{(t-\eta)}^{-\beta_1-1}}E_{\alpha ,-\beta_1}^{-\gamma_1}\left[ \delta {{(t-\eta)}^{\alpha }} \right]-\right.$$
$$\left.-\sum\limits_{k=1}^{+\infty }{\frac{{}{{(x+2na)}^{k}}}{k!}{{(t-\eta)}^{-{{\beta }_{1}}k-\beta_1-1}}E_{\alpha ,-{{\beta }_{1}}k-\beta_1}^{-{{\gamma }_{1}}k-\gamma_1}\left[ \delta {{(t-\eta)}^{\alpha }} \right]}\right)+$$
$$+\frac{\partial}{\partial x}\sum_{n=0}^{+\infty}\left(-{{(t-\eta)}^{-\beta_1-1}}E_{\alpha ,-\beta_1}^{-\gamma_1}\left[ \delta {{(t-\eta)}^{\alpha }} \right]-\right.$$
$$\left.-\sum\limits_{k=1}^{+\infty }{\frac{{{\left( -1 \right)}^{k}}{{(x+2na)}^{k}}}{k!}{{(t-\eta)}^{-{{\beta }_{1}}k-\beta_1-1}}E_{\alpha ,-{{\beta }_{1}}k-\beta_1}^{-{{\gamma }_{1}}k-\gamma_1}\left[ \delta {{(t-\eta)}^{\alpha }} \right]}\right)=$$
$$=-\sum_{n=-\infty}^{-1}\sum\limits_{k=0}^{+\infty }{\frac{{}{{(x+2na)}^{k}}}{k!}{{(t-\eta)}^{-{{\beta }_{1}}k-2\beta_1-1}}E_{\alpha ,-{{\beta }_{1}}k-2\beta_1}^{-{{\gamma }_{1}}k-2\gamma_1}\left[ \delta {{(t-\eta)}^{\alpha }} \right]}+$$
$$+\sum_{n=0}^{+\infty}\sum\limits_{k=0}^{+\infty }{\frac{{{\left( -1 \right)}^{k}}{{(x+2na)}^{k}}}{k!}{{(t-\eta)}^{-{{\beta }_{1}}k-2\beta_1-1}}E_{\alpha ,-{{\beta }_{1}}k-2\beta_1}^{-{{\gamma }_{1}}k-2\gamma_1}\left[ \delta {{(t-\eta)}^{\alpha }} \right]}=$$
$$=\sum_{n=-\infty}^{+\infty}sign(x+2na)\sum_{k=0}^{+\infty}\frac{(-1)^k\left|x+2na\right|^k(t-\eta)^{-\beta_1k-\beta-1}}{k!}E_{\alpha, -\beta_1k-\beta}^{-\gamma-\gamma_1k}\left[\delta(t-\eta)^{\alpha }\right];$$

$$2)\,\,\,\frac{\partial^2}{\partial x^2}G\left( t,x,\eta,s \right)=$$
$$=\frac{\partial^2}{\partial x^2}\left(\frac{{{ (t-\eta) }^{{{\beta }_{1}-1}}}}{2}\sum\limits_{n=-\infty }^{\infty }{\left[ {{E}_{12}}\left( \left. \begin{matrix}
-{{\gamma }_{1}},1,{{\gamma }_{1}};\,\,\,\,\,\,\,\,\,\,\,\,\,\,\,\,\,\,\,\,\,\,\,\,\,\,\,  \\
   -{{\beta }_{1}},\alpha ,{{\beta }_{1}};-{{\gamma }_{1}},{{\gamma }_{1}};1,1;1,1  \\
\end{matrix} \right|\begin{matrix}
   -\left| x-s+2an \right|{{(t-\eta)}^{-{{\beta }_{1}}}}  \\
   \delta {{(t-\eta)}^{\alpha }}  \\
\end{matrix} \right) \right.}-\right.$$
$$\left.\left. -{{E}_{12}}\left( \left. \begin{matrix}
   -{{\gamma }_{1}},1,{{\gamma }_{1}};\,\,\,\,\,\,\,\,\,\,\,\,\,\,\,\,\,\,\,\,\,\,\,\,\,\,\,  \\
   -{{\beta }_{1}},\alpha ,{{\beta }_{1}};-{{\gamma }_{1}},{{\gamma }_{1}};1,1;1,1  \\
\end{matrix} \right|\begin{matrix}
   -\left| x+s+2an \right|{{(t-\eta)}^{-{{\beta }_{1}}}}  \\
   \delta {{(t-\eta)}^{\alpha }}  \\
\end{matrix} \right) \right]\right)=$$
$$=\frac{\partial^2}{\partial x^2}\left( \sum\limits_{n=-\infty }^{+\infty }{\sum\limits_{k=0}^{+\infty }{\frac{{{\left( -1 \right)}^{k}}\left[ {{\left| x-s+2na \right|}^{k}}-{{\left| x+s+2na \right|}^{k}} \right]}{2\Gamma(k+1)}}}(t-\eta)^{-\beta_1k+{{\beta }_{1}}-1}E_{\alpha, \beta_1-\beta_1k}^{{\gamma }_{1}-{{\gamma }_{1}}k}\left[\delta (t-\eta)^\alpha\right]\right)=$$
$$=\frac{\partial^2}{\partial x^2}\left( \sum\limits_{n=-\infty }^{-1 }{\sum\limits_{k=0}^{+\infty }{\frac{\left[ {{( x-s+2na)}^{k}}-{{(x+s+2na)}^{k}} \right]}{2\Gamma(k+1)}}}(t-\eta)^{-\beta_1k+{{\beta }_{1}-1}}E_{\alpha, \beta_1-\beta_1k}^{{\gamma }_{1}-{{\gamma }_{1}}k}\left[\delta (t-\eta)^\alpha\right]\right)+$$
$$+\frac{\partial^2}{\partial x^2}\left( {\sum\limits_{k=0}^{+\infty }{\frac{{(-1)^k( x-s)}^{k}}{2\Gamma(k+1)}}}(t-\eta)^{-\beta_1k+{{\beta }_{1}}-1}E_{\alpha, \beta_1-\beta_1k}^{{\gamma }_{1}-{{\gamma }_{1}}k}\left[\delta (t-\eta)^\alpha\right]+\right.$$
$$+{\sum\limits_{k=0}^{+\infty }{\frac{{(-1)^k( s-x)}^{k}}{2\Gamma(k+1)}}}(t-\eta)^{-\beta_1k+{{\beta }_{1}}-1}E_{\alpha, \beta_1-\beta_1k}^{{\gamma }_{1}-{{\gamma }_{1}}k}\left[\delta (t-\eta)^\alpha\right]-$$
$$-\left.{\sum\limits_{k=0}^{+\infty }{\frac{{(-1)^k( x+s)}^{k}}{2\Gamma(k+1)}}}(t-\eta)^{-\beta_1k+{{\beta }_{1}}-1}E_{\alpha, \beta_1-\beta_1k}^{{\gamma }_{1}-{{\gamma }_{1}}k}\left[\delta (t-\eta)^\alpha\right] \right)+$$
$$+\frac{\partial^2}{\partial x^2}\left( \sum\limits_{n=1 }^{+\infty }{\sum\limits_{k=0}^{+\infty }{\frac{(-1)^k\left[ {{( x-s+2na)}^{k}}-{{(x+s+2na)}^{k}} \right]}{2\Gamma(k+1)}}}(t-\eta)^{-\beta_1k+{{\beta }_{1}}-1}E_{\alpha, \beta_1-\beta_1k}^{{\gamma }_{1}-{{\gamma }_{1}}k}\left[\delta (t-\eta)^\alpha\right]\right)=$$
$$=\sum\limits_{n=-\infty }^{+\infty }{\sum\limits_{k=0}^{+\infty }{\frac{{{\left( -1 \right)}^{k}}\left[ {{\left| x-s+2na \right|}^{k}}-{{\left| x+s+2na \right|}^{k}} \right]}{2\Gamma(k+1)}}}(t-\eta)^{-\beta_{1}k-\beta_1-1}E_{\alpha, {-\beta_1-\beta_1k}}^{-\gamma _{1}-{{\gamma }_{1}}k}\left[\delta (t-\eta)^\alpha\right].$$

Next, we evaluate the limits of $\frac{\partial^2}{\partial x^2}{{G}_{s}}\left( t,x,\eta ,0 \right),$ $\frac{\partial^2}{\partial x^2}{{G}_{s}}\left( t,x,\eta ,a \right),$ and $\frac{\partial^2}{\partial x^2}{G}\left( t,x,\eta ,s \right)$ as $\eta \to t:$

$$\lim_{\eta \to t} \frac{\partial^2}{\partial x^2}{{G}_{s}}\left( t,x,\eta ,0 \right)=$$
$$=\lim_{\eta \to t} \sum_{n=-\infty}^{+\infty}sign(x+2na)\sum_{k=0}^{+\infty}\frac{(-1)^k\left|x+2na\right|^k(t-\eta)^{-\beta_1k-\beta-1}}{k!}E_{\alpha, -\beta_1k-\beta}^{-\gamma-\gamma_1k}\left[\delta(t-\eta)^{\alpha }\right].$$
Considering that 
$$\lim_{\eta \to t}E_{\alpha, -\beta_1k-\beta}^{-\gamma-\gamma_1k}\left[\delta(t-\eta)^{\alpha }\right]=\frac{1}{\Gamma(-\beta_1k-\beta)},$$
the remaining expression can be written in terms of the Wright function:
$$\lim_{\eta \to t} \sum_{n=-\infty}^{+\infty}sign(x+2na)t^{-\beta-1}e_{1,\beta_1}^{1,-\beta}\left(-\frac{\left|x+2na\right|}{(t-\eta)^{\beta_1}}\right).$$
It is easy to see that this also equals zero according to the Lemma 3.1:
$$\lim_{\eta \to t} \frac{\partial^2}{\partial x^2}{{G}_{s}}\left( t,x,\eta ,0 \right)=0.$$
As above, it follows that 
$$\lim_{\eta \to t}\frac{\partial^2}{\partial x^2}{{G}_{s}}\left( t,x,\eta ,a \right)=0$$
and
$$\lim_{\eta \to t}\frac{\partial^2}{\partial x^2}{G}\left( t,x,\eta ,s \right)=0.$$

Hence, $\frac{\partial^2}{\partial x^2}{{G}_{s}}\left( t,x,\eta ,0 \right),$ $\frac{\partial^2}{\partial x^2}{{G}_{s}}\left( t,x,\eta ,a \right),$ and $\frac{\partial^2}{\partial x^2}{G}\left( t,x,\eta ,s \right)$ are continuous $\eta \to t.$ 
There are no singularities at the remaining points. Therefore, ${{u}_{xx}}\left( t,x \right)\in C\left( D \right).$

Let us confirm that ${}^{PRL}D_{0t}^{\alpha ,\,\beta ,\,\gamma ,\,\delta }u\left( t,x \right)\in C\left( D \right).$ For this, we apply the Prabhakar fractional derivative to both sides of the solution $u(t,x)$
by $t:$
$${}^{PRL}D_{0t}^{\alpha ,\,\beta ,\,\gamma ,\,\delta }u\left( t,x \right)={}^{PRL}D_{0t}^{\alpha ,\,\beta ,\,\gamma ,\,\delta }\left[\int\limits_{0}^{t}{{{\varphi }_{0}}\left( \eta  \right){{G}_{s}}\left( t,x,\eta ,0 \right)d\eta }-\int\limits_{0}^{t}{{{\varphi }_{1}}\left( \eta  \right){{G}_{s}}\left( t,x,\eta ,a \right)d\eta }+\right.$$
\begin{equation}\label{eq13}
   \left.+\,\int\limits_{0}^{a}{\tau \left( s \right)G\left( t,x,0,s \right)ds}+\int\limits_{0}^{t}{\int\limits_{0}^{a}{f\left( \eta ,s \right)G\left( t,x,\eta ,s \right)dsd\eta }}\right]. 
\end{equation}
Let us begin by calculating the term involving the function $\varphi_0(\eta):$
$${}^{PRL}D_{0t}^{\alpha ,\,\beta ,\,\gamma ,\,\delta }\left(\int\limits_{0}^{t}{{{\varphi }_{0}}\left( \eta  \right){{G}_{s}}\left( t,x,\eta ,0 \right)d\eta }\right)={\frac{\partial^{}}{\partial t}}{}^{P}I_{0 t}^{\alpha, 1-\beta,-\gamma, \delta} \left(\int\limits_{0}^{t}{{{\varphi }_{0}}\left( \eta  \right){{G}_{s}}\left( t,x,\eta ,0 \right)d\eta }\right)=$$
$$=\frac{\partial}{\partial t}\int\limits_0^t(t-y)^{-\beta} E_{\alpha, 1-\beta}^{-\gamma}\left[\mathcal{\delta}(t-y)^\alpha\right]dy\int\limits_{0}^{y}{{{\varphi }_{0}}\left( \eta  \right){{G}_{s}}\left( y,x,\eta ,0 \right)d\eta }=$$
$$=\frac{\partial}{\partial t}\int\limits_0^t{{\varphi }_{0}}\left( \eta  \right)d\eta\int\limits_{\eta}^{t}(t-y)^{-\beta} E_{\alpha, 1-\beta}^{-\gamma}\left[\mathcal{\delta}(t-y)^\alpha\right]{{{G}_{s}}\left( y,x,\eta ,0 \right)dy. }$$
Now we compute the integral with respect to $y$ separately:
$$\int\limits_{\eta}^{t}(t-y)^{-\beta} E_{\alpha, 1-\beta}^{-\gamma}\left[\mathcal{\delta}(t-y)^\alpha\right]{{G}_{s}}\left( y,x,\eta ,0 \right)dy=\sum_{n=-\infty}^{+\infty}sign(x+2na)\sum_{k=0}^{+\infty}\frac{(-1)^k\left|x+2na\right|^k}{k!}\times$$
$$\times\int\limits_{\eta}^{t}(t-y)^{-\beta} (y-\eta)^{-\beta_1k-1}E_{\alpha, 1-\beta}^{-\gamma}\left[\mathcal{\delta}(t-y)^\alpha\right]E_{\alpha, -\beta_1k}^{-\gamma_1k}\left[\mathcal{\delta}(y-\eta)^\alpha\right]dy.$$
Next, we use \eqref{eq2.27} for the generalized Mittag-Leffler functions and perform the substitution $y=(t-\eta)z+\eta$:
$$\sum_{n=-\infty}^{+\infty}sign(x+2na)\sum_{k=0}^{+\infty}\frac{(-1)^k\left|x+2na\right|^k}{k!}\sum_{i=0}^{+\infty}\sum_{j=0}^{i}\frac{\delta^i(-\gamma)_j(-\gamma_1k)_{i-j}}{j!(i-j)!\Gamma(\alpha j+1-\beta)\Gamma(\alpha i-\alpha j-\beta_1k)}\times$$
$$\times\int\limits_{\eta}^{t}(t-y)^{\alpha j-\beta} (y-\eta)^{\alpha i-\alpha j-\beta_1k-1}dy=$$
$$=\sum_{n=-\infty}^{+\infty}sign(x+2na)\sum_{k=0}^{+\infty}\frac{(-1)^k\left|x+2na\right|^k}{k!}\sum_{i=0}^{+\infty}\frac{\delta^i (t-\eta)^{\alpha i-\beta_1k-\beta}}{\Gamma(\alpha i-\beta_1k-\beta+1)}\sum_{j=0}^{i}\frac{(-\gamma)_j(-\gamma_1k)_{i-j}}{j!(i-j)!}.$$
In accordance with \eqref{eq2.28}, we obtain
$$\sum_{n=-\infty}^{+\infty}sign(x+2na)\sum_{k=0}^{+\infty}\frac{(-1)^k\left|x+2na\right|^k(t-\eta)^{-\beta_1k-\beta}}{k!}\sum_{i=0}^{+\infty}\frac{ (-\gamma-\gamma_1k)_i\delta^i(t-\eta)^{\alpha i}}{i!\Gamma(\alpha i-\beta_1k-\beta+1)}=$$
$$=\sum_{n=-\infty}^{+\infty}sign(x+2na)\sum_{k=0}^{+\infty}\frac{(-1)^k\left|x+2na\right|^k(t-\eta)^{-\beta_1k-\beta}}{k!}E_{\alpha, 1-\beta_1k-\beta}^{-\gamma-\gamma_1k}\left[\delta(t-\eta)^{\alpha }\right].$$

Next, we insert the derived result into the preceding equality and proceed with the computation:
$$\frac{\partial}{\partial t}\int\limits_0^t{{\varphi }_{0}}\left( \eta  \right)\sum_{n=-\infty}^{+\infty}sign(x+2na)\sum_{k=0}^{+\infty}\frac{(-1)^k\left|x+2na\right|^k(t-\eta)^{-\beta_1k-\beta}}{k!}E_{\alpha, 1-\beta_1k-\beta}^{-\gamma-\gamma_1k}\left[\delta(t-\eta)^{\alpha }\right]d\eta=$$
$$=\lim_{\eta \to t}\left[{{\varphi }_{0}}\left( \eta  \right)\sum_{n=-\infty}^{+\infty}sign(x+2na)\sum_{k=0}^{+\infty}\frac{(-1)^k\left|x+2na\right|^k(t-\eta)^{-\beta_1k-\beta}}{k!}E_{\alpha, 1-\beta_1k-\beta}^{-\gamma-\gamma_1k}\left[\delta(t-\eta)^{\alpha }\right]\right]+$$
$$+\int\limits_0^t{{\varphi }_{0}}\left( \eta  \right)\sum_{n=-\infty}^{+\infty}sign(x+2na)\sum_{k=0}^{+\infty}\frac{(-1)^k\left|x+2na\right|^k(t-\eta)^{-\beta_1k-\beta-1}}{k!}E_{\alpha, -\beta_1k-\beta}^{-\gamma-\gamma_1k}\left[\delta(t-\eta)^{\alpha }\right]d\eta.$$
Now let us carefully evaluate the limit that arises from taking the derivative:
$$\lim_{\eta \to t}\left[{{\varphi }_{0}}\left( \eta  \right)\sum_{n=-\infty}^{+\infty}sign(x+2na)\sum_{k=0}^{+\infty}\frac{(-1)^k\left|x+2na\right|^k(t-\eta)^{-\beta_1k-\beta}}{k!}\right]\times$$
$$\times\lim_{\eta \to t}E_{\alpha, 1-\beta_1k-\beta}^{-\gamma-\gamma_1k}\left[\delta(t-\eta)^{\alpha }\right]=$$
$$=\lim_{\eta \to t}\left[{{\varphi }_{0}}\left( \eta  \right)\sum_{n=-\infty}^{+\infty}sign(x+2na)(t-\eta)^{-\beta}\sum_{k=0}^{+\infty}\frac{(-1)^k\left|x+2na\right|^k(t-\eta)^{-\beta_1k}}{k!\Gamma(1-\beta_1k-\beta)}\right].$$
We rewrite the final expression using the Wright-type function and make use of its  property \cite{Pskhu 2}
$$\lim_{\left|z\right| \to \infty}z^2 e_{\alpha,{\beta}}^{\alpha-k,{\delta }}\left( z \right)=\frac{1}{\Gamma(-k-\alpha)\Gamma(\delta+2\beta)}:$$

$$\lim_{\eta \to t}\left[{{\varphi }_{0}}\left( \eta  \right)\sum_{n=-\infty}^{+\infty}sign(x+2na)(t-\eta)^{-\beta}e_{1,\beta_1}^{1,1-\beta}\left(-\frac{\left|x+2na\right|}{(t-\eta)^{\beta_1}}\right)\right]=$$
$$=\lim_{\eta \to t}\left[{{\varphi }_{0}}\left( \eta  \right)\sum_{n=-\infty}^{+\infty}\frac{sign(x+2na)}{(x+2na)^2}\left(-\frac{\left|x+2na\right|}{(t-\eta)^{\beta_1}}\right)^2e_{1,\beta_1}^{1,1-\beta}\left(-\frac{\left|x+2na\right|}{(t-\eta)^{\beta_1}}\right)\right]=0.$$
Therefore,
$${}^{PRL}D_{0t}^{\alpha ,\,\beta ,\,\gamma ,\,\delta }\left(\int\limits_{0}^{t}{{{\varphi }_{0}}\left( \eta  \right){{G}_{s}}\left( t,x,\eta ,0 \right)d\eta }\right)=\int\limits_0^t{{\varphi }_{0}}\left( \eta  \right)\sum_{n=-\infty}^{+\infty}sign(x+2na)\times$$
\begin{equation}\label{eq14}
\times\sum_{k=0}^{+\infty}\frac{(-1)^k\left|x+2na\right|^k(t-\eta)^{-\beta_1k-\beta-1}}{k!}E_{\alpha, -\beta_1k-\beta}^{-\gamma-\gamma_1k}\left[\delta(t-\eta)^{\alpha }\right]d\eta.
\end{equation}

Similarly, for the term involving $\varphi_1(\eta)$ in \eqref{eq13}, if we repeat the same calculations as above, we obtain the following:
$${}^{PRL}D_{0t}^{\alpha ,\,\beta ,\,\gamma ,\,\delta }\left(\int\limits_{0}^{t}{{{\varphi }_{1}}\left( \eta  \right){{G}_{s}}\left( t,x,\eta ,a \right)d\eta }\right)=\int\limits_0^t{{\varphi }_{0}}\left( \eta  \right)\sum_{n=-\infty}^{+\infty}sign(x+(2n+1)a)\times$$
\begin{equation}\label{eq15}
\times\sum_{k=0}^{+\infty}\frac{(-1)^k\left|x+(2n+1)a\right|^k(t-\eta)^{-\beta_1k-\beta-1}}{k!}E_{\alpha, -\beta_1k-\beta}^{-\gamma-\gamma_1k}\left[\delta(t-\eta)^{\alpha }\right]d\eta.
\end{equation}

Next, we proceed to the term involving $\tau(s):$
$${}^{PRL}D_{0t}^{\alpha ,\,\beta ,\,\gamma ,\,\delta }\left(\int\limits_{0}^{a}{\tau \left( s \right) G\left( t,x,0,s \right)ds}\right)=\frac{\partial}{\partial t}{}^{P}I_{0 t}^{\alpha, 1-\beta,-\gamma, \delta}\left(\int\limits_{0}^{a}{\tau \left( s \right) G\left( t,x,0,s \right)ds}\right)=$$
$$=\frac{\partial}{\partial t}\int\limits_0^t(t-y)^{-\beta} E_{\alpha, 1-\beta}^{-\gamma}\left[\mathcal{\delta}(t-y)^\alpha\right]dy\int\limits_{0}^{a}{\tau \left( s \right) G\left( y,x,0,s\right)ds}=$$
\begin{equation}\label{eq16}
    =\frac{\partial}{\partial t}\int\limits_0^a\tau \left( s \right)\,ds\int\limits_{0}^{t}{(t-y)^{-\beta} E_{\alpha, 1-\beta}^{-\gamma}\left[\mathcal{\delta}(t-y)^\alpha\right] G\left( y,x,0,s\right)dy}.
\end{equation}
Now let us evaluate the following integral:
$$\int\limits_{0}^{t}{(t-y)^{-\beta} E_{\alpha, 1-\beta}^{-\gamma}\left[\mathcal{\delta}(t-y)^\alpha\right]G\left( y,x,0,s\right)dy}=\int\limits_{0}^{t}{(t-y)^{-\beta} E_{\alpha, 1-\beta}^{-\gamma}\left[\mathcal{\delta}(t-y)^\alpha\right]\times}$$
$$\times\sum\limits_{n=-\infty }^{+\infty }{\sum\limits_{k=0}^{+\infty }{\frac{{{\left( -1 \right)}^{k}}\left[ {{\left| x-s+2na \right|}^{k}}-{{\left| x+s+2na \right|}^{k}} \right]}{2\Gamma(k+1)}}}y^{-\beta_1k+{{\beta }_{1}}-1}E_{\alpha, \beta_1-\beta_1k}^{{\gamma }_{1}-{{\gamma }_{1}}k}\left[\delta y^\alpha\right]dy=$$
$$=\int\limits_{0}^{t}{(t-y)^{-\beta} \sum\limits_{n=-\infty }^{+\infty }{\sum\limits_{k=0}^{+\infty }{\frac{{{\left( -1 \right)}^{k}}\left[ {{\left| x-s+2na \right|}^{k}}-{{\left| x+s+2na \right|}^{k}} \right]}{2\Gamma(k+1)}}}y^{-{{\beta }_{1}k}+\beta_1-1}\times}$$
$$\times \sum_{j=0}^{+\infty}\frac{(-\gamma)_j\delta^j(t-y)^{\alpha j}}{j!\Gamma(\alpha j+1-\beta)}\sum_{i=0}^{+\infty}\frac{({{\gamma }_{1}}-{{\gamma }_{1}}k)_i\delta^iy^{\alpha i}}{i!\,\Gamma(\alpha i+\beta_1-\beta_1k)}dy.$$
We apply \eqref{eq2.27} to the sums with respect to $i$ and $j$ and perform the substitution $y=tz$:
$$\sum\limits_{n=-\infty }^{+\infty }{\sum\limits_{k=0}^{+\infty }{\frac{{{\left( -1 \right)}^{k}}\left[ {{\left| x-s+2na \right|}^{k}}-{{\left| x+s+2na \right|}^{k}} \right]}{2\Gamma(k+1)}}}\times$$
$$\times \sum_{j=0}^{+\infty}\sum_{i=0}^{j}\frac{(-\gamma)_i({{\gamma }_{1}}-{{\gamma }_{1}}k)_{j-i}\delta^j}{i!(j-i)!\Gamma(\alpha i+1-\beta)\Gamma(\alpha j -\alpha i+\beta_1-\beta_1k)}\int\limits_{0}^{t}{(t-y)^{\alpha i-\beta} y^{\alpha j-\alpha i-{{\beta }_{1}k}+\beta_1-1}dy}=$$
$$=\sum\limits_{n=-\infty }^{+\infty }{\sum\limits_{k=0}^{+\infty }{\frac{{{\left( -1 \right)}^{k}}\left[ {{\left| x-s+2na \right|}^{k}}-{{\left| x+s+2na \right|}^{k}} \right]}{2\Gamma(k+1)}}} \sum_{j=0}^{+\infty}\sum_{i=0}^{j}\frac{(-\gamma)_i({{\gamma }_{1}}-{{\gamma }_{1}}k)_{j-i}\delta^jt^{\alpha j-\beta_1-\beta_{1}k}}{i!(j-i)!\Gamma(\alpha j -\beta_1-\beta_1k+1)}.$$
Using \eqref{eq2.28}, we obtain
$$\sum\limits_{n=-\infty }^{+\infty }{\sum\limits_{k=0}^{+\infty }{\frac{{{\left( -1 \right)}^{k}}\left[ {{\left| x-s+2na \right|}^{k}}-{{\left| x+s+2na \right|}^{k}} \right]}{2\Gamma(k+1)}}}\sum_{j=0}^{+\infty}\frac{({-{\gamma }_{1}}-{{\gamma }_{1}}k)_{j}\delta^j t^{\alpha j-\beta_1-\beta_1k}}{j!\Gamma(\alpha j -\beta_1-\beta_1k+1)}.$$
Now, according to \eqref{eq15}, we take the derivative of the last equality with respect to $t$ and express the obtained result in terms ofhe generalized Mittag–Leffler function:
$$\frac{\partial}{\partial t}\sum\limits_{n=-\infty }^{+\infty }{\sum\limits_{k=0}^{+\infty }{\frac{{{\left( -1 \right)}^{k}}\left[ {{\left| x-s+2na \right|}^{k}}-{{\left| x+s+2na \right|}^{k}} \right]}{2\Gamma(k+1)}}}\sum_{j=0}^{+\infty}\frac{({-{\gamma }_{1}}-{{\gamma }_{1}}k)_{j}\delta^j t^{\alpha j-\beta_1-\beta_1k}}{j!\Gamma(\alpha j -\beta_1-\beta_1k+1)}=$$
$$=\sum\limits_{n=-\infty }^{+\infty }{\sum\limits_{k=0}^{+\infty }{\frac{{{\left( -1 \right)}^{k}}\left[ {{\left| x-s+2na \right|}^{k}}-{{\left| x+s+2na \right|}^{k}} \right]}{2\Gamma(k+1)}}}\sum_{j=0}^{+\infty}\frac{({-{\gamma }_{1}}-{{\gamma }_{1}}k)_{j}\delta^j t^{\alpha j-\beta_1-\beta_1k-1}}{j!\Gamma(\alpha j -\beta_1-\beta_1k)}=$$
$$=\sum\limits_{n=-\infty }^{+\infty }{\sum\limits_{k=0}^{+\infty }{\frac{{{\left( -1 \right)}^{k}}\left[ {{\left| x-s+2na \right|}^{k}}-{{\left| x+s+2na \right|}^{k}} \right]}{2\Gamma(k+1)}}}t^{-\beta_1-\beta_{1}k-1}E_{\alpha, {-\beta_1-\beta_1k}}^{-\gamma _{1}-{{\gamma }_{1}}k}\left[\delta t^\alpha\right].$$
Thus, we determine that 
$${}^{PRL}D_{0t}^{\alpha ,\,\beta ,\,\gamma ,\,\delta }\left(\int\limits_{0}^{a}{\tau \left( s\right) G\left( t,x,0,s \right)ds}\right)=$$
\begin{equation}\label{eq17}
    =\int\limits_{0}^{a}\tau(s)\sum\limits_{n=-\infty }^{+\infty }{\sum\limits_{k=0}^{+\infty }{\frac{{{\left( -1 \right)}^{k}}\left[ {{\left| x-s+2na \right|}^{k}}-{{\left| x+s+2na \right|}^{k}} \right]}{2\Gamma(k+1)}}}t^{-\beta_1-\beta_{1}k-1}E_{\alpha, {-\beta_1-\beta_1k}}^{-\gamma _{1}-{{\gamma }_{1}}k}\left[\delta t^\alpha\right]ds.
\end{equation}

Finally, we compute the last term of equation \eqref{eq13}:
$${}^{PRL}D_{0t}^{\alpha ,\,\beta ,\,\gamma ,\,\delta }\left(\int\limits_{0}^{t}{\int\limits_{0}^{a}{f\left( \eta ,s \right)G\left( t,x,\eta ,s \right)dsd\eta }}\right)=$$
$$={\frac{\partial^{}}{\partial t}}{}^{P}I_{0 t}^{\alpha, 1-\beta,-\gamma, \delta} \left(\int\limits_{0}^{t}{\int\limits_{0}^{a}{f\left( \eta ,s \right)G\left( t,x,\eta ,s \right)dsd\eta }}\right)=$$
$$=\frac{\partial}{\partial t}\int\limits_0^t(t-y)^{-\beta} E_{\alpha, 1-\beta}^{-\gamma}\left[\mathcal{\delta}(t-y)^\alpha\right]dy\int\limits_{0}^{a}ds{\int\limits_{0}^{y}{f\left( \eta ,s \right)G\left( y,x,\eta ,s \right)d\eta }}=$$
$$=\frac{\partial}{\partial t}\int\limits_{0}^{a}ds\int\limits_0^tf\left( \eta ,s \right)d\eta\int\limits_{\eta}^{t}(t-y)^{-\beta} E_{\alpha, 1-\beta}^{-\gamma}\left[\mathcal{\delta}(t-y)^\alpha\right]{G\left( y,x,\eta ,s \right)dy. }$$
Here we replace the function $G\left( y,x,\eta ,s \right)$ by its representation given in \eqref{eq6*}:
$$\frac{\partial}{\partial t}\int\limits_{0}^{a}ds\int\limits_0^tf\left( \eta ,s \right)\left[\int\limits_{\eta}^{t}(t-y)^{-\beta} E_{\alpha, 1-\beta}^{-\gamma}\left[\mathcal{\delta}(t-y)^\alpha\right]\times\right.$$
$$\times{\frac{{(y-\eta)^{{{\beta }_{1}}-1}}}{2}{ {{E}_{12}}\left( \left. \begin{matrix}
   -{{\gamma }_{1}},1,{{\gamma }_{1}};\,\,\,\,\,\,\,\,\,\,\,\,\,\,\,\,\,\,\,\,\,\,\,\,\,\,\,  \\
   -{{\beta }_{1}},\alpha ,{{\beta }_{1}};-{{\gamma }_{1}},{{\gamma }_{1}};1,1;1,1  \\
\end{matrix} \right|\begin{matrix}
   -\left| x-s \right|{(y-\eta)^{-{{\beta }_{1}}}}  \\
   \delta {(y-\eta)^{\alpha }}  \\
\end{matrix} \right) }dy }+$$
\begin{equation}\label{eq19*}
    \left.+\int\limits_{\eta}^{t}(t-y)^{-\beta} E_{\alpha, 1-\beta}^{-\gamma}\left[\mathcal{\delta}(t-y)^\alpha\right]K(y,x,\eta,s)dy\right]d\eta.
\end{equation}
Now we compute the first part of \eqref{eq19*}:
$$\frac{\partial}{\partial t}\int\limits_{0}^{a}ds\int\limits_0^tf\left( \eta ,s \right)d\eta\int\limits_{\eta}^{t}(t-y)^{-\beta} E_{\alpha, 1-\beta}^{-\gamma}\left[\mathcal{\delta}(t-y)^\alpha\right]\times$$
$$\times{\frac{{(y-\eta)^{{{\beta }_{1}}-1}}}{2}{ {{E}_{12}}\left( \left. \begin{matrix}
   -{{\gamma }_{1}},1,{{\gamma }_{1}};\,\,\,\,\,\,\,\,\,\,\,\,\,\,\,\,\,\,\,\,\,\,\,\,\,\,\,  \\
   -{{\beta }_{1}},\alpha ,{{\beta }_{1}};-{{\gamma }_{1}},{{\gamma }_{1}};1,1;1,1  \\
\end{matrix} \right|\begin{matrix}
   -\left| x-s \right|{(y-\eta)^{-{{\beta }_{1}}}}  \\
   \delta {(y-\eta)^{\alpha }}  \\
\end{matrix} \right) }dy }=$$
$$=\frac{\partial}{\partial t}\int\limits_{0}^{a}ds\int\limits_0^tf\left( \eta ,s \right){\sum\limits_{k=0}^{+\infty }{\frac{{{\left( -1 \right)}^{k}} {{\left| x-s \right|}^{k}}}{2\Gamma(k+1)}}}d\eta\times$$
$$\times\int\limits_{\eta}^{t}(t-y)^{-\beta} (y-\eta)^{-\beta_1k+{{\beta }_{1}}-1}E_{\alpha, 1-\beta}^{-\gamma}\left[\mathcal{\delta}(t-y)^\alpha\right]E_{\alpha, \beta_1-\beta_1k}^{{\gamma }_{1}-{{\gamma }_{1}}k}\left[\delta (y-\eta)^\alpha\right]dy.$$
Let us evaluate the integral with respect to $y.$ We apply \eqref{eq2.27} to the Prabhakar functions and perform the substitution $y=(t-\eta)z+\eta.$ Using \eqref{eq2.28}, we obtain
$$\int\limits_{\eta}^{t}(t-y)^{-\beta} (y-\eta)^{-\beta_1k+{{\beta }_{1}}-1}E_{\alpha, 1-\beta}^{-\gamma}\left[\mathcal{\delta}(t-y)^\alpha\right]E_{\alpha, \beta_1-\beta_1k}^{{\gamma }_{1}-{{\gamma }_{1}}k}\left[\delta (y-\eta)^\alpha\right]dy=$$
$$=\sum_{j=0}^{+\infty}\sum_{i=0}^{j}\frac{(-\gamma)_i({{\gamma }_{1}}-{{\gamma }_{1}}k)_{j-i}\delta^j}{i!(j-i)!\Gamma(\alpha i+1-\beta)\Gamma(\alpha j -\alpha i+\beta_1-\beta_1k)}\int\limits_{\eta}^{t}{(t-y)^{\alpha i-\beta} (y-\eta)^{\alpha j-\alpha i-{{\beta }_{1}k}+\beta_1-1}dy}=$$
$$=\sum_{j=0}^{+\infty}\sum_{i=0}^{j}\frac{(-\gamma)_i({{\gamma }_{1}}-{{\gamma }_{1}}k)_{j-i}\delta^j(t-\eta)^{\alpha j-\beta_1-\beta_{1}k}}{i!(j-i)!\Gamma(\alpha j -\beta_1-\beta_1k+1)}=\sum_{j=0}^{+\infty}\frac{({-{\gamma }_{1}}-{{\gamma }_{1}}k)_{j}\delta^j (t-\eta)^{\alpha j-\beta_1-\beta_1k}}{j!\Gamma(\alpha j -\beta_1-\beta_1k+1)}.$$
We continue the computation using the obtained result:
$$\frac{\partial}{\partial t}\int\limits_{0}^{a}ds\int\limits_0^tf\left( \eta ,s \right){\sum\limits_{k=0}^{+\infty }{\frac{{{\left( -1 \right)}^{k}} {{\left| x-s \right|}^{k}}}{2\Gamma(k+1)}}}\sum_{j=0}^{+\infty}\frac{({-{\gamma }_{1}}-{{\gamma }_{1}}k)_{j}\delta^j (t-\eta)^{\alpha j-\beta_1-\beta_1k}}{j!\Gamma(\alpha j -\beta_1-\beta_1k+1)}d\eta=$$
$$=\frac{1}{2}\frac{\partial}{\partial t}\int\limits_{0}^{x}ds\int\limits_0^tf\left( \eta ,s \right)(t-\eta)^{-\beta_1}{\sum\limits_{k=0}^{+\infty }{\frac{1}{\Gamma(k+1)}}}\left(-\frac{x-s}{(t-\eta)^{\beta_1}}\right)^{k}E_{\alpha, {1-\beta_1-\beta_1k}}^{-\gamma _{1}-{{\gamma }_{1}}k}\left[\delta (t-\eta)^\alpha\right]d\eta+$$
$$+\frac{1}{2}\frac{\partial}{\partial t}\int\limits_{x}^{a}ds\int\limits_0^tf\left( \eta ,s \right)(t-\eta)^{-\beta_1}{\sum\limits_{k=0}^{+\infty }{\frac{1}{\Gamma(k+1)}}}\left(-\frac{s-x}{(t-\eta)^{\beta_1}}\right)^{k}E_{\alpha, {1-\beta_1-\beta_1k}}^{-\gamma _{1}-{{\gamma }_{1}}k}\left[\delta (t-\eta)^\alpha\right]d\eta.$$
Next, we differentiate the resulting expression with respect to the parameter $t:$
$$\frac{1}{2}\lim_{\eta \to t}\int\limits_{0}^{x}f\left( \eta ,s \right)(t-\eta)^{-\beta_1}{\sum\limits_{k=0}^{+\infty }{\frac{(-1)^k}{\Gamma(k+1)}}}\left(\frac{x-s}{(t-\eta)^{\beta_1}}\right)^{k}E_{\alpha, {1-\beta_1-\beta_1k}}^{-\gamma _{1}-{{\gamma }_{1}}k}\left[\delta (t-\eta)^\alpha\right]ds+$$
$$+\frac{1}{2}\lim_{\eta \to t}\int\limits_{x}^{a}f\left( \eta ,s \right)(t-\eta)^{-\beta_1}{\sum\limits_{k=0}^{+\infty }{\frac{(-1)^k}{\Gamma(k+1)}}}\left(\frac{s-x}{(t-\eta)^{\beta_1}}\right)^{k}E_{\alpha, {1-\beta_1-\beta_1k}}^{-\gamma _{1}-{{\gamma }_{1}}k}\left[\delta (t-\eta)^\alpha\right]ds+$$
$$+\frac{1}{2}\int\limits_{0}^{a}ds\int\limits_0^tf\left( \eta ,s \right){\sum\limits_{k=0}^{+\infty }{\frac{{{\left( -1 \right)}^{k}} {{\left| x-s \right|}^{k}}}{\Gamma(k+1)}}}(t-\eta)^{-\beta_1-\beta_{1}k-1}E_{\alpha, {-\beta_1-\beta_1k}}^{-\gamma _{1}-{{\gamma }_{1}}k}\left[\delta (t-\eta)^\alpha\right]d\eta.$$
In the first and second integrals of the final expression, we perform the substitutions $z_1=\frac{x-s}{(t-\eta)^{\beta_1}}$ and $z_1=\frac{s-x}{(t-\eta)^{\beta_1}}$, respectively:
$$\frac{1}{2}\lim_{\eta \to t}\int\limits_{0}^{\frac{x}{(t-\eta)^{\beta_1}}}f\left( \eta ,x-(t-\eta)^{\beta_1}z_1 \right){\sum\limits_{k=0}^{+\infty }{\frac{(-1)^k{z_1}^{k}}{\Gamma(k+1)}}}E_{\alpha, {1-\beta_1-\beta_1k}}^{-\gamma _{1}-{{\gamma }_{1}}k}\left[\delta (t-\eta)^\alpha\right]dz_1+$$
$$+\frac{1}{2}\lim_{\eta \to t}\int\limits_{0}^{\frac{a-x}{(t-\eta)^{\beta_1}}}f\left( \eta ,x+(t-\eta)^{\beta_1}z_2 \right){\sum\limits_{k=0}^{+\infty }{\frac{(-1)^k{z_2}^{k}}{\Gamma(k+1)}}}E_{\alpha, {1-\beta_1-\beta_1k}}^{-\gamma _{1}-{{\gamma }_{1}}k}\left[\delta (t-\eta)^\alpha\right]dz_2+$$
$$+\frac{1}{2}\int\limits_{0}^{a}ds\int\limits_0^tf\left( \eta ,s \right){\sum\limits_{k=0}^{+\infty }{\frac{{{\left( -1 \right)}^{k}} {{\left| x-s \right|}^{k}}}{\Gamma(k+1)}}}(t-\eta)^{-\beta_1-\beta_{1}k-1}E_{\alpha, {-\beta_1-\beta_1k}}^{-\gamma _{1}-{{\gamma }_{1}}k}\left[\delta (t-\eta)^\alpha\right]d\eta=$$
$$=\frac{f(t,x)}{2}\int\limits_{0}^{+\infty}{\sum\limits_{k=0}^{+\infty }{\frac{(-1)^k{z_1}^{k}}{\Gamma(k+1)\Gamma({1-\beta_1-\beta_1k})}}}dz_1+$$
$$+\frac{f(t,x)}{2}\int\limits_{0}^{+\infty}{\sum\limits_{k=0}^{+\infty }{\frac{(-1)^k{z_2}^{k}}{\Gamma(k+1)\Gamma({1-\beta_1-\beta_1k})}}}dz_2+$$
$$+\frac{1}{2}\int\limits_{0}^{a}ds\int\limits_0^tf\left( \eta ,s \right){\sum\limits_{k=0}^{+\infty }{\frac{{{\left( -1 \right)}^{k}} {{\left| x-s \right|}^{k}}}{\Gamma(k+1)}}}(t-\eta)^{-\beta_1-\beta_{1}k-1}E_{\alpha, {-\beta_1-\beta_1k}}^{-\gamma _{1}-{{\gamma }_{1}}k}\left[\delta (t-\eta)^\alpha\right]d\eta.$$
After computing the integrals with respect to $z_1$ and $z_2$
and applying property 
$$\lim_{|z| \to \infty} ze_{\alpha,\beta}^{\mu,\delta}(z)=-\frac{1}{\Gamma(\mu-\alpha)\Gamma(\delta+\beta)}$$
of the Wright function, we arrive at the following result:
\begin{equation}\label{eq20*}
    f(t,x)+\frac{1}{2}\int\limits_{0}^{a}ds\int\limits_0^tf\left( \eta ,s \right){\sum\limits_{k=0}^{+\infty }{\frac{{{\left( -1 \right)}^{k}} {{\left| x-s \right|}^{k}}}{\Gamma(k+1)}}}(t-\eta)^{-\beta_1-\beta_{1}k-1}E_{\alpha, {-\beta_1-\beta_1k}}^{-\gamma _{1}-{{\gamma }_{1}}k}\left[\delta (t-\eta)^\alpha\right]d\eta.
\end{equation}
Next, we evaluate the second term of \eqref{eq19*} 
\begin{equation}\label{eq21*}
    \frac{\partial}{\partial t}\int\limits_{0}^{a}ds\int\limits_0^tf\left( \eta ,s \right)d\eta\int\limits_{\eta}^{t}(t-y)^{-\beta} E_{\alpha, 1-\beta}^{-\gamma}\left[\mathcal{\delta}(t-y)^\alpha\right]K(y,x,\eta,s)dy
\end{equation}
using the same approach. Firstly, we calculate the following integral:
$$\int\limits_{\eta}^{t}(t-y)^{-\beta} E_{\alpha, 1-\beta}^{-\gamma}\left[\mathcal{\delta}(t-y)^\alpha\right]{K\left( y,x,\eta ,s \right)dy}=\int\limits_{\eta}^{t}{(t-y)^{-\beta} E_{\alpha, 1-\beta}^{-\gamma}\left[\mathcal{\delta}(t-y)^\alpha\right]\times}$$
$$\times\sum\limits_{n=-\infty }^{+\infty }{}'{\sum\limits_{k=0}^{+\infty }{\frac{{{\left( -1 \right)}^{k}}{\left| x-s+2na \right|}^{k}}{2\Gamma(k+1)}}}(y-\eta)^{-\beta_1k+{{\beta }_{1}}-1}E_{\alpha, \beta_1-\beta_1k}^{{\gamma }_{1}-{{\gamma }_{1}}k}\left[\delta (y-\eta)^\alpha\right]dy-$$
$$-\int\limits_{\eta}^{t}{(t-y)^{-\beta} E_{\alpha, 1-\beta}^{-\gamma}\left[\mathcal{\delta}(t-y)^\alpha\right]\times}$$
$$\times\sum\limits_{n=-\infty }^{+\infty }{\sum\limits_{k=0}^{+\infty }{\frac{{{\left( -1 \right)}^{k}}{\left| x+s+2na \right|}^{k}}{2\Gamma(k+1)}}}(y-\eta)^{-\beta_1k+{{\beta }_{1}}-1}E_{\alpha, \beta_1-\beta_1k}^{{\gamma }_{1}-{{\gamma }_{1}}k}\left[\delta (y-\eta)^\alpha\right]dy=$$
$$=\sum\limits_{n=-\infty }^{+\infty }{}'{\sum\limits_{k=0}^{+\infty }{\frac{{{\left( -1 \right)}^{k}}{\left| x-s+2na \right|}^{k}}{2\Gamma(k+1)}}}\sum_{j=0}^{+\infty}\frac{({-{\gamma }_{1}}-{{\gamma }_{1}}k)_{j}\delta^j (t-\eta)^{\alpha j-\beta_1-\beta_1k}}{j!\Gamma(\alpha j -\beta_1-\beta_1k+1)}-$$
$$-\sum\limits_{n=-\infty }^{+\infty }{\sum\limits_{k=0}^{+\infty }{\frac{{{\left( -1 \right)}^{k}}{\left| x+s+2na \right|}^{k}}{2\Gamma(k+1)}}}\sum_{j=0}^{+\infty}\frac{({-{\gamma }_{1}}-{{\gamma }_{1}}k)_{j}\delta^j (t-\eta)^{\alpha j-\beta_1-\beta_1k}}{j!\Gamma(\alpha j -\beta_1-\beta_1k+1)}$$
We substitute the obtained result into \eqref{eq21*} and differentiate with respect to the parameter $t$:
$$\frac{\partial}{\partial t}\int\limits_{0}^{a}ds\int\limits_0^tf\left( \eta ,s \right)\left[\sum\limits_{n=-\infty }^{+\infty }{}'{\sum\limits_{k=0}^{+\infty }{\frac{{{\left( -1 \right)}^{k}}{\left| x-s+2na \right|}^{k}}{2\Gamma(k+1)}}}\sum_{j=0}^{+\infty}\frac{({-{\gamma }_{1}}-{{\gamma }_{1}}k)_{j}\delta^j (t-\eta)^{\alpha j-\beta_1-\beta_1k}}{j!\Gamma(\alpha j -\beta_1-\beta_1k+1)}-\right.$$
$$\left.-\sum\limits_{n=-\infty }^{+\infty }{\sum\limits_{k=0}^{+\infty }{\frac{{{\left( -1 \right)}^{k}}{\left| x+s+2na \right|}^{k}}{2\Gamma(k+1)}}}\sum_{j=0}^{+\infty}\frac{({-{\gamma }_{1}}-{{\gamma }_{1}}k)_{j}\delta^j (t-\eta)^{\alpha j-\beta_1-\beta_1k}}{j!\Gamma(\alpha j -\beta_1-\beta_1k+1)}\right]d\eta=$$
$$=\int\limits_{0}^{a}\lim_{\eta \to t}\left[f\left( \eta ,s \right)\sum\limits_{n=-\infty }^{+\infty }{}'{\sum\limits_{k=0}^{+\infty }{\frac{{{\left( -1 \right)}^{k}}{\left| x-s+2na \right|}^{k}}{2\Gamma(k+1)}}}t^{-\beta_1-\beta_1k}E_{\alpha,1-\beta_1-\beta_1k}^{-\gamma_1-\gamma_1k}\left[\delta(t-\eta)^\alpha\right]-\right.$$
$$\left.-f\left( \eta ,s \right)\sum\limits_{n=-\infty }^{+\infty }{\sum\limits_{k=0}^{+\infty }{\frac{{{\left( -1 \right)}^{k}}{\left| x+s+2na \right|}^{k}}{2\Gamma(k+1)}}}t^{-\beta_1-\beta_1k}E_{\alpha,1-\beta_1-\beta_1k}^{-\gamma_1-\gamma_1k}\left[\delta(t-\eta)^\alpha\right]\right]ds+$$
$$+\int\limits_{0}^{a}ds\int\limits_0^tf\left( \eta ,s \right)\left[\sum\limits_{n=-\infty }^{+\infty }{}'{\sum\limits_{k=0}^{+\infty }{\frac{{{\left( -1 \right)}^{k}{\left| x-s+2na \right|}^{k}}}{2\Gamma(k+1)}}}(t-\eta)^{-\beta_1-\beta_{1}k-1}E_{\alpha, {-\beta_1-\beta_1k}}^{-\gamma _{1}-{{\gamma }_{1}}k}\left[\delta (t-\eta)^\alpha\right]-\right.$$
$$\left.-\sum\limits_{n=-\infty }^{+\infty }{\sum\limits_{k=0}^{+\infty }{\frac{{{\left( -1 \right)}^{k}{\left| x+s+2na \right|}^{k}}}{2\Gamma(k+1)}}}(t-\eta)^{-\beta_1-\beta_{1}k-1}E_{\alpha, {-\beta_1-\beta_1k}}^{-\gamma _{1}-{{\gamma }_{1}}k}\left[\delta (t-\eta)^\alpha\right]\right]d\eta.$$
At this point, the value of the following limit is zero, according to the Lemma 3.1:
$$\lim_{\eta \to t}\left[\sum\limits_{n=-\infty }^{+\infty }{}'{\sum\limits_{k=0}^{+\infty }{\frac{{{\left( -1 \right)}^{k}}{\left| x-s+2na \right|}^{k}}{2\Gamma(k+1)}}}t^{-\beta_1-\beta_1k}E_{\alpha,1-\beta_1-\beta_1k}^{-\gamma_1-\gamma_1k}\left[\delta(t-\eta)^\alpha\right]-\right.$$
$$\left.-\sum\limits_{n=-\infty }^{+\infty }{\sum\limits_{k=0}^{+\infty }{\frac{{{\left( -1 \right)}^{k}}{\left| x+s+2na \right|}^{k}}{2\Gamma(k+1)}}}t^{-\beta_1-\beta_1k}E_{\alpha,1-\beta_1-\beta_1k}^{-\gamma_1-\gamma_1k}\left[\delta(t-\eta)^\alpha\right]\right]=0.$$
Hence,
$$\frac{\partial}{\partial t}\int\limits_{0}^{a}ds\int\limits_0^tf\left( \eta ,s \right)d\eta\int\limits_{\eta}^{t}(t-y)^{-\beta} E_{\alpha, 1-\beta}^{-\gamma}\left[\mathcal{\delta}(t-y)^\alpha\right]K(y,x,\eta,s)dy=$$
$$=\int\limits_{0}^{a}ds\int\limits_0^tf\left( \eta ,s \right)\left[\sum\limits_{n=-\infty }^{+\infty }{}'{\sum\limits_{k=0}^{+\infty }{\frac{{{\left( -1 \right)}^{k}{\left| x-s+2na \right|}^{k}}}{2\Gamma(k+1)}}}(t-\eta)^{-\beta_1-\beta_{1}k-1}E_{\alpha, {-\beta_1-\beta_1k}}^{-\gamma _{1}-{{\gamma }_{1}}k}\left[\delta (t-\eta)^\alpha\right]-\right.$$
\begin{equation}\label{eq22*}
    \left.-\sum\limits_{n=-\infty }^{+\infty }{\sum\limits_{k=0}^{+\infty }{\frac{{{\left( -1 \right)}^{k}{\left| x+s+2na \right|}^{k}}}{2\Gamma(k+1)}}}(t-\eta)^{-\beta_1-\beta_{1}k-1}E_{\alpha, {-\beta_1-\beta_1k}}^{-\gamma _{1}-{{\gamma }_{1}}k}\left[\delta (t-\eta)^\alpha\right]\right]d\eta.
\end{equation}

Summing expressions \eqref{eq20*} and \eqref{eq22*}, we arrive at the following result:
$${}^{PRL}D_{0t}^{\alpha ,\,\beta ,\,\gamma ,\,\delta }\left(\int\limits_{0}^{t}{\int\limits_{0}^{a}{f\left( \eta ,s \right)G\left( t,x,\eta ,s \right)dsd\eta }}\right)=$$
$$=f(t,x)+\int\limits_{0}^{a}\int\limits_0^tf\left( \eta ,s \right)\sum\limits_{n=-\infty }^{+\infty }{\sum\limits_{k=0}^{+\infty }{\frac{{{\left( -1 \right)}^{k}}\left[ {{\left| x-s+2na \right|}^{k}}-{{\left| x+s+2na \right|}^{k}} \right]}{2\Gamma(k+1)}}}\times$$
\begin{equation}\label{eq19}
    \times(t-\eta)^{-\beta_1-\beta_{1}k-1}E_{\alpha, {-\beta_1-\beta_1k}}^{-\gamma _{1}-{{\gamma }_{1}}k}\left[\delta (t-\eta)^\alpha\right]d\eta ds.
\end{equation}

If we examine the limiting behavior of \eqref{eq14}, \eqref{eq15} and \eqref{eq19} at the point $\eta=t$, it is not difficult to see that all of them are equal to $0$ according to the Lemma 3.1. In \eqref{eq17}, there is no singular point. 

Therefore, for our solution, the condition
$${}^{PRL}D_{0t}^{\alpha ,\,\beta ,\,\gamma ,\,\delta }u\left( t,x \right)\in C\left( D \right)$$
is satisfied.

We have thoroughly verified the class to which our solution belongs. Now we proceed to check whether the obtained solution indeed satisfies the equation \eqref{eq2.1} in the manuscript. 

Let us introduce the following notation:
$$u_0(t,x)=u(t,x)-u_1(t,x),\,\,\,u_1(t,x)=\int\limits_0^t\int\limits_0^af(\eta,s)G(t,x,\eta,s)ds d\eta.$$
First, we show that $\mathop{Lu_0=0}.$ Let us begin by calculating the term involving the function $\varphi_0(\eta).$ Since we have computed the value of $\frac{\partial^2}{\partial x^2}{{G}_{s}}\left( t,x,\eta ,0 \right)$ above, we can write the following accordingly:
$$\frac{\partial^2}{\partial x^2}\int\limits_{0}^{t}{{{\varphi }_{0}}\left( \eta  \right){{G}_{s}}\left( t,x,\eta ,0 \right)d\eta }=$$
$$=\int\limits_{0}^{t}{{\varphi }_{0}}(\eta)\sum_{n=-\infty}^{+\infty}sign(x+2na)\sum_{k=0}^{+\infty}\frac{(-1)^k\left|x+2na\right|^k(t-\eta)^{-\beta_1k-\beta-1}}{k!}E_{\alpha, -\beta_1k-\beta}^{-\gamma-\gamma_1k}\left[\delta(t-\eta)^{\alpha }\right]d\eta.$$
Comparing the last equality with \eqref{eq14}, we can see that they are equal. So their difference equals zero. The term containing $\varphi_1(\eta)$ can be demonstrated analogously to the term containing $\varphi_0(\eta)$.

Hence, we proceed to the term involving $\tau(\xi).$ Using the value of $\frac{\partial^2}{\partial x^2}{G}\left( t,x,\eta ,s \right)$, we obtain the following expression:
$$\frac{\partial^2}{\partial x^2}\int\limits_{0}^{a}{\tau \left( s \right)G\left( t,x,0,s \right)ds}=$$
$$=\int\limits_{0}^{a}\tau \left( s \right)\sum\limits_{n=-\infty }^{+\infty }{\sum\limits_{k=0}^{+\infty }{\frac{{{\left( -1 \right)}^{k}}\left[ {{\left| x-s+2na \right|}^{k}}-{{\left| x+s+2na \right|}^{k}} \right]}{2\Gamma(k+1)}}}t^{-\beta_{1}k-\beta_1-1}E_{\alpha, {-\beta_1-\beta_1k}}^{-\gamma _{1}-{{\gamma }_{1}}k}\left[\delta t^\alpha\right]ds.$$
We also see that the difference between \eqref{eq17} and the last equality is equal to zero. Thus we have proven that $\mathop{Lu_0=0}.$

From the value of $\frac{\partial^2}{\partial x^2}{G}\left( t,x,\eta ,s \right)$, we determine that 
$$\frac{\partial^2}{\partial x^2}\int\limits_{0}^{t}{\int\limits_{0}^{a}{f\left( \eta ,s \right)G\left( t,x,\eta ,s \right)dsd\eta }}=\int\limits_{0}^{t}\int\limits_{0}^{a}f\left( \eta ,s \right)\times$$
$$\times\sum\limits_{n=-\infty }^{+\infty }{\sum\limits_{k=0}^{+\infty }{\frac{{{\left( -1 \right)}^{k}}\left[ {{\left| x-s+2na \right|}^{k}}-{{\left| x+s+2na \right|}^{k}} \right]}{2\Gamma(k+1)}}}\times$$
\begin{equation}\label{eq20}
    \times(t-\eta)^{-\beta_{1}k-\beta_1-1}E_{\alpha, {-\beta_1-\beta_1k}}^{-\gamma _{1}-{{\gamma }_{1}}k}\left[\delta (t-\eta)^\alpha\right]dsd\eta.
\end{equation}
By subtracting equation \eqref{eq20} from equation \eqref{eq19}, it follows that $Lu_1=f(t,x).$

Next, we verify that the solution \eqref{eq2.4} satisfies the boundary conditions \eqref{eq2.2}. For this purpose, we let $x \to 0$ in \eqref{eq2.4}. From \eqref{eq2.5}, it follows that 
$$\lim _{x \to 0}{G}\left( t,x,\eta ,s \right)=0, \,\,\,\,\,\,\lim _{x \to 0}{G}_{s}\left( t,x,\eta ,a \right)=0,\,\,\,\,\,\lim _{x \to 0} G\left( t,x,0,s\right)=0.$$

As a result, we get
$$\lim _{x \to 0}u\left( t,x \right)=\lim _{x \to 0}\int\limits_{0}^{t}{{{\varphi }_{0}}\left( \eta  \right){{G}_{s}}\left( t,x,\eta ,0 \right)d\eta }.$$
Now, let us evaluate this equality directly:
$$\lim _{x \to 0}\int\limits_{0}^{t}{{{\varphi }_{0}}\left( \eta  \right){{G}_{s}}\left( t,x,\eta ,0 \right)d\eta }=\lim _{x \to 0}\int\limits_{0}^{t}{{{\varphi }_{0}}\left( \eta  \right)\sum\limits_{n=-\infty }^{\infty }{sign\left( x+2na \right)\omega \left( t-\eta ,\left| x+2na \right| \right)\,}d\eta }=$$
$$=\lim _{x \to 0}\int\limits_{0}^{t}{{{\varphi }_{0}}\left( \eta  \right)\omega \left( t-\eta ,x \right)d\eta}-\lim _{x \to 0}\int\limits_{0}^{t}{{{\varphi }_{0}}\left( \eta  \right)\sum\limits_{n=1}^{\infty }{\left[ \omega \left( t-\eta ,2na-x \right)-\omega \left( t-\eta ,2na+x \right) \right]\,}d\eta }=$$
$$=\lim _{x \to 0}\int\limits_{0}^{t}{{{\varphi }_{0}}\left( \eta  \right)\omega \left( t-\eta ,x \right)d\eta}.$$

We replace the function $\omega$ with its value and apply the substitution 
$s=\frac{x}{(t-\eta)^{\beta_1}}$:
$$\lim _{x \to 0}u\left( t,x \right)=\lim _{x \to 0}\int\limits_{0}^{t}{{{\varphi }_{0}}\left( \eta  \right)\omega \left( t-\eta ,x \right)d\eta}=$$
$$=\lim _{x \to 0}\int\limits_{0}^{t}{{{\varphi }_{0}}\left( \eta  \right)\sum\limits_{n=0}^{+\infty }{\frac{{{\left( -1 \right)}^{n}}{{x}^{n}}}{n!}{{(t-\eta)}^{-{{\beta }_{1}}n-1}}E_{\alpha ,-{{\beta }_{1}}n}^{-{{\gamma }_{1}}n}\left[ \delta {{(t-\eta)}^{\alpha }} \right]}d\eta}=$$
$$=\lim _{x \to 0}\int\limits_{0}^{+\infty}{{{\varphi }_{0}}\left[ t-\left(\frac{x}{s}\right)^{\frac{1}{\beta_1}}  \right]\sum\limits_{n=0}^{+\infty }{\frac{{{\left( -1 \right)}^{n}}{{s}^{n-1}}}{\beta_1 n!}E_{\alpha ,-{{\beta }_{1}}n}^{-{{\gamma }_{1}}n}\left[ \delta \left(\frac{x}{s}\right)^{\frac{\alpha}{\beta_1}} \right]}ds}=$$
$$=\varphi_0(t)\int\limits_{0}^{+\infty}\sum\limits_{n=0}^{+\infty }{\frac{{{\left( -1 \right)}^{n}}{{s}^{n-1}}}{ n!\beta_1 \Gamma(-{{\beta }_{1}}n)}}ds=-\varphi_0(t)\left.\sum\limits_{n=0}^{+\infty }{\frac{{{\left( -1 \right)}^{n}}{{s}^{n}}}{ \Gamma(n+1)\Gamma(1-{{\beta }_{1}}n)}}\right|_0^{+\infty}=$$
$$=-\varphi_0(t)\lim_{s \to +\infty}e_{1,\beta_1}^{1,1}(-s)+\varphi_0(t)=\varphi_0(t).$$

Similarly, we have
$$\lim _{x \to a}u\left( t,x \right)=\varphi_1(t).$$

We now check that our solution \eqref{eq2.4} fulfills the initial condition \eqref{eq2.3}:
$$\lim_{t \to 0} {}^{P}I_{0t}^{\alpha ,\,1-\beta ,\,-\gamma ,\,\delta }u\left( t,x \right)=\lim_{t \to 0}{}^{P}I_{0t}^{\alpha ,\,1-\beta ,\,-\gamma ,\,\delta }\int\limits_{0}^{a}{\tau \left( s \right) G\left( t,x,0,s \right)ds}=$$
$$=\lim_{t \to 0}{}^{P}I_{0t}^{\alpha ,\,1-\beta ,\,-\gamma ,\,\delta }\int\limits_{0}^{a}{\tau \left( s \right) \sum\limits_{k=0}^{+\infty }{\frac{{{\left( -1 \right)}^{k}} {{\left| x-s \right|}^{k}}}{2\Gamma(k+1)}}}t^{-\beta_1k+{{\beta }_{1}}-1}E_{\alpha, \beta_1-\beta_1k}^{{\gamma }_{1}-{{\gamma }_{1}}k}\left[\delta t^\alpha\right]ds+$$
$$+\lim_{t \to 0}{}^{P}I_{0t}^{\alpha ,\,1-\beta ,\,-\gamma ,\,\delta }\int\limits_{0}^{a}{\tau \left( s \right) K\left( t,x,0,s \right)ds}=$$
$$=\lim_{t \to 0}{}^{P}I_{0t}^{\alpha ,\,1-\beta ,\,-\gamma ,\,\delta }\int\limits_{0}^{x}{\tau \left( s \right) \sum\limits_{k=0}^{+\infty }{\frac{{{\left( -1 \right)}^{k}} {{\left( x-s \right)}^{k}}}{2\Gamma(k+1)}}}t^{-\beta_1k+{{\beta }_{1}}-1}E_{\alpha, \beta_1-\beta_1k}^{{\gamma }_{1}-{{\gamma }_{1}}k}\left[\delta t^\alpha\right]ds+$$
$$+\lim_{t \to 0}{}^{P}I_{0t}^{\alpha ,\,1-\beta ,\,-\gamma ,\,\delta }\int\limits_{x}^{a}{\tau \left( s \right) \sum\limits_{k=0}^{+\infty }{\frac{{{\left( -1 \right)}^{k}} {{\left( s-x \right)}^{k}}}{2\Gamma(k+1)}}}t^{-\beta_1k+{{\beta }_{1}}-1}E_{\alpha, \beta_1-\beta_1k}^{{\gamma }_{1}-{{\gamma }_{1}}k}\left[\delta t^\alpha\right]ds+$$
$$+\lim_{t \to 0}{}^{P}I_{0t}^{\alpha ,\,1-\beta ,\,-\gamma ,\,\delta }\int\limits_{0}^{a}{\tau \left( s \right) K\left( t,x,0,s \right)ds}.$$
From the last expression, we evaluate the term containing the integral with respect to $s$ from $0$ to $x:$
$$\lim_{t \to 0}{}^{P}I_{0t}^{\alpha ,\,1-\beta ,\,-\gamma ,\,\delta }\int\limits_{0}^{x}{\tau \left( s \right) \sum\limits_{k=0}^{+\infty }{\frac{{{\left( -1 \right)}^{k}} {{\left( x-s \right)}^{k}}}{2\Gamma(k+1)}}}t^{-\beta_1k+{{\beta }_{1}}-1}E_{\alpha, \beta_1-\beta_1k}^{{\gamma }_{1}-{{\gamma }_{1}}k}\left[\delta t^\alpha\right]ds=$$
$$=\lim_{t \to 0}\int\limits_{0}^{t}{{{\left( t-y \right)}^{-\beta }}}E_{\alpha ,\,1-\beta }^{-\gamma }\left[ \delta {{\left( t-y \right)}^{\alpha }} \right]dy\int\limits_{0}^{x}{\tau \left( s \right) \sum\limits_{k=0}^{+\infty }{\frac{{{\left( -1 \right)}^{k}} {{\left( x-s \right)}^{k}}}{2\Gamma(k+1)}}}y^{-\beta_1k+{{\beta }_{1}}-1}E_{\alpha, \beta_1-\beta_1k}^{{\gamma }_{1}-{{\gamma }_{1}}k}\left[\delta y^\alpha\right]ds=$$
$$=\lim_{t \to 0}\int\limits_{0}^{x}{\tau \left( s \right) \sum\limits_{k=0}^{+\infty }{\frac{{{\left( -1 \right)}^{k}} {{\left( x-s \right)}^{k}}}{2\Gamma(k+1)}}}ds\int\limits_{0}^{t}y^{-\beta_1k+{{\beta }_{1}}-1}E_{\alpha, \beta_1-\beta_1k}^{{\gamma }_{1}-{{\gamma }_{1}}k}\left[\delta y^\alpha\right]{{{\left( t-y \right)}^{-\beta }}}E_{\alpha ,\,1-\beta }^{-\gamma }\left[ \delta {{\left( t-y \right)}^{\alpha }} \right]dy.$$
After computing the integral with respect to $y$, we obtain 
$$\int\limits_{0}^{t}y^{-\beta_1k+{{\beta }_{1}}-1}E_{\alpha, \beta_1-\beta_1k}^{{\gamma }_{1}-{{\gamma }_{1}}k}\left[\delta y^\alpha\right]{{{\left( t-y \right)}^{-\beta }}}E_{\alpha ,\,1-\beta }^{-\gamma }\left[ \delta {{\left( t-y \right)}^{\alpha }} \right]dy=t^{-\beta_1-\beta_1k}E_{\alpha,1-\beta_1-\beta_1k}^{-\gamma_1-\gamma_1k}[\delta t^\alpha].$$
Next, we proceed to evaluate the integral with respect to $s.$ We make the substitution $\frac{x-s}{t^{\beta_1}}=z:$
$$\lim_{t \to 0}\int\limits_{0}^{x}{\tau \left( s \right) \sum\limits_{k=0}^{+\infty }{\frac{{{\left( -1 \right)}^{k}} }{2\Gamma(k+1)}}}\left(\frac{x-s}{t^{\beta_1}}\right)^{k}t^{-\beta_1}E_{\alpha,1-\beta_1-\beta_1k}^{-\gamma_1-\gamma_1k}[\delta t^\alpha]ds=$$
$$=\frac{1}{2}\lim_{t \to 0}\int\limits_{0}^{\frac{x}{t^{\beta_1}}}{\tau \left( x-t^{\beta_1} z\right) \sum\limits_{k=0}^{+\infty }{\frac{{{\left( -1 \right)}^{k}z^{k}} }{\Gamma(k+1)}}}E_{\alpha,1-\beta_1-\beta_1k}^{-\gamma_1-\gamma_1k}[\delta t^\alpha]dz.$$
Now we compute the limit and use \eqref{eq6}:
$$\lim_{t \to 0}{}^{P}I_{0t}^{\alpha ,\,1-\beta ,\,-\gamma ,\,\delta }\int\limits_{0}^{x}{\tau \left( s \right) \sum\limits_{k=0}^{+\infty }{\frac{{{\left( -1 \right)}^{k}} {{\left( x-s \right)}^{k}}}{2\Gamma(k+1)}}}t^{-\beta_1k+{{\beta }_{1}}-1}E_{\alpha, \beta_1-\beta_1k}^{{\gamma }_{1}-{{\gamma }_{1}}k}\left[\delta t^\alpha\right]ds=$$
$$=\frac{\tau(x)}{2}\int\limits_{0}^{+\infty}{ \sum\limits_{k=0}^{+\infty }{\frac{{{\left( -1 \right)}^{k}z^{k}} }{\Gamma(k+1)\Gamma(1-\beta_1-\beta_1k)}}}dz=$$
\begin{equation}\label{eq25}
    =-\frac{\tau(x)}{2}\left.(-z)e_{1,\beta_1}^{2,1-\beta_1}(-z)\right|_0^{+\infty}=\frac{\tau(x)}{2}.
\end{equation}
Similarly,
\begin{equation}\label{eq26}
    \lim_{t \to 0}{}^{P}I_{0t}^{\alpha ,\,1-\beta ,\,-\gamma ,\,\delta }\int\limits_{x}^{a}{\tau \left( s \right) \sum\limits_{k=0}^{+\infty }{\frac{{{\left( -1 \right)}^{k}} {{\left( s-x \right)}^{k}}}{2\Gamma(k+1)}}}t^{-\beta_1k+{{\beta }_{1}}-1}E_{\alpha, \beta_1-\beta_1k}^{{\gamma }_{1}-{{\gamma }_{1}}k}\left[\delta t^\alpha\right]ds=\frac{\tau(x)}{2}.
\end{equation}
Now we turn to the following term:
$$\lim_{t \to 0}{}^{P}I_{0t}^{\alpha ,\,1-\beta ,\,-\gamma ,\,\delta }\int\limits_{0}^{a}{\tau \left( s \right) K\left( t,x,0,s \right)ds}=$$
$$=\lim_{t \to 0}\int\limits_{0}^{t}{{{\left( t-y \right)}^{-\beta }}}E_{\alpha ,\,1-\beta }^{-\gamma }\left[ \delta {{\left( t-y \right)}^{\alpha }} \right]dy\times$$
$$\times\int\limits_{0}^{a}\tau \left( s \right)\sum\limits_{n=-\infty }^{+\infty }{}'{\sum\limits_{k=0}^{+\infty }{\frac{{{\left( -1 \right)}^{k}}{\left| x-s+2na \right|}^{k}}{2\Gamma(k+1)}}}y^{-\beta_1k+{{\beta }_{1}}-1}E_{\alpha, \beta_1-\beta_1k}^{{\gamma }_{1}-{{\gamma }_{1}}k}\left[\delta y^\alpha\right]ds-$$
$$-\lim_{t \to 0}\int\limits_{0}^{t}{{{\left( t-y \right)}^{-\beta }}}E_{\alpha ,\,1-\beta }^{-\gamma }\left[ \delta {{\left( t-y \right)}^{\alpha }} \right]dy\times$$
$$\times\int\limits_{0}^{a}\tau \left( s \right)\sum\limits_{n=-\infty }^{+\infty }{\sum\limits_{k=0}^{+\infty }{\frac{{{\left( -1 \right)}^{k}}{\left| x+s+2na \right|}^{k}}{2\Gamma(k+1)}}}y^{-\beta_1k+{{\beta }_{1}}-1}E_{\alpha, \beta_1-\beta_1k}^{{\gamma }_{1}-{{\gamma }_{1}}k}\left[\delta y^\alpha\right]ds=$$
$$=\lim_{t \to 0}\int\limits_{0}^{a}{\tau \left( s \right) \sum\limits_{n=-\infty }^{+\infty }{}'{\sum\limits_{k=0}^{+\infty }{\frac{{{\left( -1 \right)}^{k}}{\left| x-s+2na \right|}^{k}}{2\Gamma(k+1)}}}}t^{-\beta_1-\beta_1k}E_{\alpha,1-\beta_1-\beta_1k}^{-\gamma_1-\gamma_1k}[\delta t^\alpha]ds-$$
$$-\lim_{t \to 0}\int\limits_{0}^{a}{\tau \left( s \right) \sum\limits_{n=-\infty }^{+\infty }{\sum\limits_{k=0}^{+\infty }{\frac{{{\left( -1 \right)}^{k}}{\left| x+s+2na \right|}^{k}}{2\Gamma(k+1)}}}}t^{-\beta_1-\beta_1k}E_{\alpha,1-\beta_1-\beta_1k}^{-\gamma_1-\gamma_1k}[\delta t^\alpha]ds.$$
When evaluating the limits in the last equality, by using the Lemma 3.1, we see that their values are equal to 0.

Consequently, according to \eqref{eq25} and \eqref{eq26}, 
$$\underset{t\to 0}{\mathop{\lim }}\,{}^{P}I_{0t}^{\alpha ,\,1-\beta ,\,-\gamma ,\,\delta }u\left( t,x \right)=\tau \left( x \right),\,\,\,  0\le x\le a.	$$

\subsection{\bf Second appendix} \label{secA2}
 
Let us prove the equality \eqref{eq2.6}. To proceed, we express the left-hand side of \eqref{eq2.6} in the following form:
$$\underset{t\to 0}{\mathop{\lim }}\,I_{0t}^{\alpha ,\,1-\frac{\beta }{2},\,-\frac{\gamma }{2},\,\delta }u\left( t,x \right)=\underset{t\to 0}{\mathop{\lim }}\,\int\limits_{0}^{t}{{{\left( t-z \right)}^{-\frac{\beta }{2}}}}{{z}^{\beta -1}}E_{\alpha ,1-\frac{\beta }{2}}^{-\frac{\gamma }{2}}\left[ \delta {{\left( t-z \right)}^{\alpha }} \right]{{z}^{1-\beta }}u\left( z,x \right)dz.$$
From the fact that ${{t}^{1-\beta }}u\left( t,x \right)\in C\left( \overline{D} \right),$ it follows that the expression ${{t}^{1-\beta }}u\left( t,x \right)$ is bounded.
$$\underset{t\to 0}{\mathop{\lim }}\,\left|\int\limits_{0}^{t}{{{\left( t-z \right)}^{-\frac{\beta }{2}}}}{{z}^{\beta -1}}E_{\alpha ,1-\frac{\beta }{2}}^{-\frac{\gamma }{2}}\left[ \delta {{\left( t-z \right)}^{\alpha }} \right]{{z}^{1-\beta }}u\left( z,x \right)dz\right|\le $$
$$\le \left\| {{z}^{1-\beta }}u\left( z,x \right) \right\|\underset{t\to 0}{\mathop{\lim }}\,\left|\int\limits_{0}^{t}{{{\left( t-z \right)}^{-\frac{\beta }{2}}}}{{z}^{\beta -1}}E_{\alpha ,1-\frac{\beta }{2}}^{-\frac{\gamma }{2}}\left[ \delta {{\left( t-z \right)}^{\alpha }} \right]dz\right|=$$
$$=\left\| {{z}^{1-\beta }}u\left( z,x \right) \right\|\underset{t\to 0}{\mathop{\lim }}\left|\sum\limits_{k=0}^{+\infty }{\frac{{{\left( -\frac{\gamma }{2} \right)}_{k}}{{\delta }^{k}}}{k!\Gamma \left( \alpha k+1-\frac{\beta }{2} \right)}}\,\int\limits_{0}^{t}{{{\left( t-z \right)}^{\alpha k-\frac{\beta }{2}}}}{{z}^{\beta -1}}dz\right|.$$
Next, we perform the substitution $z=ts:$
$$\left\| {{z}^{1-\beta }}u\left( z,x \right) \right\|\underset{t\to 0}{\mathop{\lim }}\left|\sum\limits_{k=0}^{+\infty }{\frac{{{\left( -\frac{\gamma }{2} \right)}_{k}}{{\delta }^{k}}}{k!\Gamma \left( \alpha k+1-\frac{\beta }{2} \right)}}\,{{t}^{\alpha k+\frac{\beta }{2}}}\int\limits_{0}^{1}{{{\left( 1-s \right)}^{\alpha k-\frac{\beta }{2}}}}{{s}^{\beta -1}}ds\right|=$$
$$=\left\| {{z}^{1-\beta }}u\left( z,x \right) \right\|\underset{t\to 0}{\mathop{\lim }}\left|{t}^{\frac{\beta }{2}}\Gamma \left( \beta  \right)\sum\limits_{k=0}^{+\infty }{\frac{{{\left( -\frac{\gamma }{2} \right)}_{k}}{{\delta }^{k}}{{t}^{\alpha k}}}{k!\Gamma \left( \alpha k+1+\frac{\beta }{2} \right)}}\,\right|=0.$$
Hence, \eqref{eq2.6} is proved.

\subsection{\bf Third appendix} \label{secA3}

Let us first prove the following formula:
\begin{equation}\label{eq1}
    {}^{P}I_{0 t}^{\alpha, \beta, \gamma, \delta}f(t)=\frac{d}{dt}{}^{P}I_{0 t}^{\alpha, 1+\beta, \gamma, \delta}f(t).
\end{equation}
We expand the Prabhakar integral on the right-hand side of the equality. Then, we differentiate the integral with respect to the parameter $t:$
$$\frac{d}{dt}{}^{P}I_{0 t}^{\alpha, 1+\beta, \gamma, \delta}f(t)=\frac{d}{dt}\int\limits_{0}^{t}(t-s)^{\beta}E_{\alpha, 1+\beta}^{\gamma}\left[\delta(t-s)^\alpha\right]f(s)ds=$$
$$=\lim_{s\to t}\left[(t-s)^{\beta}E_{\alpha, 1+\beta}^{\gamma}\left[\delta(t-s)^\alpha\right]f(s)\right]+$$
$$+\int\limits_{0}^{t}\frac{\partial}{\partial t}\left[(t-s)^{\beta}E_{\alpha, 1+\beta}^{\gamma}\left[\delta(t-s)^\alpha\right]f(s)\right]ds=$$
$$=\int\limits_{0}^{t}(t-s)^{\beta-1}E_{\alpha, \beta}^{\gamma}\left[\delta(t-s)^\alpha\right]f(s)ds={}^{P}I_{0 t}^{\alpha, \beta, \gamma, \delta}f(t).$$

Now we prove that the following semigroup property holds for the fractional order Prabhakar derivative:
\begin{equation}\label{eq2}
    {}^{PRL}D_{0t}^{\alpha ,\,\frac{\beta}{2} ,\,\frac{\gamma}{2} ,\,\delta }\left({}^{PRL}D_{0t}^{\alpha ,\,\frac{\beta}{2} ,\,\frac{\gamma}{2} ,\,\delta }u\left( t,x \right)\right)={}^{PRL}D_{0t}^{\alpha ,\beta,\,\gamma ,\,\delta }u(t,x).
\end{equation}
For the left-hand side of \eqref{eq2}, the condition  
${}^{PRL}D_{0t}^{\alpha,\frac{\beta}{2},\frac{\gamma}{2},\delta}u(t,x)\in C(D)$
needs to be imposed. However, this follows from the condition stated in Definition 2.1, that is,
${}^{PRL}D_{0t}^{\alpha,\beta,\gamma,\delta}u(t,x)\in C(D).$

We express the left-hand side of \eqref{eq2} in terms of the Prabhakar integral, and for convenience, introduce the notation ${}^{P}I_{0 t}^{\alpha, 1-\frac{\beta}{2},-\frac{\gamma}{2}, \delta}u(t,x)=\varphi(t,x):$
$${}^{PRL}D_{0t}^{\alpha ,\,\frac{\beta}{2} ,\,\frac{\gamma}{2} ,\,\delta }\left({}^{PRL}D_{0t}^{\alpha ,\,\frac{\beta}{2} ,\,\frac{\gamma}{2} ,\,\delta }u\left( t,x \right)\right)=$$
$$={\frac{\partial}{\partial t}}{}^{P}I_{0 t}^{\alpha, 1-\frac{\beta}{2},-\frac{\gamma}{2}, \delta}\left({\frac{\partial}{\partial t}}{}^{P}I_{0 t}^{\alpha, 1-\frac{\beta}{2},-\frac{\gamma}{2}, \delta}u(t,x)\right)=$$
$$={\frac{\partial}{\partial t}}{}^{P}I_{0 t}^{\alpha, 1-\frac{\beta}{2},-\frac{\gamma}{2}, \delta}\left({\frac{\partial\varphi(t,x)}{\partial t}}\right).$$
Now we apply the formula \eqref{eq1}:
$${\frac{\partial^2}{\partial t^2}}{}^{P}I_{0 t}^{\alpha, 2-\frac{\beta}{2},-\frac{\gamma}{2}, \delta}\left({\frac{\partial\varphi(t,x)}{\partial t}}\right)=$$
$$={\frac{\partial^2}{\partial t^2}}\int\limits_{0}^{t}(t-s)^{1-\frac{\beta}{2}}E_{\alpha, 2-\frac{\beta}{2}}^{-\frac{\gamma}{2}}\left[\delta(t-s)^\alpha\right]{\frac{\partial\varphi(s,x)}{\partial s}}ds.$$

Now we integrate the integral by parts:
$$\frac{\partial^2}{\partial t^2}\left[\left.(t-s)^{1-\frac{\beta}{2}}E_{\alpha, 2-\frac{\beta}{2}}^{-\frac{\gamma}{2}}\left[\delta(t-s)^\alpha\right]\varphi(s,x)\right|_0^t+\right.$$
$$\left.+\int\limits_{0}^{t}(t-s)^{-\frac{\beta}{2}}E_{\alpha, 1-\frac{\beta}{2}}^{-\frac{\gamma}{2}}\left[\delta(t-s)^\alpha\right]\varphi(s,x)ds\right]=$$
$$=\frac{\partial^2}{\partial t^2}\left[\lim_{s\to t}\left((t-s)^{1-\frac{\beta}{2}}E_{\alpha, 2-\frac{\beta}{2}}^{-\frac{\gamma}{2}}\left[\delta(t-s)^\alpha\right]\varphi(s,x)\right)-\right.$$
$$-t^{1-\frac{\beta}{2}}E_{\alpha, 2-\frac{\beta}{2}}^{-\frac{\gamma}{2}}\left[\delta t^\alpha\right]\lim_{s\to 0}\varphi(s,x)+$$
$$\left.+\int\limits_{0}^{t}(t-s)^{-\frac{\beta}{2}}E_{\alpha, 1-\frac{\beta}{2}}^{-\frac{\gamma}{2}}\left[\delta(t-s)^\alpha\right]\varphi(s,x)ds\right].$$
Since $1-\frac{\beta}{2}>0,$
$$\lim_{s\to t}\left((t-s)^{1-\frac{\beta}{2}}E_{\alpha, 2-\frac{\beta}{2}}^{-\frac{\gamma}{2}}\left[\delta(t-s)^\alpha\right]\varphi(s,x)\right)=0.$$
According to \eqref{eq2.6}, we obtain the following
$$\lim_{s\to 0}\varphi(s,x)=\lim_{s\to 0}{}^{P}I_{0 s}^{\alpha, 1-\frac{\beta}{2},-\frac{\gamma}{2}, \delta}u(s,x)=0.$$
The above equality holds due to the fact that $t^{1-\beta}u(t,x)\in C(\overline{D}).$ This condition is included in Definition 2.1.

Since the above limits are equal to zero, we obtain the following:
$${}^{PRL}D_{0t}^{\alpha ,\,\frac{\beta}{2} ,\,\frac{\gamma}{2} ,\,\delta }\left({}^{PRL}D_{0t}^{\alpha ,\,\frac{\beta}{2} ,\,\frac{\gamma}{2} ,\,\delta }u\left( t,x \right)\right)=$$
$$=\frac{\partial^2}{\partial t^2}\int\limits_{0}^{t}(t-s)^{-\frac{\beta}{2}}E_{\alpha, 1-\frac{\beta}{2}}^{-\frac{\gamma}{2}}\left[\delta(t-s)^\alpha\right]\varphi(s,x)ds.$$
Next, we replace the function $\varphi(s,x)$ with its integral form and change the order of integration:
$$\frac{\partial^2}{\partial t^2}\int\limits_{0}^{t}(t-s)^{-\frac{\beta}{2}}E_{\alpha, 1-\frac{\beta}{2}}^{-\frac{\gamma}{2}}\left[\delta(t-s)^\alpha\right]ds\int\limits_{0}^{s}(s-z)^{-\frac{\beta}{2}}E_{\alpha, 1-\frac{\beta}{2}}^{-\frac{\gamma}{2}}\left[\delta(s-z)^\alpha\right]u(z,x)dz=$$
$$=\frac{\partial^2}{\partial t^2}\int\limits_{0}^{t}u(z,x)dz\int\limits_{z}^{t}(t-s)^{-\frac{\beta}{2}}(s-z)^{-\frac{\beta}{2}}E_{\alpha, 1-\frac{\beta}{2}}^{-\frac{\gamma}{2}}\left[\delta(t-s)^\alpha\right]E_{\alpha, 1-\frac{\beta}{2}}^{-\frac{\gamma}{2}}\left[\delta(s-z)^\alpha\right]ds=$$
$$=\frac{\partial^2}{\partial t^2}\int\limits_{0}^{t}u(z,x)dz\int\limits_{z}^{t}(t-s)^{-\frac{\beta}{2}}(s-z)^{-\frac{\beta}{2}}\times$$
$$\times \sum_{k=0}^{+\infty}\frac{\left(-\frac{\gamma}{2}\right)_k \delta^k (t-s)^{\alpha k}}{k!\Gamma\left(\alpha k+1-\frac{\beta}{2}\right)} \sum_{m=0}^{+\infty}\frac{\left(-\frac{\gamma}{2}\right)_m \delta^m (s-z)^{\alpha m}}{m!\Gamma\left(\alpha m+1-\frac{\beta}{2}\right)}ds.$$
Next, we apply the formula \eqref{eq2.27}:
$$\sum_{k=0}^{+\infty}\sum_{m=0}^{k}\frac{\delta^k\left(-\frac{\gamma}{2}\right)_m \left(-\frac{\gamma}{2}\right)_{k-m}  }{m!(k-m)!\Gamma\left(\alpha m+1-\frac{\beta}{2}\right)\Gamma\left(\alpha k-\alpha m+1-\frac{\beta}{2}\right)}\times$$
$$\times\frac{\partial^2}{\partial t^2}\int\limits_{0}^{t}u(z,x)dz\int\limits_{z}^{t}(t-s)^{\alpha m-\frac{\beta}{2}}(s-z)^{\alpha k-\alpha m-\frac{\beta}{2}}ds.$$
Now, in the integral with respect to $s$, we make the substitution $s=(t-z)y+z$ and use \eqref{eq2.28}:
$$\sum_{k=0}^{+\infty}\frac{\delta^k}{\Gamma\left(\alpha k+2-\beta\right)}\sum_{m=0}^{k}\frac{\left(-\frac{\gamma}{2}\right)_m \left(-\frac{\gamma}{2}\right)_{k-m}  }{m!(k-m)!}\frac{\partial^2}{\partial t^2}\int\limits_{0}^{t}(t-z)^{\alpha k+1-\beta}u(z,x)dz=$$
$$=\sum_{k=0}^{+\infty}\frac{ (-\gamma)_k \delta^k}{k! \Gamma\left(\alpha k+2-\beta\right)}\frac{\partial^2}{\partial t^2}\int\limits_{0}^{t}(t-z)^{\alpha k+1-\beta}u(z,x)dz=$$
$$=\frac{\partial^2}{\partial t^2}\int\limits_{0}^{t}(t-z)^{\alpha k+1-\beta}E_{\alpha, 2-\beta}^{-\gamma}\left[\delta(t-z)^\alpha\right]u(z,x)dz=$$
$$=\frac{\partial^2}{\partial t^2}{}^{P}I_{0 t}^{\alpha, 2-\beta, -\gamma, \delta}u(t,x)=\frac{\partial}{\partial t}{}^{P}I_{0 t}^{\alpha, 1-\beta, -\gamma, \delta}u(t,x)=$$
$$={}^{PRL}D_{0t}^{\alpha ,\beta,\,\gamma ,\,\delta }u(t,x).$$
It has been proved that property \eqref{eq2} holds.

Let us expand the left-hand side of the following formula term by term:
$$\left( {}^{PRL}D_{0t}^{\alpha ,\,\frac{\beta }{2},\,\frac{\gamma }{2},\,\delta }-\frac{\partial }{\partial x} \right)\left( {}^{PRL}D_{0t}^{\alpha ,\,\frac{\beta }{2},\,\frac{\gamma }{2},\,\delta }+\frac{\partial }{\partial x} \right)u\left( t,x \right)=$$
$$={}^{PRL}D_{0t}^{\alpha ,\,\frac{\beta }{2},\,\frac{\gamma }{2},\,\delta }\left({}^{PRL}D_{0t}^{\alpha ,\,\frac{\beta }{2},\,\frac{\gamma }{2},\,\delta }u(t,x)\right)+{}^{PRL}D_{0t}^{\alpha ,\,\frac{\beta }{2},\,\frac{\gamma }{2},\,\delta }\left(\frac{\partial}{\partial x}u(t,x)\right)-$$
$$-\frac{\partial}{\partial x}\left({}^{PRL}D_{0t}^{\alpha ,\,\frac{\beta }{2},\,\frac{\gamma }{2},\,\delta }u(t,x)\right)-\frac{\partial^2}{\partial x^2}u(t,x).$$

The first and fourth terms on the right-hand side of the last equality are clear. Now, we show the existence of the second and third terms and that they are equal to each other.
Regarding the second term, it is required that $u_x(t,x)\in C(D)$, which follows from the condition $u_{xx}(t,x)\in C(D)$ stated in Definition 2.1.

Let us consider the following difference:
$${}^{PRL}D_{0t}^{\alpha ,\,\frac{\beta }{2},\,\frac{\gamma }{2},\,\delta }\left(\frac{\partial}{\partial x}u(t,x)\right)-\frac{\partial}{\partial x}\left({}^{PRL}D_{0t}^{\alpha ,\,\frac{\beta }{2},\,\frac{\gamma }{2},\,\delta }u(t,x)\right).$$
Since the Prabhakar fractional derivative is taken with respect to the variable $t$, while the derivative is with respect to $x$, in the second term we can move the derivative with respect to $x$ inside. As a result, we can write
$${}^{PRL}D_{0t}^{\alpha ,\,\frac{\beta }{2},\,\frac{\gamma }{2},\,\delta }\left(\frac{\partial}{\partial x}u(t,x)\right)=\frac{\partial}{\partial x}\left({}^{PRL}D_{0t}^{\alpha ,\,\frac{\beta }{2},\,\frac{\gamma }{2},\,\delta }u(t,x)\right)$$
and the difference under consideration turns out to be zero.

\subsection{\bf Fourth appendix} \label{secA4}

Now we solve the problem \eqref{eq2.9}-\eqref{eq2.10}-\eqref{eq2.11}. Let $V\left( t,x;\eta ,\xi  \right)$ be a function of the variables $\eta $ and $\xi$ in a domain $\left\{ \left( \eta ,\xi  \right):\,\,0<\eta <t,\,\,0<\xi <x \right\},$ which satisfies the following equation for any fixed $\left( t,x \right):$
\begin{equation}\label{eq3.1}
    {}^{PC}D_{\eta t}^{\alpha ,{{\beta }_{1}},{{\gamma }_{1}},\delta }V\left( t,x;\eta ,\xi  \right)-{{V}_{\xi }}\left( t,x;\eta ,\xi  \right)=1,		
\end{equation}
also satisfies these conditions
\begin{equation}\label{eq3.2}
    {{\left. V\left( t,x;\eta ,\xi  \right) \right|}_{\eta =t}}=0,   \,\,\,{{\left. V\left( t,x;\eta ,\xi  \right) \right|}_{\xi =x}}=0,
\end{equation}
where ${ }^{P C} D_{0 t}^{\alpha, \beta, \gamma, \delta} y(t)={ }^P I_{0 t}^{\alpha, 1-\beta,-\gamma, \delta} \frac{d}{d t} y(t),$ $y \in A C^m[0, T]$ 
is the Caputo-Prabhakar fractional derivative \cite{Garra}.

First, we rewrite the equation \eqref{eq3.1} by replacing the variables $x$ and $t$ with $\xi $ and $\eta :$
$${}^{PRL}D_{0\eta }^{\alpha ,\,{{\beta }_{1}},\,{{\gamma }_{1}},\,\delta }u\left( \eta ,\xi  \right)+{{u}_{\xi }}\left( \eta ,\xi  \right)=v\left( \eta ,\xi  \right).$$
Then, we multiply both sides of the last equation by $V\left( t,x;\eta ,\xi  \right),$ integrate with respect to $\xi $ over the interval $\left[ 0,x \right],$ and with respect to $\eta $ over the interval $\left[ 0,t \right]:$
$$\int\limits_{0}^{t}{\int\limits_{0}^{x}{{}^{PRL}D_{0\eta }^{\alpha ,\,{{\beta }_{1}},\,{{\gamma }_{1}},\,\delta }u\left( \eta ,\xi  \right)V\left( t,x;\eta ,\xi  \right)d\xi d\eta }}+\int\limits_{0}^{t}{\int\limits_{0}^{x}{{{u}_{\xi }}\left( \eta ,\xi  \right)V\left( t,x;\eta ,\xi  \right)d\xi }d\eta }=$$
\begin{equation}\label{eq3.3}
    =\int\limits_{0}^{t}{\int\limits_{0}^{x}{v\left( \eta ,\xi  \right)V\left( t,x;\eta ,\xi  \right)d\xi d\eta }}.
\end{equation}
For convenience, let us denote the first integral in \eqref{eq3.3} by ${{I}_{1}}$:
$${{I}_{1}}=\int\limits_{0}^{t}{\int\limits_{0}^{x}{V\left( t,x;\eta ,\xi  \right){}^{PRL}D_{0\eta }^{\alpha ,\,{{\beta }_{1}},\,{{\gamma }_{1}},\,\delta }u\left( \eta ,\xi  \right)d\xi d\eta }}=\int\limits_{0}^{t}{\int\limits_{0}^{x}{V\left( t,x;\eta ,\xi  \right)\frac{\partial }{\partial \eta }I_{0\eta }^{\alpha ,\,1-{{\beta }_{1}},\,-{{\gamma }_{1}},\,\delta }u\left( \eta ,\xi  \right)d\xi d\eta }}.$$
Applying integration by parts with respect to $\eta $ and considering the conditions \eqref{eq2.10} and \eqref{eq3.2}, we obtain the following result:
$${{I}_{1}}=-\int\limits_{0}^{x}{d\xi }\int\limits_{0}^{t}{{{V}_{\eta }}\left( t,x;\eta ,\xi  \right)I_{0\eta }^{\alpha ,\,1-\beta ,\,-\gamma ,\,\delta }u\left( \eta ,\xi  \right)d\eta }.$$
Now, we express the Prabhakar fractional integral and change the order of integration:
$${{I}_{1}}=-\int\limits_{0}^{x}{d\xi }\int\limits_{0}^{t}{{{V}_{\eta }}\left( t,x;\eta ,\xi  \right)d\eta }\int\limits_{0}^{\eta }{{{\left( \eta -s \right)}^{-{{\beta }_{1}}}}E_{\alpha ,\,1-{{\beta }_{1}}}^{-{{\gamma }_{1}}}\left[ \delta {{\left( \eta -s \right)}^{\alpha }} \right]}\,u\left( s,\xi  \right)ds=$$
$$=-\int\limits_{0}^{x}{d\xi }\int\limits_{0}^{t}{u\left( s,\xi  \right)ds}\int\limits_{s}^{t}{{{\left( \eta -s \right)}^{-{{\beta }_{1}}}}E_{\alpha ,\,1-{{\beta }_{1}}}^{-{{\gamma }_{1}}}\left[ \delta {{\left( \eta -s \right)}^{\alpha }} \right]{{V}_{\eta }}\left( t,x;\eta ,\xi  \right)d\eta }\,=$$
$$=-\int\limits_{0}^{x}{\int\limits_{0}^{t}{u\left( \eta ,\xi  \right)I_{\eta t}^{\alpha ,1-{{\beta }_{1}},-{{\gamma }_{1}},\delta }}{{V}_{\eta }}\left( t,x;\eta ,\xi  \right)d\eta d\xi }=\int\limits_{0}^{x}{\int\limits_{0}^{t}{u\left( \eta ,\xi  \right){}^{PC}D_{\eta t}^{\alpha ,{{\beta }_{1}},{{\gamma }_{1}},\delta }}V\left( t,x;\eta ,\xi  \right)d\eta d\xi }.$$
Let us evaluate ${{I}_{2}}$ that is the second integral of \eqref{eq3.3}. We integrate by parts with respect to $\xi$ and consider the conditions \eqref{eq2.11} and \eqref{eq3.2}, then we get the following result:
$${{I}_{2}}=\int\limits_{0}^{t}{\int\limits_{0}^{x}{{{u}_{\xi }}\left( \eta ,\xi  \right)V\left( t,x;\eta ,\xi  \right)d\xi }d\eta }=-\int\limits_{0}^{t}{V\left( t,x;\eta ,0 \right){{\varphi }_{0}}\left( \eta  \right)d\eta -}\int\limits_{0}^{t}{\int\limits_{0}^{x}{u\left( \eta ,\xi  \right){{V}_{\xi }}\left( t,x;\eta ,\xi  \right)d\xi d\eta }}.$$
Substituting the results of the integrals ${{I}_{1}}$ and ${{I}_{2}}$ into the equation \eqref{eq3.3}, we obtain the following expression: 
$$\int\limits_{0}^{t}{\int\limits_{0}^{x}{u\left( \eta ,\xi  \right)\left[ {}^{PC}D_{\eta t}^{\alpha ,{{\beta }_{1}},{{\gamma }_{1}},\delta }V\left( t,x;\eta ,\xi  \right)-{{V}_{\xi }}\left( t,x;\eta ,\xi  \right) \right]d\xi d\eta }}=$$
$$=\int\limits_{0}^{t}{{{\varphi }_{0}}\left( \eta  \right)V\left( t,x;\eta ,0 \right)d\eta }+\int\limits_{0}^{t}{\int\limits_{0}^{x}{v\left( \eta ,\xi  \right)V\left( t,x;\eta ,\xi  \right)d\xi d\eta }}.$$
According to \eqref{eq3.1}, we get the following result:
$$\int\limits_{0}^{t}{\int\limits_{0}^{x}{u\left( \eta ,\xi  \right)d\xi d\eta }}=\int\limits_{0}^{t}{{{\varphi }_{0}}\left( \eta  \right)V\left( t,x;\eta ,0 \right)d\eta }+\int\limits_{0}^{t}{\int\limits_{0}^{x}{v\left( \eta ,\xi  \right)V\left( t,x;\eta ,\xi  \right)d\xi d\eta }}.$$
To find the unknown function $u\left( t,x \right)$ from the final expression, we first differentiate with respect to $x$, then with respect to $t$, and arrive at the following result:
\begin{equation}\label{eq3.4}
u\left( t,x \right)=\int\limits_{0}^{t}{{{\varphi }_{0}}\left( \eta  \right){{V}_{xt}}\left( t,x;\eta ,0 \right)d\eta }+\int\limits_{0}^{t}{\int\limits_{0}^{x}{v\left( \eta ,\xi  \right){{V}_{xt}}\left( t,x;\eta ,\xi  \right)d\xi d\eta }}.
\end{equation}
Our goal is to find the function $V\left( t,x;\eta ,\xi  \right).$ Next, we will assume that the function $V\left( t,x;\eta ,\xi  \right)$ can be represented as a function of the difference of the arguments $t,\,\,\,\eta$ and $x,\,\,\,\xi $:
$$V\left( t,x;\eta ,\xi  \right)=\widetilde{V}\left( t-\eta ,x-\xi  \right)=\widetilde{V}\left( y,s \right).$$
From \eqref{eq3.2} we get
$${{\left. V\left( t,x;\eta ,\xi  \right) \right|}_{\eta =t}}=\widetilde{V}\left( 0,s \right)=0,\,\,\,  {{\left. V\left( t,x;\eta ,\xi  \right) \right|}_{\xi =x}}=\widetilde{V}\left( y,0 \right)=0.$$
Now, let’s also express equation \eqref{eq3.1} in terms of the function $\widetilde{V}\left( y,s \right)$. In this step, we utilize the relation between the Prabhakar derivative of Riemann-Liouville type and the Prabhakar derivative of Caputo type. 

The following formula holds between the Prabhakar derivative of Riemann-Liouville type and the Prabhakar derivative of Caputo type for any function $f\in A{{C}^{m}}\left( a,b \right)$ \cite{Fernandez}:
$${}^{PC}D_{at}^{\alpha ,\beta ,\gamma ,\delta }f\left( t \right)={}^{PRL}D_{at}^{\alpha ,\beta ,\gamma ,\delta }\left[ f\left( t \right)-\sum\limits_{j=0}^{m-1}{\frac{{{f}^{\left( j \right)}}\left( a \right)}{j!}{{\left( t-a \right)}^{j}}} \right],$$
where $\alpha ,\beta ,\delta ,\gamma \in \mathbb{C}$ with $\operatorname{Re}\left( \alpha  \right)>0,$ $\operatorname{Re}\left( \beta  \right)\ge 0$ and $m=\left[ \operatorname{Re}\left( \beta  \right) \right]+1.$

In our case, $m=\left[ \operatorname{Re}\left( {{\beta }_{1}} \right) \right]+1=1,$ ${{\left. V\left( t,x;\eta ,\xi  \right) \right|}_{t=\eta }}=0.$ Hence, the following equality is true:
$${}^{PC}D_{\eta t}^{\alpha ,{{\beta }_{1}},{{\gamma }_{1}},\delta }V\left( t,x;\eta ,\xi  \right)={}^{PRL}D_{\eta t}^{\alpha ,{{\beta }_{1}},{{\gamma }_{1}},\delta }V\left( t,x;\eta ,\xi  \right).$$

Now, we consider the effect of the operator to the function $V\left( t,x;\eta ,\xi  \right)$:
$${}^{PC}D_{\eta t}^{\alpha ,{{\beta }_{1}},{{\gamma }_{1}},\delta }V\left( t,x;\eta ,\xi  \right)={}^{PRL}D_{\eta t}^{\alpha ,{{\beta }_{1}},{{\gamma }_{1}},\delta }V\left( t,x;\eta ,\xi  \right)=-\frac{\partial }{\partial \eta }{}^{P}I_{\eta t}^{\alpha ,\,1-{{\beta }_{1}},\,-{{\gamma }_{1}},\,\delta }V\left( t,x;\eta ,\xi  \right)=$$
$$=-\frac{\partial }{\partial \eta }\int\limits_{\eta }^{t}{{{\left( z-\eta  \right)}^{-{{\beta }_{1}}}}}E_{\alpha ,{{\beta }_{1}}}^{{{\gamma }_{1}}}\left[ \delta {{\left( z-\eta  \right)}^{\alpha }} \right]V\left( t,x;z,\xi  \right)dz=\left\{ z=t-h \right\}=$$
$$=\frac{\partial }{\partial \eta }\int\limits_{t-\eta }^{0}{{{\left( t-\eta -h \right)}^{-{{\beta }_{1}}}}}E_{\alpha ,{{\beta }_{1}}}^{{{\gamma }_{1}}}\left[ \delta {{\left( t-\eta -h \right)}^{\alpha }} \right]V\left( t,x;t-h,\xi  \right)dh=$$
$$=\frac{\partial }{\partial \eta }\int\limits_{t-\eta }^{0}{{{\left( t-\eta -h \right)}^{-{{\beta }_{1}}}}}E_{\alpha ,{{\beta }_{1}}}^{{{\gamma }_{1}}}\left[ \delta {{\left( t-\eta -h \right)}^{\alpha }} \right]\widetilde{V}\left( h,s \right)dh=\left\{ t-\eta =y \right\}=$$
$$=-\frac{\partial }{\partial y}\int\limits_{y}^{0}{{{\left( y-h \right)}^{-{{\beta }_{1}}}}}E_{\alpha ,{{\beta }_{1}}}^{{{\gamma }_{1}}}\left[ \delta {{\left( y-h \right)}^{\alpha }} \right]\widetilde{V}\left( h,s \right)dh=$$
$$=\frac{\partial }{\partial y}\int\limits_{0}^{y}{{{\left( y-h \right)}^{-{{\beta }_{1}}}}}E_{\alpha ,{{\beta }_{1}}}^{{{\gamma }_{1}}}\left[ \delta {{\left( y-h \right)}^{\alpha }} \right]\widetilde{V}\left( h,s \right)dh={}^{PRL}D_{0y}^{\alpha ,{{\beta }_{1}},{{\gamma }_{1}},\delta }\widetilde{V}\left( y,s \right)={}^{PC}D_{0y}^{\alpha ,{{\beta }_{1}},{{\gamma }_{1}},\delta }\widetilde{V}\left( y,s \right).$$
Now, we calculate the first-order derivative:
$${{V}_{\xi }}\left( t,x;\eta ,\xi  \right)={{\widetilde{V}}_{s}}\left( y,s \right)\cdot {{s}_{\xi }}=-{{\widetilde{V}}_{s}}\left( y,s \right).$$
In that case, it follows from \eqref{eq3.1} and \eqref{eq3.2} that $\widetilde{V}\left( y,s \right)$ is the solution of the following problem:
\begin{equation}\label{eq3.5}
    {}^{PC}D_{0y}^{\alpha ,{{\beta }_{1}},{{\gamma }_{1}},\delta }\widetilde{V}\left( y,s \right)+{{\widetilde{V}}_{s}}\left( y,s \right)=1,
\end{equation}
\begin{equation}\label{eq3.6}
\widetilde{V}\left( 0,s \right)=0, \,\,\,\, \widetilde{V}\left( y,0 \right)=0.
\end{equation}
We apply the Laplace transform to the both sides of the equation \eqref{eq3.5} with respect to $y:$
$${{L}_{y}}\left[ {}^{PC}D_{0y}^{\alpha ,{{\beta }_{1}},{{\gamma }_{1}},\delta }\widetilde{V}\left( y,s \right)+{{\widetilde{V}}_{s}}\left( y,s \right) \right]={{L}_{y}}\left[ {}^{PC}D_{0y}^{\alpha ,{{\beta }_{1}},{{\gamma }_{1}},\delta }\widetilde{V}\left( y,s \right) \right]+{{L}_{y}}\left[ {{\widetilde{V}}_{s}}\left( y,s \right) \right]={{L}_{y}}\left[ 1 \right].$$
If we denote ${{L}_{y}}\left[ \widetilde{V}\left( y,s \right) \right]=\mu \left( s;p \right),$ then we get
\begin{equation}\label{eq3.7}
    {{L}_{y}}\left[ {}^{PC}D_{0y}^{\alpha ,{{\beta }_{1}},{{\gamma }_{1}},\delta }\widetilde{V}\left( y,s \right) \right]+{{\mu }_{s}}\left( s;p \right)=\frac{1}{p}.	
\end{equation}
Now, we apply the Laplace transform to the operator:
$${{L}_{y}}\left[ {}^{PC}D_{0y}^{\alpha ,{{\beta }_{1}},{{\gamma }_{1}},\delta }\widetilde{V}\left( y,s \right) \right]=\int\limits_{0}^{+\infty }{{{e}^{-py}}}\int\limits_{0}^{y}{{{\left( y-z \right)}^{-{{\beta }_{1}}}}}E_{\alpha ,1-{{\beta }_{1}}}^{-{{\gamma }_{1}}}\left[ \delta {{\left( y-z \right)}^{\alpha }} \right]{{\widetilde{V}}_{z}}\left( z,s \right)dzdy=$$
$$=\int\limits_{0}^{+\infty }{{{\widetilde{V}}_{z}}\left( z,s \right)dz}\int\limits_{z}^{+\infty }{{{e}^{-py}}{{\left( y-z \right)}^{-{{\beta }_{1}}}}}E_{\alpha ,1-{{\beta }_{1}}}^{-{{\gamma }_{1}}}\left[ \delta {{\left( y-z \right)}^{\alpha }} \right]dy=\left\{ y=l+z \right\}=$$
$$=\int\limits_{0}^{+\infty }{{{e}^{-pz}}{{\widetilde{V}}_{z}}\left( z,s \right)dz}\int\limits_{0}^{+\infty }{{{e}^{-pl}}{{l}^{-{{\beta }_{1}}}}}E_{\alpha ,1-{{\beta }_{1}}}^{-{{\gamma }_{1}}}\left[ \delta {{l}^{\alpha }} \right]dl.$$
First, we evaluate the integral with respect to $l$ and we use this formula 
\begin{equation}\label{eq3.8}
    \sum\limits_{k=0}^{+\infty }{\frac{{{\left( a \right)}_{k}}{{x}^{k}}}{k!}}={{\left( 1-x \right)}^{-a}}
\end{equation}
to obtain the result:
$$\int\limits_{0}^{+\infty }{{{e}^{-pl}}{{l}^{-{{\beta }_{1}}}}}E_{\alpha ,-{{\beta }_{1}}}^{-{{\gamma }_{1}}}\left[ \delta {{l}^{\alpha }} \right]dl=\sum\limits_{k=0}^{+\infty }{\frac{{{\left( -{{\gamma }_{1}} \right)}_{k}}{{\delta }^{k}}}{\Gamma \left( \alpha k+1-{{\beta }_{1}} \right)k!}}\int\limits_{0}^{+\infty }{{{e}^{-pl}}{{l}^{\alpha k-{{\beta }_{1}}}}}dl=$$
$$=\left\{ pl=v \right\}={{p}^{{{\beta }_{1}}-1}}\sum\limits_{k=0}^{+\infty }{\frac{{{\left( -\gamma  \right)}_{k}}{{\delta }^{k}}{{p}^{-\alpha k}}}{\Gamma \left( \alpha k+1-{{\beta }_{1}} \right)k!}}\int\limits_{0}^{+\infty }{{{e}^{-v}}{{v}^{\alpha k-{{\beta }_{1}}}}}dv=$$
$$={{p}^{{{\beta }_{1}}-1}}\sum\limits_{k=0}^{+\infty }{\frac{{{\left( -{{\gamma }_{1}} \right)}_{k}}{{\left( \delta {{p}^{-\alpha }} \right)}^{k}}}{k!}}={{p}^{{{\beta }_{1}}-1}}{{\left( 1-\frac{\delta }{{{p}^{\alpha }}} \right)}^{{{\gamma }_{1}}}}.$$
Now, we integrate by parts with respect to $z:$
$$\int\limits_{0}^{+\infty }{{{e}^{-pz}}{{\widetilde{V}}_{z}}\left( z,s \right)dz}=p\int\limits_{0}^{+\infty }{{{e}^{-pz}}\widetilde{V}\left( z,s \right)dz}=p{{L}_{y}}\left[ \widetilde{V}\left( y,s \right) \right].$$
As a result, we have identified 
$${{L}_{y}}\left[ {}^{PC}D_{0t}^{\alpha ,{{\beta }_{1}},{{\gamma }_{1}},\delta }\widetilde{V}\left( y,s \right) \right]={{p}^{{{\beta }_{1}}}}{{\left( 1-\frac{\delta }{{{p}^{\alpha }}} \right)}^{{{\gamma }_{1}}}}{{L}_{y}}\left[ \widetilde{V}\left( y,s \right) \right],$$ so we can rewrite \eqref{eq3.7} and get the linear differential equation:
$${{\mu }_{s}}\left( s;p \right)+\lambda \left( p \right)\,\mu \left( s;p \right)=\frac{1}{p},$$
where $\lambda \left( p \right)={{p}^{{{\beta }_{1}}}}{{\left( 1-\frac{\delta }{{{p}^{\alpha }}} \right)}^{{{\gamma }_{1}}}}.$

We solve this linear differential equation and determine that $\mu \left( s;p \right)=C{{e}^{-\lambda s}}+\frac{1}{p\lambda }.$  Since $\widetilde{V}\left( y,0 \right)=0,$ it follows that $\mu \left( 0;p \right)=0\,.$ By using $\mu \left( 0;p \right)=0\,,$ we can find that $C=-\frac{1}{p\lambda }.$ As a result, we obtain the following
$${{L}_{y}}\left[ \widetilde{V}\left( y,s \right) \right]=\mu \left( s;p \right)=-\frac{{{e}^{-\lambda s}}}{p\lambda }+\frac{1}{p\lambda }.$$
Now, we seek the inverse Laplace transform:
\begin{equation}\label{eq3.9}
    L_{y}^{-1}\left[ \mu  \right]=L_{y}^{-1}\left[ -\frac{{{e}^{-\lambda s}}}{p\lambda }+\frac{1}{p\lambda } \right]=-L_{y}^{-1}\left[ \frac{{{e}^{-\lambda s}}}{p\lambda } \right]+L_{y}^{-1}\left[ \frac{1}{p\lambda } \right].
\end{equation}
We apply the Laplace transform to the function 
${{\widetilde{V}}_{1}}\left( y,s \right)={{y}^{{{\beta }_{1}}}}\sum\limits_{n=0}^{+\infty }{\frac{{{\left( -s \right)}^{n}}{{y}^{-{{\beta }_{1}}n}}}{n!}E_{\alpha ,-{{\beta }_{1}}n+1+{{\beta }_{1}}}^{-{{\gamma }_{1}}n+{{\gamma }_{1}}}\left[ \delta {{y}^{\alpha }} \right]}:$
$${{L}_{y}}\left[ {{\widetilde{V}}_{1}}\left( y,s \right) \right]={{L}_{y}}\left[ {{y}^{{{\beta }_{1}}}}\sum\limits_{n=0}^{+\infty }{\frac{{{\left( -s \right)}^{n}}{{y}^{-{{\beta }_{1}}n}}}{n!}E_{\alpha ,-{{\beta }_{1}}n+1+{{\beta }_{1}}}^{-{{\gamma }_{1}}n+{{\gamma }_{1}}}\left[ \delta {{y}^{\alpha }} \right]} \right]=$$
$$=\int\limits_{0}^{+\infty }{{{e}^{-py}}}{{y}^{{{\beta }_{1}}}}\sum\limits_{n=0}^{+\infty }{\frac{{{\left( -s \right)}^{n}}{{y}^{-{{\beta }_{1}}n}}}{n!}E_{\alpha ,-{{\beta }_{1}}n+1+{{\beta }_{1}}}^{-{{\gamma }_{1}}n+{{\gamma }_{1}}}\left[ \delta {{y}^{\alpha }} \right]}\,dy=$$
$$=\sum\limits_{n=0}^{+\infty }{\frac{{{\left( -s \right)}^{n}}}{n!}\sum\limits_{k=0}^{+\infty }{\frac{{{\left( -{{\gamma }_{1}}n+{{\gamma }_{1}} \right)}_{k}}{{\delta }^{k}}}{k!\Gamma \left( \alpha k-{{\beta }_{1}}n+1+{{\beta }_{1}} \right)}}}\int\limits_{0}^{+\infty }{{{e}^{-py}}}{{y}^{{{\beta }_{1}}-{{\beta }_{1}}n+\alpha k}}\,dy=$$
$$=\left\{ y=\frac{\theta }{p} \right\}=\sum\limits_{n=0}^{+\infty }{\frac{{{\left( -s \right)}^{n}}}{n!}\sum\limits_{k=0}^{+\infty }{\frac{{{\left( -{{\gamma }_{1}}n+{{\gamma }_{1}} \right)}_{k}}{{\delta }^{k}}{{p}^{-\alpha k+{{\beta }_{1}}n-{{\beta }_{1}}-1}}}{k!}}}=$$
$$={{p}^{-{{\beta }_{1}}-1}}\sum\limits_{n=0}^{+\infty }{\frac{{{\left( -s \right)}^{n}}{{p}^{{{\beta }_{1}}n}}}{n!}\sum\limits_{k=0}^{+\infty }{\frac{{{\left( -{{\gamma }_{1}}n+{{\gamma }_{1}} \right)}_{k}}{{\left[ \delta {{p}^{-\alpha }} \right]}^{k}}}{k!}}}.$$
Using the formula \eqref{eq3.8}, we get
$${{L}_{y}}\left[ {{\widetilde{V}}_{1}}\left( y,s \right) \right]={{p}^{-{{\beta }_{1}}-1}}\sum\limits_{n=0}^{+\infty }{\frac{{{\left( -s \right)}^{n}}{{p}^{{{\beta }_{1}}n}}}{n!}{{\left( 1-\delta {{p}^{-\alpha }} \right)}^{{{\gamma }_{1}}n-{{\gamma }_{1}}}}}=$$
$$=\frac{1}{p}\sum\limits_{n=0}^{+\infty }{\frac{{{\left( -s \right)}^{n}}{{\left[ {{p}^{{{\beta }_{1}}}}{{\left( 1-\frac{\delta }{{{p}^{\alpha }}} \right)}^{{{\gamma }_{1}}}} \right]}^{n}}}{n!\left[ {{p}^{{{\beta }_{1}}}}{{\left( 1-\frac{\delta }{{{p}^{\alpha }}} \right)}^{{{\gamma }_{1}}}} \right]}}=\frac{1}{p\lambda }\sum\limits_{n=0}^{+\infty }{\frac{{{\left( -\lambda s \right)}^{n}}}{n!}}=\frac{{{e}^{-\lambda s}}}{p\lambda }.$$
For $s=0$, we have
$${{L}_{y}}\left[ {{\widetilde{V}}_{1}}\left( y,0 \right) \right]={{L}_{y}}\left[ {{y}^{\beta }}E_{\alpha ,1+{{\beta }_{1}}}^{{{\gamma }_{1}}}\left[ \delta {{y}^{\alpha }} \right] \right]=\frac{1}{p\lambda }.$$
According to \eqref{eq3.9} and the results ${{L}_{y}}\left[ {{\widetilde{V}}_{1}}\left( y,s \right) \right]$ and ${{L}_{y}}\left[ {{\widetilde{V}}_{1}}\left( y,0 \right) \right]$, we obtain
$$\widetilde{V}\left( y,s \right)=L_{y}^{-1}\left[ \mu  \right]=-{{y}^{{{\beta }_{1}}}}\sum\limits_{n=0}^{+\infty }{\frac{{{\left( -s \right)}^{n}}{{y}^{-{{\beta }_{1}}n}}}{n!}E_{\alpha ,-{{\beta }_{1}}n+1+{{\beta }_{1}}}^{-{{\gamma }_{1}}n+{{\gamma }_{1}}}\left[ \delta {{y}^{\alpha }} \right]}+{{y}^{{{\beta }_{1}}}}E_{\alpha ,1+{{\beta }_{1}}}^{{{\gamma }_{1}}}\left[ \delta {{y}^{\alpha }} \right]=$$
$$=-{{y}^{{{\beta }_{1}}}}\sum\limits_{n=1}^{+\infty }{\frac{{{\left( -s \right)}^{n}}{{y}^{-{{\beta }_{1}}n}}}{n!}E_{\alpha ,-{{\beta }_{1}}n+1+{{\beta }_{1}}}^{-{{\gamma }_{1}}n+{{\gamma }_{1}}}\left[ \delta {{y}^{\alpha }} \right]}=\sum\limits_{n=0}^{+\infty }{\frac{{{\left( -1 \right)}^{n}}{{s}^{n+1}}{{y}^{-{{\beta }_{1}}n}}}{\left( n+1 \right)!}E_{\alpha ,-{{\beta }_{1}}n+1}^{-{{\gamma }_{1}}n}\left[ \delta {{y}^{\alpha }} \right]}.$$
Now, let us return to the function $V\left( t,x;\eta ,\xi  \right):$
$$V\left( t,x;\eta ,\xi  \right)=\widetilde{V}\left( t-\eta ,x-\xi  \right)=\sum\limits_{n=0}^{+\infty }{\frac{{{\left( -1 \right)}^{n}}{{\left( x-\xi  \right)}^{n+1}}{{\left( t-\eta  \right)}^{-{{\beta }_{1}}n}}}{\left( n+1 \right)!}E_{\alpha ,-{{\beta }_{1}}n+1}^{-{{\gamma }_{1}}n}\left[ \delta {{\left( t-\eta  \right)}^{\alpha }} \right]}.$$
To find the solution \eqref{eq3.4}, we take the derivative of the function $V\left( t,x;\eta ,\xi  \right)$ with respect to $x$:
$${{V}_{x}}\left( t,x;\eta ,\xi  \right)=\sum\limits_{n=0}^{+\infty }{\frac{{{\left( -1 \right)}^{n}}{{\left( x-\xi  \right)}^{n}}}{n!}{{\left( t-\eta  \right)}^{-{{\beta }_{1}}n}}E_{\alpha ,-{{\beta }_{1}}n+1}^{-{{\gamma }_{1}}n}\left[ \delta {{\left( t-\eta  \right)}^{\alpha }} \right]}=$$
$$=1+\sum\limits_{n=1}^{+\infty }{\frac{{{\left( -1 \right)}^{n}}{{\left( x-\xi  \right)}^{n}}}{n!}{{\left( t-\eta  \right)}^{-{{\beta }_{1}}n}}E_{\alpha ,-{{\beta }_{1}}n+1}^{-{{\gamma }_{1}}n}\left[ \delta {{\left( t-\eta  \right)}^{\alpha }} \right]},$$
then  with respect to $t$:
$${{V}_{xt}}\left( t,x;\eta ,\xi  \right)=\sum\limits_{n=1}^{+\infty }{\frac{{{\left( -1 \right)}^{n}}{{\left( x-\xi  \right)}^{n}}}{n!}{{\left( t-\eta  \right)}^{-{{\beta }_{1}}n-1}}E_{\alpha ,-{{\beta }_{1}}n}^{-{{\gamma }_{1}}n}\left[ \delta {{\left( t-\eta  \right)}^{\alpha }} \right]}=$$
$$=\sum\limits_{n=0}^{+\infty }{\frac{{{\left( -1 \right)}^{n}}{{\left( x-\xi  \right)}^{n}}}{n!}{{\left( t-\eta  \right)}^{-{{\beta }_{1}}n-1}}E_{\alpha ,-{{\beta }_{1}}n}^{-{{\gamma }_{1}}n}\left[ \delta {{\left( t-\eta  \right)}^{\alpha }} \right]}.$$
Finally, using \eqref{eq3.4}, we obtain the solution \eqref{eq2.12}.

\subsection{\bf Fifth appendix}\label{secA5}

We evaluate the limit of the function $\omega(t,x)$ at $t=0$:

$$\lim_{t\to 0}\omega \left( t,x \right)=\lim_{t\to 0}\sum\limits_{n=0}^{+\infty }{\frac{{{\left( -1 \right)}^{n}}{{x}^{n}}}{n!}{{t}^{-{{\beta }_{1}}n-1}}E_{\alpha ,-{{\beta }_{1}}n}^{-{{\gamma }_{1}}n}\left[ \delta {{t}^{\alpha }} \right]}=$$
$$=\lim_{t\to 0}\sum\limits_{n=0}^{+\infty }{\frac{{{\left( -1 \right)}^{n}}{{x}^{n}}}{n!}{{t}^{-{{\beta }_{1}}n-1}}}\lim_{t\to 0}\sum_{m=0}^{+\infty}\frac{(-\gamma_1n)_m \delta^mt^{\alpha m}}{m!\Gamma(\alpha m-\beta_1n)}=$$
$$=\lim_{t\to 0}\sum\limits_{n=0}^{+\infty }{\frac{{{\left( -1 \right)}^{n}}{{x}^{n}}}{n!}{{t}^{-{{\beta }_{1}}n-1}}}\lim_{t\to 0}\left[\frac{1}{\Gamma(-\beta_1n)}+\frac{(-\gamma_1n)_1\delta t^\alpha}{\Gamma(\alpha-\beta_1n)}+...+\frac{(-\gamma_1n)_m \delta^mt^{\alpha m}}{m!\Gamma(\alpha m-\beta_1n)}+...\right]=$$
$$=\lim_{t\to 0}t^{-1}\sum\limits_{n=0}^{+\infty }\frac{1}{\Gamma(n+1)\Gamma(-\beta_1n)}\left(-\frac{x}{t^{\beta_1}}\right)^n=\lim_{t\to 0}t^{-1}e_{1,\beta_1}^{1,0}\left(-\frac{x}{t^{\beta_1}}\right).$$

According to the Lemma 3.1, we also obtain the following:

$$\lim_{t\to 0}t^{-1}e_{1,\beta_1}^{1,0}\left(-\frac{x}{t^{\beta_1}}\right)\le \lim_{t\to 0}\left|t^{-1}e_{1,\beta_1}^{1,0}\left(-\frac{x}{t^{\beta_1}}\right)\right|\le$$
$$\le\frac{1}{\beta_1 \pi }{{\left( \frac{\cos \beta_1 \theta \pi }{\beta_1 \pi } \right)}^{-\frac{1}{\beta_1 }}}\Gamma \left( \frac{1}{\beta_1 } \right){{\left( x-{{x}_{0}} \right)}^{-\frac{1}{\beta_1 }}}\lim_{t\to 0}e_{1,\beta_1 }^{1,1}\left( -\frac{{{x}_{0}}}{{{t}^{\beta_1 }}} \right).$$

Since $\lim_{t\to 0}e_{1,{{\beta }_{1}}}^{1,1}\left( -\frac{{{x}_{0}}}{{{ t }^{{{\beta }_{1}}}}} \right)=0,$ it follows that 
$$\lim_{t\to 0}\omega \left( t,x \right)=0.$$

Thus, the function $\omega(t,x)$ does not have a singularity at $t=0$.

\subsection{\bf Sixth appendix}\label{secA6}

In the following, we examine the problem  \eqref{eq2.14}-\eqref{eq2.15}-\eqref{eq2.16}. This problem is solved in the same way as the problem \eqref{eq2.9}-\eqref{eq2.10}-\eqref{eq2.11} above. Similarly, we consider a function $F\left( t,x;\eta ,\xi  \right),$ $\left\{ \left( \eta ,\xi  \right):\,\,0<\eta <t,\,\,x<\xi <a \right\},$ which satisfies the following equation for any fixed $\left( t,x \right):$ 
\begin{equation}\label{eq3.10}
    {}^{PC}D_{\eta t}^{\alpha ,{{\beta }_{1}},{{\gamma }_{1}},\delta }F\left( t,x;\eta ,\xi  \right)+{{F}_{\xi }}\left( t,x;\eta ,\xi  \right)=1,
\end{equation}
also satisfies these conditions
\begin{equation}\label{eq3.11}
    {{\left. F\left( t,x;\eta ,\xi  \right) \right|}_{\eta =t}}=0,   \,\,\,\,{{\left. F\left( t,x;\eta ,\xi  \right) \right|}_{\xi =x}}=0.
\end{equation}
First, we rewrite the equation \eqref{eq2.14} by replacing the variables $x$ and $t$ with $\xi$ and $\eta .$ Then, we multiply both sides of the last equation by $F\left( t,x;\eta ,\xi  \right),$ integrate with respect to $\xi $ over the interval $\left[ x,a \right],$ and with respect to $\eta $ over the interval $\left[ 0,t \right]:$
$$\int\limits_{0}^{t}{\int\limits_{x}^{a}{F\left( t,x;\eta ,\xi  \right){}^{PRL}D_{0\eta }^{\alpha ,\,{{\beta }_{1}},\,{{\gamma }_{1}},\,\delta }v\left( \eta ,\xi  \right)d\xi d\eta }}-\int\limits_{0}^{t}{\int\limits_{x}^{a}{{{v}_{\xi }}\left( \eta ,\xi  \right)F\left( t,x;\eta ,\xi  \right)d\xi }d\eta }=$$
\begin{equation}\label{eq3.12}
    =\int\limits_{0}^{t}{\int\limits_{x}^{a}{f\left( \eta ,\xi  \right)F\left( t,x;\eta ,\xi  \right)d\xi d\eta }}.
\end{equation}
Let us denote the first integral in \eqref{eq3.12} by ${{L}_{1}}$ and integrate by parts with respect to $\eta $ considering the conditions \eqref{eq2.15} and \eqref{eq3.11}:
$${{L}_{1}}=\int\limits_{0}^{t}{\int\limits_{x}^{a}{F\left( t,x;\eta ,\xi  \right){}^{PRL}D_{0\eta }^{\alpha ,\,{{\beta }_{1}},\,{{\gamma }_{1}},\,\delta }v\left( \eta ,\xi  \right)d\xi d\eta }}=$$
$$=\int\limits_{0}^{t}{\int\limits_{x}^{a}{F\left( t,x;\eta ,\xi  \right)\frac{\partial }{\partial \eta }I_{0\eta }^{\alpha ,\,1-{{\beta }_{1}},\,-{{\gamma }_{1}},\,\delta }v\left( \eta ,\xi  \right)d\xi d\eta }}=$$
$$=-\int\limits_{x}^{a}{F\left( t,x;0,\xi  \right)\tau \left( \xi  \right)d\xi -}\int\limits_{x}^{a}{d\xi }\int\limits_{0}^{t}{{{F}_{\eta }}\left( t,x;\eta ,\xi  \right)I_{0\eta }^{\alpha ,\,1-{{\beta }_{1}},\,-{{\gamma }_{1}},\,\delta }v\left( \eta ,\xi  \right)d\eta }.$$
By repeating the calculations performed above for ${{I}_{1}}$, we obtain the following result:
$${{L}_{1}}=-\int\limits_{x}^{a}{F\left( t,x;0,\xi  \right)\tau \left( \xi  \right)d\xi +\int\limits_{x}^{a}{\int\limits_{0}^{t}{v\left( \eta ,\xi  \right)}}}\,{}^{PC}D_{\eta t}^{\alpha ,{{\beta }_{1}},{{\gamma }_{1}},\delta }F\left( t,x;\eta ,\xi  \right)d\eta d\xi .$$
The second integral in the equation \eqref{eq3.12} is computed in the same manner as ${{I}_{2}}$:
$${{L}_{2}}=\int\limits_{0}^{t}{F\left( t,x;\eta ,a \right)\psi \left( \eta  \right)d\eta -}\int\limits_{0}^{t}{\int\limits_{x}^{a}{v\left( \eta ,\xi  \right){{F}_{\xi }}\left( t,x;\eta ,\xi  \right)d\xi d\eta }}.$$
According to the results of the integrals ${{L}_{1}}$ and ${{L}_{2}}$ we rewrite \eqref{eq3.12}: 
$$\int\limits_{0}^{t}{\int\limits_{x}^{a}{v\left( \eta ,\xi  \right)\left[ {}^{PC}D_{\eta t}^{\alpha ,{{\beta }_{1}},{{\gamma }_{1}},\delta }F\left( t,x;\eta ,\xi  \right)+{{F}_{\xi }}\left( t,x;\eta ,\xi  \right) \right]d\xi d\eta }}=$$
$$=\int\limits_{0}^{t}{\psi \left( \eta  \right)F\left( t,x;\eta ,a \right)d\eta }+\,\int\limits_{x}^{a}{\tau \left( \xi  \right)F\left( t,x;0,\xi  \right)d\xi }+\int\limits_{0}^{t}{\int\limits_{x}^{a}{f\left( \eta ,\xi  \right)F\left( t,x;\eta ,\xi  \right)d\xi d\eta }}.$$
According to \eqref{eq3.10}, we get the following result:
$$\int\limits_{0}^{t}{\int\limits_{x}^{a}{v\left( \eta ,\xi  \right)d\xi d\eta }}=\int\limits_{0}^{t}{\psi \left( \eta  \right)F\left( t,x;\eta ,a \right)d\eta }+\,$$
$$+\,\int\limits_{x}^{a}{\tau \left( \xi  \right)F\left( t,x;0,\xi  \right)d\xi }+\int\limits_{0}^{t}{\int\limits_{x}^{a}{f\left( \eta ,\xi  \right)F\left( t,x;\eta ,\xi  \right)d\xi d\eta }}.$$
Now, we differentiate with respect to $x$, then with respect to $t$, and arrive at the following result:
$$v\left( t,x \right)=-\int\limits_{0}^{t}{\psi \left( \eta  \right){{F}_{xt}}\left( t,x;\eta ,a \right)d\eta }-\,$$
\begin{equation}\label{eq3.13}
-\int\limits_{x}^{a}{\tau \left( \xi  \right){{F}_{xt}}\left( t,x;0,\xi  \right)d\xi }-\int\limits_{0}^{t}{\int\limits_{x}^{a}{f\left( \eta ,\xi  \right){{F}_{xt}}\left( t,x;\eta ,\xi  \right)d\xi d\eta }}.
\end{equation}
Our goal is to find the function $F\left( t,x;\eta ,\xi  \right).$ Next, we will assume that the function $F\left( t,x;\eta ,\xi  \right)$ can be represented as a function of the difference of the arguments $t,\,\,\,\eta $ and $\xi ,\,\,\,x$:
$$F\left( t,x;\eta ,\xi  \right)=\widetilde{F}\left( t-\eta ,\xi -x \right)=\widetilde{F}\left( y,s \right).$$
In that case, it follows from \eqref{eq3.10} and \eqref{eq3.11} that $\widetilde{F}\left( y,s \right)$ is the solution of the following problem:
\begin{equation}\label{eq3.14}
{}^{PC}D_{0y}^{\alpha ,{{\beta }_{1}},{{\gamma }_{1}},\delta }\widetilde{F}\left( y,s \right)+{{\widetilde{F}}_{s}}\left( y,s \right)=1,
\end{equation}
\begin{equation}\label{eq3.15}
\widetilde{F}\left( 0,s \right)=0, \,\,\,\, \widetilde{F}\left( y,0 \right)=0.
\end{equation}
The problem \eqref{eq3.14}-\eqref{eq3.15} is exactly the same as the problem \eqref{eq3.5}-\eqref{eq3.6}, therefore, the solution can easily be written as follows:
$$\widetilde{F}\left( y,s \right)=\sum\limits_{n=0}^{+\infty }{\frac{{{\left( -1 \right)}^{n}}{{s}^{n+1}}{{y}^{-{{\beta }_{1}}n}}}{\left( n+1 \right)!}E_{\alpha ,-{{\beta }_{1}}n+1}^{-{{\gamma }_{1}}n}\left[ \delta {{y}^{\alpha }} \right]}.$$
The function $F\left( t,x;\eta ,\xi  \right)$ then takes the following form:
$$F\left( t,x;\eta ,\xi  \right)=\widetilde{F}\left( t-\eta ,\xi -x \right)=\sum\limits_{n=0}^{+\infty }{\frac{{{\left( -1 \right)}^{n}}{{\left( \xi -x \right)}^{n+1}}{{\left( t-\eta  \right)}^{-{{\beta }_{1}}n}}}{\left( n+1 \right)!}E_{\alpha ,-{{\beta }_{1}}n+1}^{-{{\gamma }_{1}}n}\left[ \delta {{\left( t-\eta  \right)}^{\alpha }} \right]}.$$
To find the solution \eqref{eq3.13}, we take the derivatives of the function $F\left( t,x;\eta ,\xi  \right)$ with respect to $x$ and $t$:
$${{F}_{x}}\left( t,x;\eta ,\xi  \right)=-\sum\limits_{n=0}^{+\infty }{\frac{{{\left( -1 \right)}^{n}}{{\left( \xi -x \right)}^{n}}}{n!}{{\left( t-\eta  \right)}^{-{{\beta }_{1}}n}}E_{\alpha ,-{{\beta }_{1}}n+1}^{-{{\gamma }_{1}}n}\left[ \delta {{\left( t-\eta  \right)}^{\alpha }} \right]}=$$
$$=-1-\sum\limits_{n=1}^{+\infty }{\frac{{{\left( -1 \right)}^{n}}{{\left( \xi -x \right)}^{n}}}{n!}{{\left( t-\eta  \right)}^{-{{\beta }_{1}}n}}E_{\alpha ,-{{\beta }_{1}}n+1}^{-{{\gamma }_{1}}n}\left[ \delta {{\left( t-\eta  \right)}^{\alpha }} \right]};$$
$${{F}_{xt}}\left( t,x;\eta ,\xi  \right)=-\sum\limits_{n=1}^{+\infty }{\frac{{{\left( -1 \right)}^{n}}{{\left( \xi -x \right)}^{n}}}{n!}{{\left( t-\eta  \right)}^{-{{\beta }_{1}}n-1}}E_{\alpha ,-{{\beta }_{1}}n}^{-{{\gamma }_{1}}n}\left[ \delta {{\left( t-\eta  \right)}^{\alpha }} \right]}=$$
$$=-\sum\limits_{n=0}^{+\infty }{\frac{{{\left( -1 \right)}^{n}}{{\left( \xi -x \right)}^{n}}}{n!}{{\left( t-\eta  \right)}^{-{{\beta }_{1}}n-1}}E_{\alpha ,-{{\beta }_{1}}n}^{-{{\gamma }_{1}}n}\left[ \delta {{\left( t-\eta  \right)}^{\alpha }} \right]}.$$
We rewrite \eqref{eq3.13} and obtain the solution \eqref{eq2.17}.

\subsection{\bf Seventh appendix} \label{secA7}

Now we prove the equality \eqref{eq2.19}. To do this, using the form \eqref{eq2.13} of $\omega \left( t,x \right)$, we rewrite the left-hand side of \eqref{eq2.19}:
$$\int\limits_{0}^{y}{\omega \left( y-t,{{x}_{1}} \right)}\,\omega \left( t,{{x}_{2}} \right)dt=\int\limits_{0}^{y}{\sum\limits_{n=0}^{+\infty }{\frac{{{\left( -1 \right)}^{n}}{{x}_{1}}^{n}}{n!}{{\left( y-t \right)}^{-{{\beta }_{1}}n-1}}E_{\alpha ,-{{\beta }_{1}}n}^{-{{\gamma }_{1}}n}\left[ \delta {{\left( y-t \right)}^{\alpha }} \right]}}\times $$
$$\times \sum\limits_{m=0}^{+\infty }{\frac{{{\left( -1 \right)}^{m}}{{x}_{2}}^{m}}{m!}{{t}^{-{{\beta }_{1}}m-1}}E_{\alpha ,-{{\beta }_{1}}m}^{-{{\gamma }_{1}}m}\left[ \delta {{t}^{\alpha }} \right]}\,dt.$$
We use formula \eqref{eq2.27} two consecutive times:
$$\sum\limits_{n=0}^{+\infty }{\sum\limits_{m=0}^{n}{\frac{{{\left( -1 \right)}^{n}}x_{1}^{m}x_{2}^{n-m}}{m!\left( n-m \right)!}\int\limits_{0}^{y}{{{\left( y-t \right)}^{-{{\beta }_{1}}m-1}}{{t}^{-{{\beta }_{1}}n+{{\beta }_{1}}m-1}}}\times }}$$
$$\times \sum\limits_{k=0}^{+\infty }{\frac{{{\left( -{{\gamma }_{1}}m \right)}_{k}}{{\delta }^{k}}{{\left( y-t \right)}^{\alpha k}}}{\Gamma \left( \alpha k-{{\beta }_{1}}m \right)k!}}\sum\limits_{z=0}^{+\infty }{\frac{{{\left( -{{\gamma }_{1}}n+{{\gamma }_{1}}m \right)}_{z}}{{\delta }^{z}}{{t}^{\alpha z}}}{\Gamma \left( \alpha k-{{\beta }_{1}}n+{{\beta }_{1}}m \right)z!}}\,dt=$$
$$=\sum\limits_{n=0}^{+\infty }{\sum\limits_{m=0}^{n}{\frac{{{\left( -1 \right)}^{n}}x_{1}^{m}x_{2}^{n-m}}{m!\left( n-m \right)!}}}\sum\limits_{k=0}^{+\infty }{\sum\limits_{z=0}^{k}{\frac{{{\left( -{{\gamma }_{1}}m \right)}_{z}}{{\left( -{{\gamma }_{1}}n+{{\gamma }_{1}}m \right)}_{k-z}}{{\delta }^{k}}}{\Gamma \left( \alpha z-{{\beta }_{1}}m \right)\Gamma \left( \alpha k-\alpha z-{{\beta }_{1}}n+{{\beta }_{1}}m \right)z!\left( k-z \right)!}}}\times $$
$$\times \int\limits_{0}^{y}{{{\left( y-t \right)}^{\alpha z-{{\beta }_{1}}m-1}}{{t}^{\alpha k-\alpha z-{{\beta }_{1}}n+{{\beta }_{1}}m-1}}}dt.$$
Utilizing the formulas \eqref{eq2.28} and \eqref{eq2.29}, we get
$$\int\limits_{0}^{y}{\omega \left( y-t,{{x}_{1}} \right)}\,\omega \left( t,{{x}_{2}} \right)dt=\sum\limits_{n=0}^{+\infty }{\frac{{{\left( -1 \right)}^{n}}{{\left( {{x}_{1}}+{{x}_{2}} \right)}^{n}}}{n!}\sum\limits_{k=0}^{+\infty }{\frac{{{\left( -{{\gamma }_{1}}n \right)}_{k}}{{\delta }^{k}}{{y}^{\alpha k-{{\beta }_{1}}n-1}}}{k!\Gamma \left( \alpha k-{{\beta }_{1}}n \right)}}}=$$
$$=\sum\limits_{n=0}^{+\infty }{\frac{{{\left( -1 \right)}^{n}}{{\left( {{x}_{1}}+{{x}_{2}} \right)}^{n}}}{n!}{{y}^{-{{\beta }_{1}}n-1}}E_{\alpha ,-{{\beta }_{1}}n}^{-{{\gamma }_{1}}n}\left[ \delta {{y}^{\alpha }} \right]}=\omega \left( y,{{x}_{1}}+{{x}_{2}} \right).$$
The formula \eqref{eq2.19} is proven.

\subsection{\bf Eighth appendix} \label{secA8}

Let us prove 
$$\int\limits_{a-x}^{a+x}{\omega \left( t-y,s \right)ds}=I_{yt}^{\alpha ,\,{{\beta }_{1}},\,{{\gamma }_{1}},\,\delta }\left[ \omega \left( t-y,a-x \right)-\omega \left( t-y,a+x \right) \right].$$
Firstly, we see the left-hand side:
$$\int\limits_{a-x}^{a+x}{\omega \left( t-y,s \right)ds}=\int\limits_{a-x}^{a+x}{\sum\limits_{n=0}^{+\infty }{\frac{{{\left( -1 \right)}^{n}}{{s}^{n}}}{n!}{{\left( t-y \right)}^{-{{\beta }_{1}}n-1}}E_{\alpha ,-{{\beta }_{1}}n}^{-{{\gamma }_{1}}n}\left[ \delta {{\left( t-y \right)}^{\alpha }} \right]}\,ds}=$$
$$=\left. \sum\limits_{n=0}^{+\infty }{\frac{{{\left( -1 \right)}^{n}}{{s}^{n+1}}}{\left( n+1 \right)!}{{\left( t-y \right)}^{-{{\beta }_{1}}n-1}}E_{\alpha ,-{{\beta }_{1}}n}^{-{{\gamma }_{1}}n}\left[ \delta {{\left( t-y \right)}^{\alpha }} \right]} \right|_{s=a-x}^{s=a+x}=$$
$$=\left. \sum\limits_{n=-1}^{+\infty }{\frac{{{\left( -1 \right)}^{n}}{{s}^{n+1}}}{\left( n+1 \right)!}{{\left( t-y \right)}^{-{{\beta }_{1}}n-1}}E_{\alpha ,-{{\beta }_{1}}n}^{-{{\gamma }_{1}}n}\left[ \delta {{\left( t-y \right)}^{\alpha }} \right]} \right|_{s=a-x}^{s=a+x}+\left. {{\left( t-y \right)}^{{{\beta }_{1}}-1}}E_{\alpha ,{{\beta }_{1}}}^{{{\gamma }_{1}}}\left[ \delta {{\left( t-y \right)}^{\alpha }} \right] \right|_{s=a-x}^{s=a+x}=$$
$$=\left. -\sum\limits_{n=0}^{+\infty }{\frac{{{\left( -1 \right)}^{n}}{{s}^{n}}}{n!}{{\left( t-y \right)}^{-{{\beta }_{1}}n+{{\beta }_{1}}-1}}E_{\alpha ,-{{\beta }_{1}}n+{{\beta }_{1}}}^{-{{\gamma }_{1}}n+{{\gamma }_{1}}}\left[ \delta {{\left( t-y \right)}^{\alpha }} \right]} \right|_{s=a-x}^{s=a+x}=$$
$$={{\left( t-y \right)}^{{{\beta }_{1}}-1}}\sum\limits_{n=0}^{+\infty }{\frac{{{\left( -1 \right)}^{n}}\left[ {{\left( a-x \right)}^{n}}-{{\left( a+x \right)}^{n}} \right]}{n!}{{\left( t-y \right)}^{-{{\beta }_{1}}n}}E_{\alpha ,-{{\beta }_{1}}n+{{\beta }_{1}}}^{-{{\gamma }_{1}}n+{{\gamma }_{1}}}\left[ \delta {{\left( t-y \right)}^{\alpha }} \right]}.$$
Now, we consider the right-hand side:
$$I_{yt}^{\alpha ,\,{{\beta }_{1}},\,{{\gamma }_{1}},\,\delta }\omega \left( t-y,a-x \right)=\int\limits_{y}^{t}{{{\left( t-s \right)}^{{{\beta }_{1}}-1}}}E_{\alpha ,\,{{\beta }_{1}}}^{{{\gamma }_{1}}}\left[ \delta {{\left( t-s \right)}^{\alpha }} \right]\omega \left( s-y,a-x \right)ds=$$
$$=\sum\limits_{n=0}^{+\infty }{\frac{{{\left( x-a \right)}^{n}}}{n!}}\int\limits_{y}^{t}{{{\left( t-s \right)}^{{{\beta }_{1}}-1}}}{{\left( s-y \right)}^{-{{\beta }_{1}}n-1}}E_{\alpha ,\,{{\beta }_{1}}}^{{{\gamma }_{1}}}\left[ \delta {{\left( t-s \right)}^{\alpha }} \right]E_{\alpha ,-{{\beta }_{1}}n}^{-{{\gamma }_{1}}n}\left[ \delta {{\left( s-y \right)}^{\alpha }} \right]ds.$$
We begin by applying formula \eqref{eq2.27}, then substitute $\left( t-y \right)z=s-y$, and lastly make use of \eqref{eq2.28}:
$$I_{yt}^{\alpha ,\,{{\beta }_{1}},\,{{\gamma }_{1}},\,\delta }\omega \left( t-y,a-x \right)=\sum\limits_{\begin{smallmatrix} 
 n=0 \\ 
 k=0 
\end{smallmatrix}}^{+\infty }{\sum\limits_{m=0}^{k}{\frac{{{\delta }^{k}}{{\left( {{\gamma }_{1}} \right)}_{m}}{{\left( -{{\gamma }_{1}}n \right)}_{k-m}}{{\left( x-a \right)}^{n}}}{m!\left( k-m \right)!n!\Gamma \left( \alpha m+{{\beta }_{1}} \right)\Gamma \left( \alpha k-\alpha m-{{\beta }_{1}}n \right)}}}\times $$
$$\times \,\int\limits_{y}^{t}{{{\left( t-s \right)}^{\alpha m+{{\beta }_{1}}-1}}}{{\left( s-y \right)}^{\alpha k-\alpha m-{{\beta }_{1}}n-1}}ds=$$
$$=\sum\limits_{\begin{smallmatrix} 
 n=0 \\ 
 k=0 
\end{smallmatrix}}^{+\infty }{\sum\limits_{m=0}^{k}{\frac{{{\delta }^{k}}{{\left( {{\gamma }_{1}} \right)}_{m}}{{\left( -{{\gamma }_{1}}n \right)}_{k-m}}{{\left( x-a \right)}^{n}}{{\left( t-y \right)}^{\alpha k-{{\beta }_{1}}n+{{\beta }_{1}}-1}}}{m!\left( k-m \right)!n!\Gamma \left( \alpha k-{{\beta }_{1}}n+{{\beta }_{1}} \right)}}}=$$
$$=\sum\limits_{n=0}^{+\infty }{\sum\limits_{k=0}^{+\infty }{\frac{{{\left( {{\gamma }_{1}}-{{\gamma }_{1}}n \right)}_{k}}{{\left( x-a \right)}^{n}}{{\delta }^{k}}{{\left( t-y \right)}^{\alpha k-{{\beta }_{1}}n+{{\beta }_{1}}-1}}}{k!n!\Gamma \left( \alpha k-{{\beta }_{1}}n+{{\beta }_{1}} \right)}}}=$$
$$={{\left( t-y \right)}^{{{\beta }_{1}}-1}}\sum\limits_{n=0}^{+\infty }{\frac{{{\left( -1 \right)}^{n}}{{\left( a-x \right)}^{n}}}{n!}{{\left( t-y \right)}^{-{{\beta }_{1}}n}}E_{\alpha ,-{{\beta }_{1}}n+{{\beta }_{1}}}^{-{{\gamma }_{1}}n+{{\gamma }_{1}}}\left[ \delta {{\left( t-y \right)}^{\alpha }} \right]}.$$
$I_{yt}^{\alpha ,\,{{\beta }_{1}},\,{{\gamma }_{1}},\,\delta }\omega \left( t-y,a+x \right)$ is calculated in the same way as above, and we arrive at the following result:
$$I_{yt}^{\alpha ,\,{{\beta }_{1}},\,{{\gamma }_{1}},\,\delta }\left[ \omega \left( t-y,a-x \right)-\omega \left( t-y,a+x \right) \right]=$$
$$={{\left( t-y \right)}^{{{\beta }_{1}}-1}}\sum\limits_{n=0}^{+\infty }{\frac{{{\left( -1 \right)}^{n}}\left[ {{\left( a-x \right)}^{n}}-{{\left( a+x \right)}^{n}} \right]}{n!}{{\left( t-y \right)}^{-{{\beta }_{1}}n}}E_{\alpha ,-{{\beta }_{1}}n+{{\beta }_{1}}}^{-{{\gamma }_{1}}n+{{\gamma }_{1}}}\left[ \delta {{\left( t-y \right)}^{\alpha }} \right]}.$$
As a result, the correctness of the expression under consideration has been proven.

\subsection{\bf Ninth appendix} \label{secA9}

Let us show that given  $KF\left( t \right)=\int\limits_{0}^{t}{F\left( y \right)\omega \left( t-y,2a \right)}\,dy$ it follows that 
$${{K}^{n}}F\left( t \right)=\int\limits_{0}^{t}{F\left( y \right)\omega \left( t-y,2na \right)}\,dy.$$
Firstly, let us calculate ${{K}^{2}}F\left( t \right)$:
$${{K}^{2}}F\left( t \right)=\int\limits_{0}^{t}{\left[ \int\limits_{0}^{y}{F\left( s \right)\omega \left( y-s,2a \right)}\,ds \right]\omega \left( t-y,2a \right)}\,dy=$$
$$=\int\limits_{0}^{t}{\omega \left( t-y,2a \right)}\,dy\int\limits_{0}^{y}{F\left( s \right)\omega \left( y-s,2a \right)}\,ds=\int\limits_{0}^{t}{F\left( s \right)}\,ds\int\limits_{s}^{t}{\omega \left( t-y,2a \right)\omega \left( y-s,2a \right)}\,dy.$$
Next, after making the substitution $z=y-s$ and referring to equality \eqref{eq2.19}, we reach the following result:
$${{K}^{2}}F\left( t \right)=\int\limits_{0}^{t}{F\left( s \right)}\,ds\int\limits_{0}^{t-s}{\omega \left( t-s-z,2a \right)\omega \left( z,2a \right)}\,dz=\int\limits_{0}^{t}{F\left( s \right)\omega \left( t-s,4a \right)}\,ds.$$
Now, let’s determine ${{K}^{3}}F\left( t \right)$:
$${{K}^{3}}F\left( t \right)=\int\limits_{0}^{t}{\left[ \int\limits_{0}^{y}{F\left( s \right)\omega \left( y-s,2a \right)}\,ds \right]\omega \left( t-y,4a \right)}\,dy=\int\limits_{0}^{t}{F\left( s \right)}\,ds\int\limits_{s}^{t}{\omega \left( t-y,4a \right)\omega \left( y-s,2a \right)}\,dy=$$
$$=\int\limits_{0}^{t}{F\left( s \right)}\,ds\int\limits_{0}^{t-s}{\omega \left( t-s-z,4a \right)\omega \left( z,2a \right)}\,dz=\int\limits_{0}^{t}{F\left( s \right)\omega \left( t-s,6a \right)}\,ds.$$
In conclusion
$$KF\left( t \right)=\int\limits_{0}^{t}{F\left( y \right)\omega \left( t-y,2a \right)}\,dy,$$
$${{K}^{2}}F\left( t \right)=\int\limits_{0}^{t}{F\left( y \right)\omega \left( t-y,4a \right)}\,dy,$$
$${{K}^{3}}F\left( t \right)=\int\limits_{0}^{t}{F\left( y \right)\omega \left( t-y,6a \right)}\,dy,$$
$$…$$
$${{K}^{n}}F\left( t \right)=\int\limits_{0}^{t}{F\left( y \right)\omega \left( t-y,2na \right)}\,dy.$$

\subsection{\bf Tenth appendix} \label{secA10}

We prove that 
$$\underset{y\to t}{\mathop{\lim }}\,\left[ W\left( t-y,a-x,a+x \right)\times  \right.$$
$$\left. \times \int\limits_{0}^{y}{{{\varphi }_{0}}\left( \eta  \right)\sum\limits_{n=0}^{+\infty }{\frac{{{\left( -1 \right)}^{n}}{{a}^{n}}}{n!}{{\left( y-\eta  \right)}^{-{{\beta }_{1}}-{{\beta }_{1}}n}}E_{\alpha ,1-{{\beta }_{1}}-{{\beta }_{1}}n}^{-{{\gamma }_{1}}-{{\gamma }_{1}}n}\left[ \delta {{\left( y-\eta  \right)}^{\alpha }} \right]}\,d\eta } \right]=0.$$
Since 
$$\underset{y\to t}{\mathop{\lim }}\,\int\limits_{0}^{y}{{{\varphi }_{0}}\left( \eta  \right)\sum\limits_{n=0}^{+\infty }{\frac{{{\left( -1 \right)}^{n}}{{a}^{n}}}{n!}{{\left( y-\eta  \right)}^{-{{\beta }_{1}}-{{\beta }_{1}}n}}E_{\alpha ,1-{{\beta }_{1}}-{{\beta }_{1}}n}^{-{{\gamma }_{1}}-{{\gamma }_{1}}n}\left[ \delta {{\left( y-\eta  \right)}^{\alpha }} \right]}\,d\eta }\ne 0,$$
$$\underset{y\to t}{\mathop{\lim }}\,W\left( t-y,a-x,a+x \right)=0.$$ 
Therefore, we compute the following limit:
$$\underset{y\to t}{\mathop{\lim }}\,W\left( t-y,a-x,a+x \right)=$$
$$=\frac{1}{2}\underset{y\to t}{\mathop{\lim }}\,\sum\limits_{n=0}^{+\infty }{\frac{{{\left( -1 \right)}^{n}}\left[ {{\left( a-x \right)}^{n}}-{{\left( a+x \right)}^{n}} \right]}{n!}{{\left( t-y \right)}^{{{\beta }_{1}}-{{\beta }_{1}}n-1}}E_{\alpha ,{{\beta }_{1}}-{{\beta }_{1}}n}^{{{\gamma }_{1}}-{{\gamma }_{1}}n}\left[ \delta {{\left( t-y \right)}^{\alpha }} \right]=}$$
$$=\frac{1}{2}\underset{y\to t}{\mathop{\lim }}\,\sum\limits_{n=0}^{+\infty }{\frac{{{\left( -1 \right)}^{n}}\left[ {{\left( a-x \right)}^{n}}-{{\left( a+x \right)}^{n}} \right]}{n!}{{\left( t-y \right)}^{{{\beta }_{1}}-{{\beta }_{1}}n-1}}\times }$$
$$\times \underset{y\to t}{\mathop{\lim }}\,\left[ \frac{1}{\Gamma \left( {{\beta }_{1}}-{{\beta }_{1}}n \right)}+\frac{{{\left( {{\gamma }_{1}}-{{\gamma }_{1}}n \right)}_{1}}\delta {{\left( t-y \right)}^{\alpha }}}{\Gamma \left( \alpha +{{\beta }_{1}}-{{\beta }_{1}}n \right)}+\frac{{{\left( {{\gamma }_{1}}-{{\gamma }_{1}}n \right)}_{2}}{{\delta }^{2}}{{\left( t-y \right)}^{2\alpha }}}{2\Gamma \left( 2\alpha +{{\beta }_{1}}-{{\beta }_{1}}n \right)}+... \right]=$$
$$=\frac{1}{2}\underset{y\to t}{\mathop{\lim }}\,\left\{ {{\left( t-y \right)}^{{{\beta }_{1}}-1}}\sum\limits_{n=0}^{+\infty }{\frac{{{\left( -1 \right)}^{n}}\left[ {{\left( a-x \right)}^{n}}-{{\left( a+x \right)}^{n}} \right]}{\Gamma \left( n+1 \right)\Gamma \left( {{\beta }_{1}}-{{\beta }_{1}}n \right)}{{\left( t-y \right)}^{-{{\beta }_{1}}n}}} \right\}=$$
$$=\frac{1}{2}\underset{y\to t}{\mathop{\lim }}\,\left\{ {{\left( t-y \right)}^{{{\beta }_{1}}-1}}\left[ e_{1,{{\beta }_{1}}}^{1,{{\beta }_{1}}}\left( -\frac{a-x}{{{\left( t-y \right)}^{{{\beta }_{1}}}}} \right)-e_{1,{{\beta }_{1}}}^{1,{{\beta }_{1}}}\left( -\frac{a+x}{{{\left( t-y \right)}^{{{\beta }_{1}}}}} \right) \right] \right\}.$$

According to the Lemma 3.1, we also obtain the following:
$$\frac{1}{2}\underset{y\to t}{\mathop{\lim }}\,{{\left( t-y \right)}^{{{\beta }_{1}}-1}}e_{1,{{\beta }_{1}}}^{1,{{\beta }_{1}}}\left( -\frac{a-x}{{{\left( t-y \right)}^{{{\beta }_{1}}}}} \right)\le \frac{1}{2}\underset{y\to t}{\mathop{\lim }}\,\left| {{\left( t-y \right)}^{{{\beta }_{1}}-1}}e_{1,{{\beta }_{1}}}^{1,{{\beta }_{1}}}\left( -\frac{a-x}{{{\left( t-y \right)}^{{{\beta }_{1}}}}} \right) \right|\le $$
$$\le \frac{1}{2{{\beta }_{1}}\pi }{{\left( \frac{\cos {{\beta }_{1}}\theta \pi }{{{\beta }_{1}}\pi } \right)}^{\frac{{{\beta }_{1}}-1}{{{\beta }_{1}}}}}\Gamma \left( \frac{1-{{\beta }_{1}}}{{{\beta }_{1}}} \right){{\left( x-{{x}_{0}} \right)}^{\frac{{{\beta }_{1}}-1}{{{\beta }_{1}}}}}\underset{y\to t}{\mathop{\lim }}\,e_{1,{{\beta }_{1}}}^{1,1}\left( -\frac{{{x}_{0}}}{{{\left( t-y \right)}^{{{\beta }_{1}}}}} \right).$$
Since $\underset{y\to t}{\mathop{\lim }}\,e_{1,{{\beta }_{1}}}^{1,1}\left( -\frac{{{x}_{0}}}{{{\left( t-y \right)}^{{{\beta }_{1}}}}} \right)=0,$ it follows that $$\underset{y\to t}{\mathop{\lim }}\,W\left( t-y,a-x,a+x \right)=0.$$

\bibliographystyle{plain}

\begin{thebibliography}{99}
\normalsize

\bibitem{Al-Refai}
Al-Refai, M., Nusseir, A., Al-Sharif, S.:
Maximum principles for fractional differential inequalities with Prabhakar derivative and their applications.
{Fractal Fract.} \textbf{6}(10), 612 (2022).
DOI: 10.3390/fractalfract6100612

\bibitem{El-Sayed}
El-Sayed, A.A., Boulaaras, S., Al-Kharousi, F.A.:
Dickson polynomial-based solutions for fractional order physics problems.
{J. Inequal. Appl.} (2025), 118:118.
DOI: 10.1186/s13660-025-03367-7

\bibitem{Fernandez}
Fernandez, A., Restrepo, J.E., Suragan, D.:
Prabhakar-type linear differential equations with variable coefficients.
{Differ. Integral Equ.} \textbf{35}(9--10), 581--610 (2022).
DOI: 10.57262/die035-0910-581

\bibitem{Garra}
Garra, R., Gorenflo, R., Polito, F., Tomovski, Z.:
Hilfer–Prabhakar derivatives and some applications.
{Appl. Math. Comput.} \textbf{242}, 576--589 (2014).
DOI: 10.1016/j.amc.2014.05.129

\bibitem{Giusti}
Giusti, A., Colombaro, I., Garra, R., et al.:
A practical guide to Prabhakar fractional calculus.
{Fract. Calc. Appl. Anal.} \textbf{23}(1), 9--54 (2020).
DOI: 10.1515/fca-2020-0034

\bibitem{guisti-col} Guisti, A., Colombaro, I.: Prabhakar-like fractional viscoelasticity.
{ Communications in Nonlinear Science and Numerical Simulation}. {\bf 56}, 138--143 (2018).
 
\bibitem{Graham}
Graham, R.L., Knuth, D.E., Patashnik, O.:
Concrete mathematics: A foundation for computer science.
Addison-Wesley, Reading, MA (1994).

\bibitem{Karimov}
Karimov, E., Hasanov, A.:
On a boundary-value problem in a bounded domain for a time-fractional diffusion equation with the Prabhakar fractional derivative.
{Bull. Karaganda Univ. Math. Ser.} \textbf{111}(3), 39--46 (2023).
DOI: 10.31489/2023m3/39-46

\bibitem{Karimov2} 
Karimov~E., Usmonov~D., Mirzaeva~M.;
Green’s Function and Solution Representation for a Boundary Value Problem Involving the Prabhakar Fractional Derivative.
Preprint, arXiv:2512.21259 [math.AP], 2025.

\bibitem{Kilbas}
Kilbas, A.A., Srivastava, H.M., Trujillo, J.J.:
Theory and Applications of Fractional Differential Equations.
Elsevier, Amsterdam (2006).

\bibitem{Knopp}
Knopp, K.:
Theory and application of infinite series.
Dover Publications, New York (1990).

\bibitem{Kolmogorov}
Kolmogorov, A.N., Fomin, S.V.:
Elements of the Theory of Functions and Functional Analysis.
Nauka, Moscow (1972).

\bibitem{Magar}
Magar, S.K., Dole, P.V., Ghadle, K.P.:
Prabhakar and Hilfer–Prabhakar fractional derivatives in the setting of $\Psi$-fractional calculus and its applications.
{Kragujevac J. Math.} \textbf{48}(4), 515--533 (2024).
DOI: 10.46793/KgJMat2404.515M

\bibitem{Odibat}
Odibat, Z., Momani, S.M.:
Fractional Green's function for fractional partial differential equations.
{J. Eur. Syst. Autom.} \textbf{42}(6--8), 639--651 (2008).
DOI: 10.3166/jesa.42.639-651

\bibitem{Podlubny}
Podlubny, I.:
Fractional Differential Equations.
Academic Press, San Diego (1999).

\bibitem{Prabhakar}
Prabhakar, T.R.:
A singular integral equation with a generalized Mittag-Leffler function in the kernel.
{Yokohama Math. J.} \textbf{19}, 7--15 (1971).

\bibitem{Prudnikov}
Prudnikov, A.P., Brychkov, Yu.A., Marichev, O.I.:
Integrals and Series: Special Functions.
Nauka, Moscow (1983).

\bibitem{Pskhu}
Pskhu, A.V.:
Green functions of the first boundary-value problem for a fractional diffusion–wave equation in multidimensional domains.
{Mathematics} \textbf{8}(4), 464 (2020).
DOI: 10.3390/math8040464

\bibitem{Pskhu 2}
Pskhu, A.V.:
Uravneniya v chastnykh proizvodnykh drobnogo poryadka.
Nauka, Moscow (2005).

\bibitem{Sarah}
Aljohani, S., Rashid, M., Kalsoom, A., Mlaiki, N.:
Exploring the influence of generalized kernels on Green’s function in fractional differential equations.
{Int. J. Anal. Appl.} \textbf{22}, 188 (2024).
DOI: 10.28924/2291-8639-22-2024-188

\bibitem{Suzuki}
Suzuki, J., Zayernouri, M., D’Elia, M.:
A survey of fractional-order models in transport and anomalous materials (Technical Report).
DOE Tech. Rep. SAND2021-11291R (2021).
DOI: 10.2172/1820001

\bibitem{T22} Tarasov ,V.E.: Fractional dynamics with depreciation and obsolescence: equations with Prabhakar fractional derivatives. {Mathematics}. {\bf 10}, 1540 (2022). 

\bibitem{Turdiev 1}
Turdiev, Kh.:
Nonlocal problem for a diffusion–wave equation involving regularized Prabhakar fractional derivative.
{Bull. Inst. Math.} \textbf{6}(6), 39--45 (2023).

\bibitem{Turdiev 2}
Turdiev, Kh.N., Usmonov, D.A.:
The Goursat problem for generalized (fractional) hyperbolic-type equation.
{Uzbek Math. J.} \textbf{69}(2), 300--306 (2025).
DOI: 10.29229/uzmj.2025-2-29

\bibitem{Usmonov}
Usmonov, D., Mirzaeva, M.:
A Cauchy problem for the sub-diffusion equation with the Prabhakar fractional derivative.
{Gulf J. Math.} \textbf{21}(2), 181--203 (2025).
DOI: 10.56947/gjom.v21i2.3708


\bibitem{wang} Wang, W., Metzler R., Tomovski Z. Distributed-order fractional diffusion equation with Hilfer-Prabhakar fractional derivative,
{Physica A: Statistical Mechanics and its Applications}. {\bf 689},  131370 (2026).
 
\end{thebibliography}

\end{document}